\documentclass[poms,final,nonblindrev]{poms_cls_no_journal} 

\OneAndAHalfSpacedXI 

\usepackage{bbm}
\usepackage{lscape}
\usepackage[inline]{enumitem}
\usepackage{amssymb}
\usepackage{multicol}
\usepackage{tikz}
\usepackage{pgfplots}
\usepackage{graphicx}
\usepackage[font=footnotesize,labelfont=bf]{caption}
\usepackage{subcaption}
\usepackage{setspace}
\usepackage{float}
\usepackage{natbib}
\usepackage{xurl}
\usepackage{bm}
\usepackage{afterpage}

\newcommand*\bV{\boldsymbol{V}}
\newcommand*\bp{\boldsymbol{p}}
\newcommand*\bs[1]{\boldsymbol{s}_{#1}}
\newcommand*\bx[1]{\boldsymbol{x}_{[#1]}}
\newcommand*\bVp{\boldsymbol{V,p}}
\newcommand*\cS{\mathcal{S}}
\newcommand*\cV{\mathcal{V}}
\newcommand*\imp{\mathcal{P}}
\newcommand*\content{\mathcal{C}}
\newcommand*\success{\rho}
\newcommand*\bsuccess{\boldsymbol{\rho}}
\newcommand*\profit{\pi}
\newcommand*\effect{\Phi}

\newcommand*\dplus[2]{\delta_{#1^+} (#2)}
\newcommand*\dminus[2]{\delta_{#1^-} (#2)}
\newcommand*\dave[2]{\delta_{#1^\mathrm{ave}} (#2)}

 \bibpunct[, ]{(}{)}{,}{a}{}{,}%
 \def\bibfont{\small}%
\TheoremsNumberedThrough     
\ECRepeatTheorems

\EquationsNumberedThrough    

\makeatletter
\let\normalsize\relax 
\let\@currsize\normalsize
\makeatother

\begin{document}


 \RUNAUTHOR{Aslan, Bakir and \c{C}avdar}

\RUNTITLE{Observational Learning in Reward-based Crowdfunding}

\TITLE{How does Observational Learning Impact Crowdfunding Outcomes for Backers, Project Creators and Platforms?}

\ARTICLEAUTHORS{%
\AUTHOR{Ayse Aslan}
\AFF{School of Computing Engineering and the Built Environment, Edinburgh Napier University, UK \\ \EMAIL{a.aslan@napier.ac.uk}} 
\AUTHOR{Ilke Bakir}
\AFF{Department of Operations, Faculty of Economics and Business, University of Groningen, The Netherlands \\ \EMAIL{i.bakir@rug.edu}}
\AUTHOR{Bahar \c{C}avdar}
\AFF{Department of Industrial and Systems Engineering, Rensselaer Polytechnic Institute, Troy, NY, USA
 \\ \EMAIL{cavdab2@rpi.edu}}
} 

\ABSTRACT{%
Reward-based crowdfunding platforms are becoming increasingly popular to finance projects proposing innovative products, e.g., Kickstarter. One important challenge of this form of financing is the uncertainty in the quality of projects. To mitigate the negative effects of this uncertainty for backers, platforms share information regarding the decisions of earlier backers visiting the project campaign pages. This allows backers not only to rely on their expertise to identify project qualities but also to learn from the decisions of their fellow backers who might be more informed. Current studies on observational learning (OL) in crowdfunding mainly focus on predicting the success chances of projects, and there is a lack of understanding of how OL affects crowdfunding dynamics for backers, project creators and platforms. This paper aims to fill this gap by using a theoretical OL model involving two projects competing for funding from backers who may have differentiated expertness in identifying project quality. By introducing various performance measures for backers, creators and platforms and comparing these measures under OL to the case without learning, we provide a thorough analysis of how OL impacts crowdfunding outcomes. We find that information sharing and OL always benefit backers, especially when the early backers are experts.  Regarding the impact of OL on creators and platforms, our analysis reveals two understudied but important aspects: the tightness of the competition for projects according to the availability of funding, and the quality difference among the proposed projects. Additionally, we investigate how OL affects the quality decisions of creators and show that OL increases the incentive for high-quality products, especially in situations where funding is scarce.
}%

\KEYWORDS{reward-based crowdfunding; observational learning; quality management; Bayesian rationality}

\maketitle

%


\section{Introduction}
\textit{Crowdfunding} is a financing mechanism used by entrepreneurs to raise funding for their projects from many individual investors (e.g., Kickstarter.com and Indiegogo.com). Crowdfunding industry has grown tremendously over the years, surpassing alternative financing forms such as venture capital \citep{Barnett2015}, and its growth is projected to continue in the future. The global market size for crowdfunding activities is estimated to reach \$1.3 billion by 2028 while growing at a compound annual growth rate (CAGR) of 4.4\% from 2021 to 2028 \citep{Manuel2022}. Many kinds of crowdfunding business models \citep[see][]{Belleflamme2015} are operated through online platforms. In a reward-based platform such as Kickstarter.com, investors are offered rewards in the form of products for their contributions to projects, while equity- and lending-based platforms may provide monetary benefits to their investors. In this paper, we focus on \textit{reward-based crowdfunding}, where investors not only provide financial support but also act as the consumers of the products. 

In reward-based crowdfunding platforms, entrepreneurs, whom we refer to as \textit{creators} from now on, launch their crowdfunding campaigns by providing information on the products they are proposing; while investors, whom we refer to as \textit{backers}, seek to make pledges to campaigns of their choice, or not. To secure funding, creators offer a menu of product rewards and corresponding pledge amounts for backers to choose from if the project successfully meets its funding goal. Various campaign structures exist, but we specifically focus on the \textit{all-or-nothing} (AoN) campaign structure, which has become the most prevalent campaign type, particularly on popular platforms such as Kickstarter.com. In an AoN campaign, creators specify a funding target and a deadline. To successfully receive the collected funds from backers, the project must reach this specified funding target by the set deadline. If the campaign falls short of the target by the deadline, all collected funds are refunded to the backers.

Since the introduction of reward-based crowdfunding platforms, many campaigns have been successful. As of the time of writing, Kickstarter.com, which is one of the prominent global platforms for reward-based crowdfunding, has funded over 26\,000 projects. Nevertheless, despite its success stories, reward-based crowdfunding  faces several challenges that require attention. One of the major obstacles in crowdfunding is the uncertainty surrounding the quality of campaigned products. 
Backers visiting a campaign's webpage read the product description, review the menu of rewards, and then must decide whether to pledge or not without being certain of the true quality and state of the final product. Because of this uncertainty, backers may inadvertently pledge support to low-quality projects, only to regret their investment after receiving products that fail to meet expectations or experiencing lengthy delivery delays. For example, according to the fulfillment report published by Kickstarter.com\footnote{https://www.kickstarter.com/fulfillment}, only 65\% of backers confirmed the statement ``the reward was delivered on time''.

The uncertainty in the quality of proposed projects can negatively impact the \textit{effectiveness} of crowdfunding platforms. In effective crowdfunding platforms, the limited funding potential of backers is matched with high-quality and promising projects; therefore, projects that otherwise could not secure the financing they need to launch their products have better chances of success than low-quality ones. One well-known success story of effective crowdfunding where a high-quality product is put on the market thanks to crowdfunding is the Pebble smartwatch. This campaign raised over \$10 million, and backers were more than content with the received reward products as their quality exceeded expectations, which has made Pebble one of the early market leaders in the smartwatch industry \citep{Metz2016}. However, crowdfunding is not always effective in matching backers with promising projects. Sometimes unworthy projects can also attract considerable funding because backers may misjudge project quality. One good example for this is the Kreyos smartwatch campaign, which launched on Indiegogo.com nearly one year after Pebble's campaign on Kickstarter.com with a similar funding target and price as Pebble, and it was able to raise over \$1.5 million. The final product was significantly delayed, and backers were disappointed because many promised features were lacking and devices were defective \citep{Newman2014}. Such examples could cause backers to regret their participation in crowdfunding, which can harm the long-term reputation of crowdfunding and discourage many from participating, ultimately reducing the market potential of crowdfunding platforms. Lack of knowledge on products that leads to the funding of low-quality projects can also hinder the success chances of creators with promising projects. This is because platforms are highly dynamic ventures where many similar projects are being campaigned at the same time that demand funding from backers. For example, if both Pebble (April 2012) and Kreyos (June 2013), which proposed very similar products with similar funding targets, had launched their campaigns at the same time, a considerable portion of the funding potential of backers could have ended up with Kreyos, instead of Pebble, which could have potentially caused Pebble to fail its campaign.

To increase the effectiveness of crowdfunding and to help backers to make more informed decisions, some measures can be taken to mitigate the uncertainty around the quality of proposed products. 
One action could be to display more detailed information during the funding campaign.
For example, Kickstarter.com shares information on accumulated funds and number of backers. Moreover, their website allows backers to communicate through comments. On Crowdfunder.co.uk, visitors can see the entire timeline of the funding raised on every project, which shows the profile of backers and how much they contributed to projects. Revealing the entire process in a transparent manner allows backers to be involved in \textit{observational learning} and use the decisions of earlier backers to infer the quality of the products. This can alleviate the information asymmetry on quality \citep{Belleflamme2015}. While observational learning is increasingly facilitated by crowdfunding platforms, its impact on crowdfunding effectiveness remains a relatively understudied area.

At any given moment, numerous live campaigns can be found, even within highly specific product categories and subcategories. For example, at the time of writing on 25 September 2023, a search query made on Kickstarter.com for live New York-based campaigns proposing ``tabletop games'' listed 14 projects. Since the availability of many similar crowdfunding projects significantly affects the pledging decisions of backers, in this paper we investigate the role of \textit{observational learning} (OL) in the effectiveness of reward-based crowdfunding platforms where multiple projects compete for funding from backers. 

In our investigation, we use a stylized model where two projects, which might differ in quality, launch at the same time with the same funding target and the same pledge amount option that promises one unit of the final product. This specific setting allows us to focus on the influence of project quality and the inherent uncertainty associated with it when backers engage in observational learning, while other factors that might affect backers' preference for projects are assumed identical. Following the convention of earlier observational learning models with rational and Bayesian belief updates \citep{Banerjee1992, Bikhchandani1998, Zhang2015, Qiu2017, Liu2022}, we assume that a number of backers arrive sequentially to a platform, each receiving their own private signals on project quality independently and observing the pledge decisions of earlier backers. 
In this study, we extend OL modeling to a situation where learning involves more than one product or project. 
We consider multiple projects to analyze the impact of competition and derive insights into how funding scarcity changes the crowdfunding dynamics. 
Moreover, we differentiate between backers based on their level of expertise in accurately assessing product quality, as explored in the empirical study by \citet{expertness1}, and provide new theoretically-driven insights on how the difference in \textit{expertness} levels among early and late backers affect crowdfunding outcomes.

Our first goal is to understand how observational learning impacts crowdfunding performance in platforms with multiple projects. For this, we utilize our OL model and examine crowdfunding dynamics under OL in a two-backer, two-project system, comparing it to a no-learning scenario. We investigate various performance measures that are relevant to backers, creators, and crowdfunding platforms. 
Specifically, we focus on the impact of OL on \begin{enumerate*}[label = \itshape(\roman*)] \item the backers' contentedness with their pledging decisions, \item the success chances of project creators, \item the platform's short-term profit, and lastly on \item the platform's effectiveness, which is an indicator of its long-term profit. \end{enumerate*} We provide numerous new insights into the (dis)advantages of observational learning and the extent to which crowdfunding processes should be transparently disclosed. Our results also shed light on the specific factors (e.g., the scarcity of funding, the quality differences among proposed projects) that determine when observational learning can deliver benefits. We validate these results on large systems using simulations.

Our second goal is to investigate how OL impacts the product quality decisions of project creators. To this end, we investigate the product quality decisions in the presence and absence of OL, and show that OL creates an incentive for developing high-quality projects, in particular under fund scarcity. 

We summarize the main contributions of this paper in two categories.
\begin{itemize}
    \item \textit{Modeling:} We extend observational learning modeling to situations that involve multiple projects or options. 
    \item \textit{New Insights for Reward-based Crowdfunding:} We explore important dimensions that play a role in controlling the impact of observational learning among backers on crowdfunding dynamics. Specifically, we highlight the significance of two key factors: the availability of funding from potential backers and the quality differences among proposed projects. Our results emphasize the importance of considering these aspects when designing policies related to information disclosure and determining the extent to which observational learning should be allowed or restricted within crowdfunding platforms. 
    
\end{itemize}

The rest of the paper is organized as follows. In Section \ref{sec:lit}, we review the related literature and highlight the contributions of our paper. Section \ref{sec:model} presents our crowdfunding model under observational learning with two competing projects. We use this model to evaluate the impact of observational learning in crowdfunding, from the perspectives of all parties involved, i.e., backers, project creators, and the platform. Section \ref{sec:meas} describes the specific measures we use to evaluate these. In Section \ref{sec:main_results}, we present our analytical results on the effect of observational learning on crowdfunding outcomes compared to the case where no learning takes place when the number of backers is limited to two. In Section \ref{sec:decs}, we investigate the equilibrium quality strategies for two competing projects. In Section \ref{sec:valid}, we consider systems with larger number of backers and use simulations to demonstrate the robustness of our main results. We conclude the paper in Section \ref{sec:conc}.

\section{Literature Review} \label{sec:lit}
Online reward-based crowdfunding platforms have been the subject of both empirical and theoretical research in recent years, mainly for the purpose of providing decision support to project creators and platforms. We classify this research stream into two categories: studies focusing on \begin{enumerate*}[label = \itshape(\roman*)] \item performance estimation and \item those providing decision-making support. We provide a summary of the literature in Table \ref{tab:lit}. \end{enumerate*}

\afterpage{
\renewcommand{\arraystretch}{1.2}
\begin{landscape}
\begin{table}
  \caption{Theoretical and Empirical Reward-based Crowdfunding Studies Providing Decision Support on Performance Estimation (PE) and Decision-Making (DM)} \label{tab:lit}
\small 
    \begin{tabular}{|l|cccccccc|}
    \multicolumn{9}{p{22cm}}{{\scriptsize (KiA: Keep-it-All, $1^{PE}:$ Estimating performance for creators' success chances, $2^{PE}:$ Estimating performance for platforms' profits, $3^{PE}:$ Estimating performance for backers' contentedness, $1^{DM}:$ On reward pricing decisions, $2^{DM}:$ On product quality decisions, $3^{DM}:$ On funding target decisions, $4^{DM}:$ On campaign duration decisions, $5^{DM}:$ On decisions over the funding mechanism, $6^{DM}:$ On decisions over the information disclosure strategy, $7^{DM}:$ On dynamic timing of promotion actions), OL: Observational learning, SL: Word-of-mouth-based social learning.}} \\
    \hline
    & & \textbf{Funding}  & \textbf{Pledging}  & \textbf{Learning} & \textbf{Backer} & \textbf{Inform.} & & \textbf{Analysis} \\
    \multicolumn{1}{|c|}{\textbf{Paper}} & \textbf{Model Type} & \textbf{Mechanism} & \textbf{Opport.}  & \textbf{Among}  & \textbf{Heterogeneity} & \textbf{Asym.} & \textbf{Competition} & \textbf{Provided} \\
    & & & \textbf{Costs} & \textbf{Backers} & & & & \\
       \hline
       \citet{MarkCF} & two-stage game & AoN & no & no & quality valuation & no & no & $1^{DM}$, $2^{DM}$ \\
       \hline
       \citet{RW} & dynamic & AoN & yes & no& quality valuation & no & no & $1^{PE}$ \\
       \hline
       \citet{bi} & two-stage game & AoN and KiA & no & no & quality valuation & no & no & $1^{DM}$, $5^{DM}$ \\
       \hline
       \citet{expertness1} & empirical & other & - & OL & expertness & yes & no & $1^{PE}$ \\ 
       \hline
       \citet{cornel} & empirical & AoN & - & OL and SL & experience & no & no & $1^{PE}$ \\
       \hline
       \citet{Li20} & empirical & AoN & yes & SL & quality valuation & no & no & $1^{PE}$, $7^{DM}$\\
       \hline
       \citet{peng} & two-stage game & AoN & no & no & package size valuation & no & no & $1^{DM}$, $3^{DM}$ \\
       \hline
       \citet{can} & two-stage game & AoN & no & no & quality valuation & no & no & $2^{DM}$ \\
       \hline
       \citet{Chakraborty2021} & game & AoN& no & no & informed or not & yes & no & $1^{DM}$, $3^{DM}$ \\
       \hline     
       \citet{Cong2017} & dynamic & AoN & no & OL & no & yes & no & $1^{DM}$, $3^{DM}$ \\
          \hline
          \citet{Li21} & game & AoN and KiA & no & no & informed or not &yes & yes & $1^{DM}$, $5^{DM}$ \\
          \hline
       \citet{Liu21} & two-stage game & AoN & no & no & quality valuation & no & no & $1^{DM}$, $2^{DM}$, $3^{DM}$ \\
       \hline
       \citet{momentum} & empirical & AoN& - & OL & no & no & no & $1^{PE}$ \\
       \hline
       \citet{efficiency} & game &  AoN and KiA & no & no & no & no & no & $5^{DM}$ \\
       \hline
       \citet{du} & dynamic & AoN & yes & no & willingness to pledge & no & no & $7^{DM}$ \\
       \hline        
       \citet{Liu2022} & two-stage game &  AoN & no & OL & no & yes & no & $1^{DM}$, $6^{DM}$ \\
         \hline
         \citet{revenue} & dynamic & AoN & no & no & quality valuation & no & no & $3^{DM}$, $4^{DM}$ \\
         \hline
       \citet{chak2} & two-stage game &  AoN & yes & no & pledge delaying costs & no & no & $1^{DM}$, $3^{DM}$ \\
       \hline
       \textbf{Our paper} & dynamic & AoN & no & OL & expertness & yes & yes & $1^{PE}$, $2^{PE}$, $3^{PE}$,  \\
       & & & & & & & & $2^{DM}$, $6^{DM}$ \\
       \hline
    \end{tabular}
\end{table}
\end{landscape}}

The empirical studies with a performance estimation focus mainly explore the factors affecting the funding dynamics and the success chances of projects, while the theoretical ones mainly use dynamic models to explain backers' sequential decisions. On the other hand, the studies with a decision-making focus provide decision tools for project creators who are faced with a multiplicity of decisions. 
Some of these decisions may involve campaign design (e.g., selecting the platform, pricing of their product rewards, design of the product menus, setting a funding target and campaign duration) or product development (e.g., investments in technology, material), which would affect product quality, and need to be made prior to the campaign launch. 
In addition to these one-time decisions, project creators sometimes are also faced with dynamic decision-making when their campaigns are ongoing, such as the timing of promotion strategies (e.g., posting updates about the project's progress on the campaign's webpage, responding to backers' comments). 
Using our theoretical dynamic model with observational learning among backers, we provide decision support on both performance estimation and decision-making, however, with a particular emphasis on performance estimation. 
As in many earlier studies, we focus on performance measures related to projects' success chances and platforms' profits, which are governed by the funding raised during the campaigns. 
However, differently from those studies, we also consider post-campaign outcomes. 
Specifically, we measure backers' contentedness with their pledging decisions after the campaign finalizes, and thereupon products are delivered if the funding goal is reached, and the true quality of the products are observed.

In the literature, although learning among backers is commonly considered in empirical studies, theoretical studies rarely incorporate this aspect. Empirical studies often capture word-of-mouth type social learning (SL) among backers with social media data (e.g., number of Facebook shares) and observational learning (OL) with the campaigns' dynamic funding data which usually track the total number of backers or amount of funding raised. 
We note two theoretical crowdfunding studies with OL \citep{Cong2017, Liu2022} similar to ours. As in these studies, we consider the information asymmetry around the true qualities of the proposed projects (between project creators and backers), and backers engage in observational learning to infer the true qualities from the decisions of earlier backers. It is important to note that the theoretical models in these studies, as well as ours, feature a much more extensive set of observable information for backers to learn from earlier backers (e.g., expertness levels -- the extent to which each backer's signal on project quality is accurate) than is typically considered in empirical studies with OL. 
Although the expertness levels are modeled in \citet{Cong2017} and \citet{Liu2022}, all backers are assumed to have the same level of expertness. In contrast, we relax this assumption in our model and distinguish backers with respect to their expertness levels, similar to the empirical study in \citet{expertness1}. By considering the backers' differences in expertness levels when explaining crowdfunding dynamics, we derive new insights, which are explained in comparison to those found in the empirical work in \citet{expertness1} in Section \ref{sec:regret_results}.  
 
One important novelty in our model is that we do not consider the crowdfunding projects in isolation as if they were monopolies; rather, we consider them in competition. Previous research in crowdfunding has largely ignored the competitive nature in crowdfunding platforms, where hundreds of active campaigns may be running simultaneously, often competing with each other. We note only one study that considers two competing projects \citep{Li21}, and similarly to ours, also includes information asymmetry on quality. However, the way this asymmetry is handled in \cite{Li21} is different from our OL approach with varied expertness levels, and more akin to the approach in \citet{Chakraborty2021}, where some backers are \textit{fully informed} and can identify the true quality with certainty, while others are \textit{uninformed}, meaning they have no intuition at all about the true project quality. In our learning model, we also allow the expertness level to fall between these two extreme scenarios. 

Given the limited funding potential from backers and the fact that not all projects can be successful, competition plays a crucial role in crowdfunding dynamics and can significantly impact the performance of crowdfunding campaigns. This competition can affect various aspects of crowdfunding, including the success chances of projects seeking funding from the same pool of backers. Our results in this paper confirm the impact of competition on crowdfunding outcomes, as we show how the nature of competition affects the impact of OL on crowdfunding dynamics. It is also worth noting that in OL models in general, which do not focus specifically on crowdfunding, learning typically involves a single option/product/project, and consumers form beliefs about its state. Therefore, another contribution of our study is to address observational learning in competition.

\section{Crowdfunding Modeling under Observational Learning with Competing Projects} \label{sec:model}

We consider a reward-based crowdfunding platform where backers visiting the platform can see the decisions of earlier backers on ongoing campaigns, along with the description of the product that is provided by the campaign creator. On such a platform, backers looking for a particular product would visit the webpages of several campaigns of interest, read their descriptions, observe the decisions of earlier backers, and then finally make their decisions on which campaign to support. 
At a given time, there may be several projects in similar categories to choose from. Therefore, to capture the pledge decisions of backers, the multiplicity of the available projects to support should be incorporated. For tractability, we restrict our model to a situation in which two creators, denoted by index $i=1,2$, who seek to promote similar products, launch their campaigns simultaneously on the platform for the same duration of time with the same funding target, in accordance with common practice \citep{revenue}.
A number of backers who are interested in pledging, denoted by index $j=1,2,...,N$, visit the platform sequentially. We assume that the number of backers required to achieve the funding target, denoted with $\Tilde{n}$, is not greater than the number of participating backers, $N$, to avoid trivial cases.

To understand the role of the product quality differential, especially when the true qualities are not perfectly known to backers, we allow project campaigns to have different qualities. Let $V^i$ denote the true quality value of project $i$. Following the convention in \citet{Zhang2015}, the quality of a project is either high or low, which are represented with 1 and 0, respectively (i.e., $V^i\in \{0,1\}$). Let 
$\cV=\{ (V^1=1,V^2=1), (V^1=1,V^2=0), (V^1=0,V^2=1), (V^1=0,V^2=0)\}$ denote the set of all possible true (or actual) quality states for two projects. 

We assume that both projects demand the same pledge amount, and therefore, the choice of backers is directly governed by how they judge the quality of each project. The requested \textit{pledge amount}, also referred to as \textit{price}, falls between zero and one. If a high-quality project is successful, each backer receives a payoff of one; and, in the case of a low-quality project being successful or a project not meeting the funding target, the backers receive a payoff of zero. Thus, only when they pledge to a high-quality project backers can obtain positive utility. 
Accordingly, backers would not be content with their decisions in two cases: \begin{enumerate*}[label = \itshape(\roman*)] \item when they pledge to a low-quality project, and \item when they decide not to pledge at all while at least one of the projects is high-quality.\end{enumerate*} It must be noted that in our model, similarly to \citet{Chakraborty2021} and \citet{MarkCF}, we do not take into account the opportunity costs that might occur when backers pledge to a project that fails to raise the target funding. As argued in \citet{Chakraborty2021}, this assumption is not very restrictive in the context of reward-based crowdfunding since backers are fully refunded in case of campaign failure.

We denote the pledging decision of backer $j$ with $x_j \in \{0,1,2\}$, where $x_j=i$ for $i\in\{1,2\}$ implies that backer $j$ supported project $i$, and $x_j=0$ implies that the backer left without investing. 
When backer $j$ visits the platform and reviews the project campaign pages, she receives a private quality signal regarding the quality of project $i$, denoted by $s_j^i$, which is independent for projects. A backer's private signal for a project can indicate either high ($H$) or low quality ($L$), and may or may not reflect the true quality state of the project, $V^i$. Let $\cS = \{ (H,H), (H,L), (L,H), (L,L) \}$ denote the set of all possible signal pairs that a backer can receive for two projects. For example, $(H, L)$ implies that according to backer $i$'s signal, the first project is high-quality and the second project is low-quality. We assume that backers are unbiased, that is, before receiving their private signal, they have no prior beliefs.

In addition to her private signal, the backer also observes the decisions of all previous backers. Hence, the information set of backer $j$ prior to making her pledging decision, denoted as $I_j$, is $ \{\bs{j}, \bx{j-1} \} $, where $\bs{j} = (s_j^1, s_j^2)$ denotes backer $j$'s private signals for the two projects and $\bx{j-1} = (x_1, \dots, x_{j-1})$ denotes the decisions of the previous backers.
 
An expert backer has a high probability of accurately assessing the quality of a proposed project, and therefore is likely to receive a high (low) private quality signal when the true quality of the project is high (low). Let $p_j$ denote the probability that backer $j$'s private signal for a project is the same as the true quality, as shown below. \vspace{-5pt}
 \begin{align}
    p_j=P(s_j^1=H| V^1=1)=P(s_j^1=L|V^1=0)=P(s_j^2=H| V^2=1)=P(s_j^2=L|V^2=0), \\
    1-p_j=P(s_j^1=H| V^1=0)=P(s_j^1=L|V^1=1)=P(s_j^2=H| V^2=0)=P(s_j^2=L|V^2=1).
\end{align}

Without loss of generality, we assume $p_j \in[0.5, 1]$. Larger values of $p_j$ imply higher \textit{expertness} levels, and $p_j=0.5$ represents a completely uninformed signal.

Although early backers' private signals cannot be observed by other backers, we assume that the distribution of private signals is known, allowing backers to update their beliefs on project quality through Bayesian rationality. Furthermore, in some crowdfunding platforms, backers can access the profiles of previous backers, enabling them to deduce the expertness of others by examining their history or looking for high-expertness labels provided by the crowdfunding platform, e.g., the ``superbacker" label on Kickstarter.com. Bayesian rationality assumption is also supported by empirical evidence suggesting that even non-expert backers are sophisticated in their ability to identify and exploit the expertness of early backers \citep{expertness1}. Moreover, the entirety of our modeling assumptions regarding the binary private signals and the incorporation of earlier decisions comply with the typical observational learning models studied in the literature \citep[see][]{Banerjee1992, Bikhchandani1998, Zhang2015, Qiu2017}.

Under Bayesian updating, the posterior belief of the $j^{\text{th}}$ backer with the information set $I_j=\{\bs{j}, \bx{j-1} \}$ that the true quality state is $\bV = (V^1, V^2) \in \cV$ is calculated as follows
\begin{equation}
    P(\bV | I_j) = P(\bV | \bs{j}, \bx{j-1}) = \frac{P(s_j^1|\bV) \; P(s_j^2|\bV) \; P(\bx{j-1}|\bV)}{\sum_{\bV' \in \cV} \, P(s_j^1| \bV') \; P(s_j^2| \bV') \; P(\bx{j-1}| \bV')},
\end{equation}
where \vspace{-10pt}
\begin{align*}
   P(\bx{j-1} | \bV) &= P(\bx{j-2}| \bV) \; P(x_{j-1}| \bx{j-2}, \bV) \\
   & = P(\bx{j-2} | \bV) \sum_{\bs{j-1} \in \cS} P(x_{j-1} | \bs{j-1}, \bx{j-2}) \; P(s_{j-1}^1|\bV) \; P(s_{j-1}^2|\bV).
\end{align*}

Backer $j$ makes her funding decision based on the posterior beliefs $P(\bV|I_j), \, \bV\in \cV$. Many studies using observational learning models in the literature have assumed that these decisions follow a \textit{pure strategy} that chooses the option with the highest expected utility \citep[see][]{Banerjee1992, Bikhchandani1998, Acemoglu2011, Zhang2015, Qiu2017, Liu2022}. These studies consider only one project; therefore, the backers' only choice is whether to pledge or not. On the other hand, in our system, backers have three options to choose from; and, the backer choice closely resembles the demand models in the assortment literature \citep{Besbes16}. Therefore, instead of a pure strategy, we consider a \textit{probabilistic pledging strategy}, as in \citet{revenue}, that is similar to the attraction-type demand models \citep{Huang2013}, such as the multinomial logit model. Under this probabilistic mechanism, backer $j$'s decision is governed by the following probabilities. \vspace{-5pt}
\begin{align}
    P(x_j=1| I_j)=&P(V^1=1, V^2=0| I_j)+0.5 P(V^1=1, V^2=1| I_j)\\
     P(x_j=2| I_j)=&P(V^1=0, V^2=1| I_j)+0.5 P(V^1=1, V^2=1| I_j)\\
      P(x_j=0| I_j)=&P(V^1=0, V^2=0| I_j)
\end{align}
These probabilities reflect the utility maximization objective in a similar fashion to that of a multinomial logit model. Note that in the case where both projects are high-quality, backers would be indifferent between them in terms of the utility they will obtain from their products. Therefore, the weight of this belief is symmetrically distributed between the two projects.

So far, we have modeled backers' decision mechanisms considering the observational learning effect and their own private signal when faced with multiple project options based on the observations from the literature. 
In the next section, we introduce new metrics to measure the performance of the crowdfunding environment and discuss how to evaluate platforms.

\section{Evaluating Crowdfunding Performance}
\label{sec:meas}

To establish a comprehensive understanding for the impact of observational learning in crowdfunding, we first introduce some metrics to evaluate the system outcomes from the perspectives of all parties involved, i.e., backers, project creators, and the crowdfunding platform. Then, we discuss how we evaluate different outcomes considering these metrics.

To track the outcomes of the crowdfunding, we let 
\begin{equation}
    n_i (\bVp) := \sum_{j=1}^N \mathbbm{1}_{\{ x_j=i | \bVp \}},\quad i=1,2
\end{equation}
denote the random variable for the total number of backers pledging to project $i$ when the true quality state is $\bV\in \cV$ and the backer expertness level is $\bp = (p_1, p_2, \dots, p_N)$.

\subsection{Evaluation Metrics}
\label{sec:metrics}

From the perspective of the backers, we define post-campaign contentedness. \emph{Backer contentedness} is defined as the expected number of backers who have no regret about their decisions  which could have happened due to missing a high-quality project opportunity or pledging to a low-quality project. Let $\content(\bVp) \in [0,N]$ denote the backer contentedness measure when the true quality state is $\bV$ and the backer expertness level is $\bp = (p_1, p_2, \dots, p_N)$.

\begin{equation}
    \content(\bVp) := \begin{cases}
    \mathbb{E}\left[n_1 (\bVp) \right] + \mathbb{E}\left[n_2 (\bVp) \right], & \bV=(1,1)\\
    \mathbb{E}\left[n_1 (\bVp) \right], & \bV=(1,0)\\
    \mathbb{E}\left[n_2 (\bVp) \right], & \bV=(0,1)\\
    N - \left( \mathbb{E} \left[ n_1 (\bVp) \right] + \mathbb{E} \left[ n_2 (\bVp) \right] \right),\; & \bV=(0,0)
    \end{cases}
    \tag{Backer contentedness}
\end{equation}
When both projects are high-quality, all backers who pledge will be content with their decisions. When only one of the projects is high-quality, only the backers supporting that project will be content. If both projects are low-quality, backers will be content only if they decide not to pledge for either project.

For project creators, crowdfunding performance is measured by the probability to meet the funding target. In AoN funding systems, projects must raise funding from at least $\Tilde{n}$ backers to reach the target and collect the funds. We use $\success_i (\bVp) \in [0,1]$ to denote the \emph{success probability} of project $i$ at quality state $\bV$ and when the backer expertness level is $\bp$, and define it as follows: 
\begin{equation}
    \success_i (\bVp) := P \Big( n_i (\bVp) \geq \Tilde{n} \Big),\quad i=1,2.
    \tag{Success probability}
\end{equation}

From the perspective of the crowdfunding platform, we define two performance measures. The first one, \emph{platform profit}, has a short-term focus as it reflects the revenue earned by the platform from service fee charges on the total funds raised from successful projects\footnote{https://tarrida.co.uk/rewards-based-crowdfunding/} \citep{Li21}. Therefore, both the total funding collected from backers and the success chances of projects are important for the platform.  We denote the platform's profit measure for a true quality state $\bV$ and expertness level $\bp$ by $ \profit (\bVp) \in [0, \gamma N]$, where $\gamma \in[0,1]$ denotes the service fee rate. $ \profit (\bVp)$ is calculated as follows:

\begin{equation}
    \profit (\bVp) := \gamma \sum_{i \in \{1, 2\}} \mathbb{E} \left[ n_i (\bVp) \right] \; \success_i (\bVp).
    \tag{Platform profit}
\end{equation}

The second metric we define from the platform's perspective has a long-term focus. For crowdfunding platforms to remain attractive for the backers in the long run, it is important that high-quality projects are funded. In a platform where low-quality projects become successful, backers may experience long delays in product delivery or receive unsatisfactory products. Such negative experiences could dampen their enthusiasm to use the crowdfunding platform, and result in abandoning crowdfunding altogether or switching to another platform. To address this long-term concern, we introduce \emph{platform effectiveness} as a metric to measure the platform's ability to attract funds to high-quality projects. We denote the platform's effectiveness under a quality state $\bV$ and expertness levels $\bp$ by $\effect (\bVp) \in [-N, N]$ and calculate it as follows:
\begin{equation}
    \effect (\bVp) := \sum_{i \in \{1, 2\}} \mathbb{E} \left[ n_i (\bVp) \right] \; \success_i (\bVp) \; \left( \mathbbm{1}_{\{V^i=1\}}-\mathbbm{1}_{\{V^i=0\}} \right).
    \tag{Platform effectiveness}
\end{equation}

\subsection{Evaluation Criteria}
Performance metrics are based on the quality state and the backer expertness. Since the latter is not known in advance, to provide more general insights, we introduce some evaluation criteria considering the entire feasible set of backer expertness.

When backers have no expertise on projects (i.e., expertness level is 0.5), or when they are fully informed about the projects (i.e., expertness level is 1.0), OL does not have any impact. In the case of no expertise, the late backer (he) can derive no useful information from the decision of an early backer (she), as he knows that she received a random, uninformative, signal and acted accordingly. In the case of fully informed backers, the ability of the late backer to rely completely on his own private signal, without needing the information on the decision of the early backer, renders OL ineffective. Therefore, to study the effectiveness of OL, we disregard these extreme cases and focus on cases where backers' expertness levels lie in between the two extremes.

To evaluate the impact of OL with respect to the performance metrics defined in Section \ref{sec:metrics} (i.e., backer contentedness ($\content$), success probabilities of the projects ($\bsuccess = (\success_1, \success_2)$), 
platform's profit ($\profit$), and platform effectiveness ($\effect$)), we compute these metrics with and without OL for given quality states, $\bV$, and backer expertness levels, $\bp = (p_1, p_2, \dots, p_N)$. Then, we determine backer expertness levels $\bp \in (0.5, 1)^N$ where observational learning has a positive or negative impact.

Let $\mu^{OL} (\bVp)$ and $\mu^{NL} (\bVp)$ denote the value of each performance metric $\mu \in \{ \content, \success_1, \success_2, \profit, \effect\}$ for a given quality state $\bV$ and expertness levels $\bp$, with and without OL, respectively. Furthermore, let $\imp_{\mu^+} (\bV)$ and $\imp_{\mu^-} (\bV)$ denote the sets of backer expertness levels, $\bp$, where OL has a positive and a negative impact, respectively, on the performance metric $\mu$.
\begin{align}
    \imp_{\mu^+} (\bV) &:=\Big\{ \bp \in (0.5, 1)^N | \mu^{OL} (\bVp) - \mu^{NL} (\bVp) \geq 0 \Big\}\\
    \imp_{\mu^-} (\bV) &:=\Big\{ \bp \in (0.5, 1)^N | \mu^{OL} (\bVp) - \mu^{NL} (\bVp) < 0 \Big\}
\end{align}

When there are some expertness levels where observational learning has a positive effect on performance metric $\mu$ in true quality state $\bV$, i.e., $\imp_{\mu^+} (\bV) \neq \emptyset$, then we say that OL has an \textit{improvement potential} for this particular metric under this particular quality state. Similarly, when $\imp_{\mu^-} (\bV) \neq \emptyset$, we say that OL has a \textit{harm potential}. Letting $|.|$ denote the cardinality of a set, when $|\imp_{\mu^+} (\bV)| \geq |\imp_{\mu^-} (\bV)|$, observational learning has a greater improvement potential than harm with respect to performance metric $\mu$ under quality state $\bV$.

In addition to improvement and harm potentials, we also quantify the magnitude of the improvement (or positive impact) and harm (or negative impact) OL can create. For a quality state $\bV$, we denote the \textit{maximum improvement potential} and \textit{maximum harm potential} of observational learning with respect to performance metric $\mu$ with $\dplus{\mu}{\bV}$ and $\dminus{\mu}{\bV}$, respectively.
\begin{align*}
    & \dplus{\mu}{\bV} := \max_{\bp \in \imp_{\mu^+} (\bV) }\Big \{ \mu^{OL} (\bVp) - \mu^{NL} (\bVp) \Big \} \tag{Max improvement potential of OL w.r.t. $\mu$} \\
    & \dminus{\mu}{\bV} := \min_{\bp \in \imp_{\mu^-} (\bV) }\Big \{ \mu^{OL} (\bVp) - \mu^{NL} (\bVp) \Big \} \tag{Max harm potential of OL w.r.t $\mu$}
\end{align*}

When all expertness levels are equally likely and independent from each other, i.e., when they are drawn from independent uniform distributions, we calculate the average impact of OL across all possible expertness levels as follows:
\begin{equation}
     \dave{\mu}{\bV} := \idotsint_{(0.5, 1)^N} 2^N \left( \mu^{OL} (\bVp) - \mu^{NL} (\bVp) \right) \,dp_N \dots \,dp_1  \tag{Average impact of OL w.r.t. $\mu$}
\end{equation}

A positive average impact with respect to $\mu$ in quality state $\bV$, i.e., $\dave{\mu}{\bV} > 0$, implies that on average OL improves the performance metric $\mu$ when the quality state is $\bV$. By comparing the maximum improvement and harm potentials, and the average impact of OL in all four quality states, we investigate how the qualities of competing projects affect the impact of OL. Moreover, we identify the quality states where OL has the highest improvement and harm potentials, as well as the largest average impact.   

Having formalized the evaluation method, in the following section, we analyze the impact of OL for different system attributes and provide insights.

\section{Impact of Observational Learning on Crowdfunding Dynamics} \label{sec:main_results}

In this section, we analyze the effect of OL on reward-based crowdfunding dynamics. For tractability, we consider two crowdfunding projects and limit the number of backers to two, as is common in reward-based crowdfunding \citep{MarkCF, Liu2022} and observational learning studies \citep{Qiu2017}. The detailed derivation of posterior beliefs and decision probabilities of this case are provided in the online supplement. 

In this small case, the first backer embodies the role of an ``early'', or a ``leader'', backer whereas the second backer takes the role of a ``late'', or a ``follower'' backer, since the second backer can observe the decision of the first and incorporate it into their own decision-making process. For clarity of exposition, we henceforth use the ``she'' pronoun to refer to the early backer, and ``he'' to refer to the late backer. 

To demonstrate the effect of observational learning, we compare the performance metrics when there is observational learning (denoted as OL), i.e., when backers observe and learn from the decisions of earlier backers, to those when there is no learning (denoted as NL).

In order to comparatively assess the impact of observational learning on creator- and platform-focused performance metrics, we consider environments where funding is scarce and abundant separately. For this purpose, we define two funding scarcity scenarios in a two-project, two-backer system. In the first scenario, which we refer to as the \textit{tight competition} case, projects can only reach their funding targets if both backers support them. In the second scenario, which we refer to as the \textit{relaxed competition} case, each project can achieve its target with just one backer.

\subsection{Impact of Observational Learning on Backer Contentedness} \label{sec:regret_results}

We characterize the effect of observational learning on backer contentedness in Proposition \ref{prop:content}.
\begin{proposition} \label{prop:content}
Observational learning has no potential to harm backer contentedness in any of the four quality states, i.e., $\imp_{\content^-} (\bV) = \emptyset,  \;\;\; \forall \bV \in \cV$, and it strictly improves backer contentedness, on average, in all four quality states, i.e., $\dave{\content}{\bV} > 0, \;\;\; \forall \bV \in \cV$. The maximum improvement potential and average impact of observational learning on backer contentedness in the each quality state are as follows: \vspace{-25pt}
\begin{multicols}{2}
\begin{equation}
\dplus{\content}{V}=\begin{cases}
			0.74, & V=(0,0)\\
                0.45, & V=(1,0), (0,1)\\
                0.25, & V= (1,1)
		 \end{cases}
   \nonumber
\end{equation}

\begin{equation}
\dave {\content}{V}=\begin{cases}
			0.08, & V=(0,0)\\
                0.05, & V=(1,0), (0,1)\\
                0.03, & V= (1,1)
		 \end{cases} \nonumber
\end{equation}
\end{multicols}
\end{proposition}

According to Proposition \ref{prop:content}, 
backers are better off with OL, and full information disclosure on other backers' pledging decisions is always preferred by backers. In particular, when the crowdfunding projects are low-quality, backers benefit considerably from OL to reduce their post-campaign regrets. This result quantifies the improved contentedness of backers when they choose crowdfunding platforms that disclose information in a transparent manner and allow observational learning as much as possible. 
\begin{figure}[h]
    \centering
         \begin{subfigure}[b]{0.32\textwidth}
            \centering
            \includegraphics[width=\textwidth]{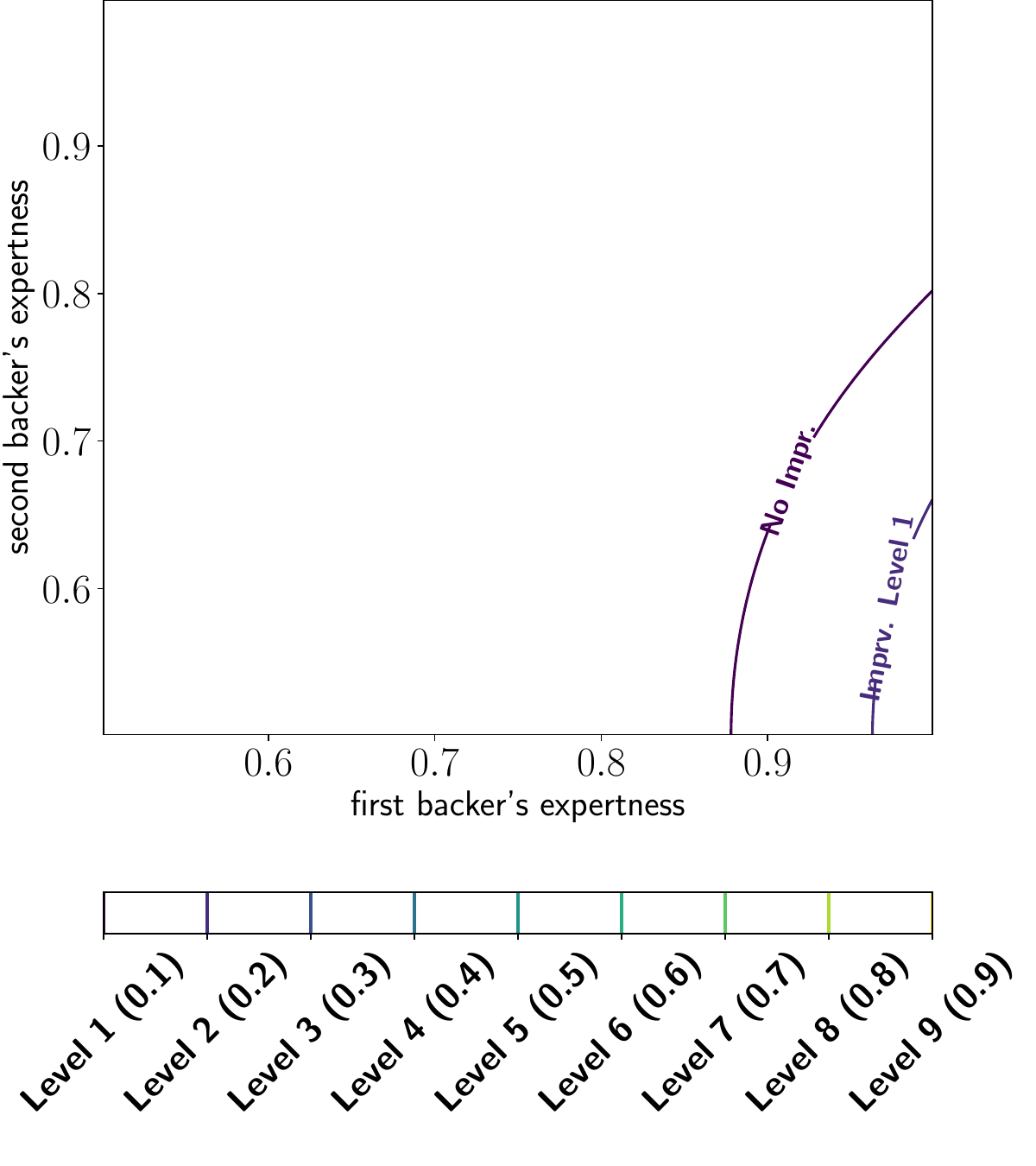}
            \caption{Both projects are high-quality}
            \label{fig:reg_11_projects}
        \end{subfigure}
          \begin{subfigure}[b]{0.32\textwidth}
            \centering 
            \includegraphics[width=\textwidth]{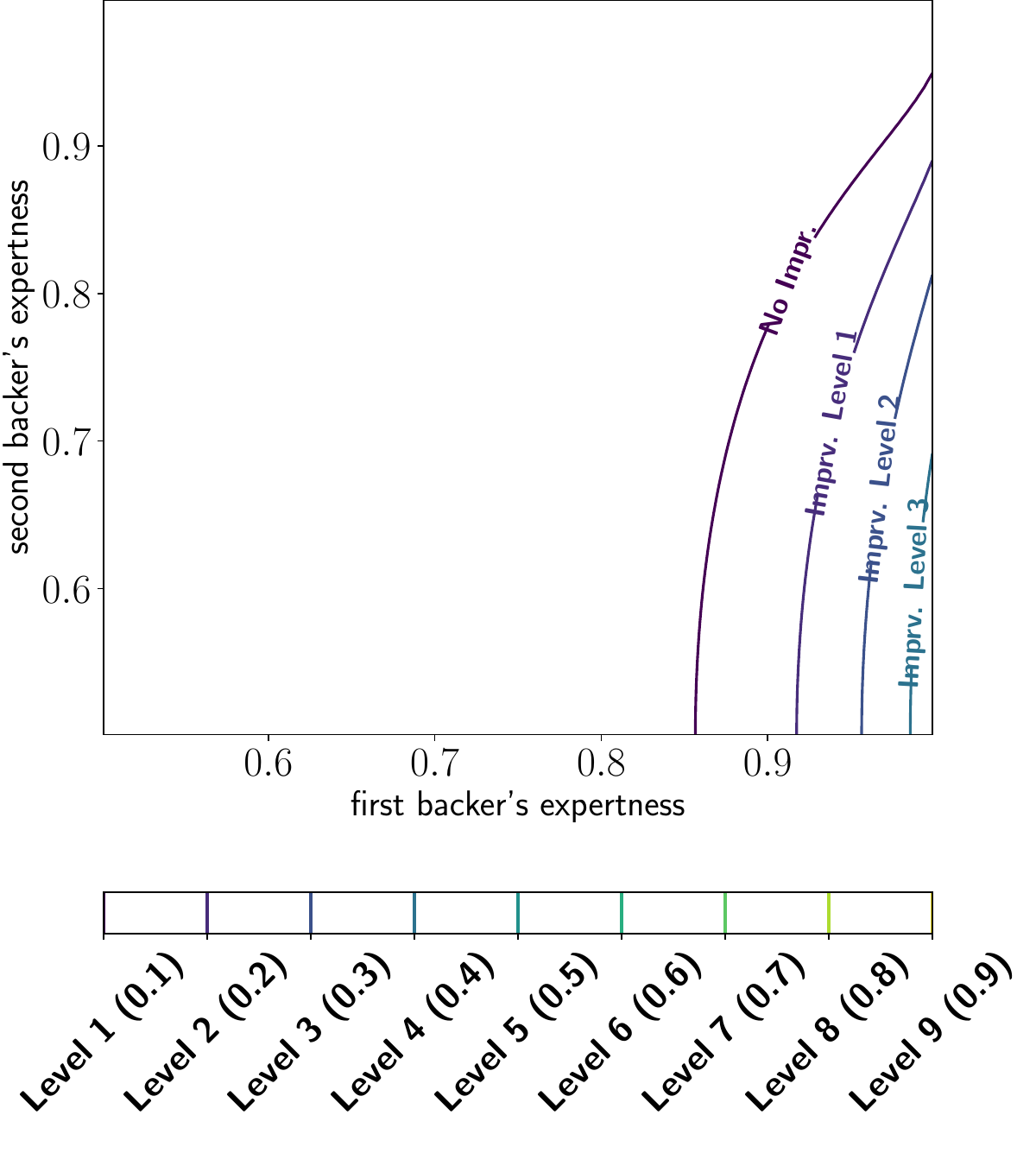}
            \caption{Only one project is high-quality}    
            \label{fig:reg_10_projects}
        \end{subfigure}
         \begin{subfigure}[b]{0.32\textwidth}
            \centering 
            \includegraphics[width=\textwidth]{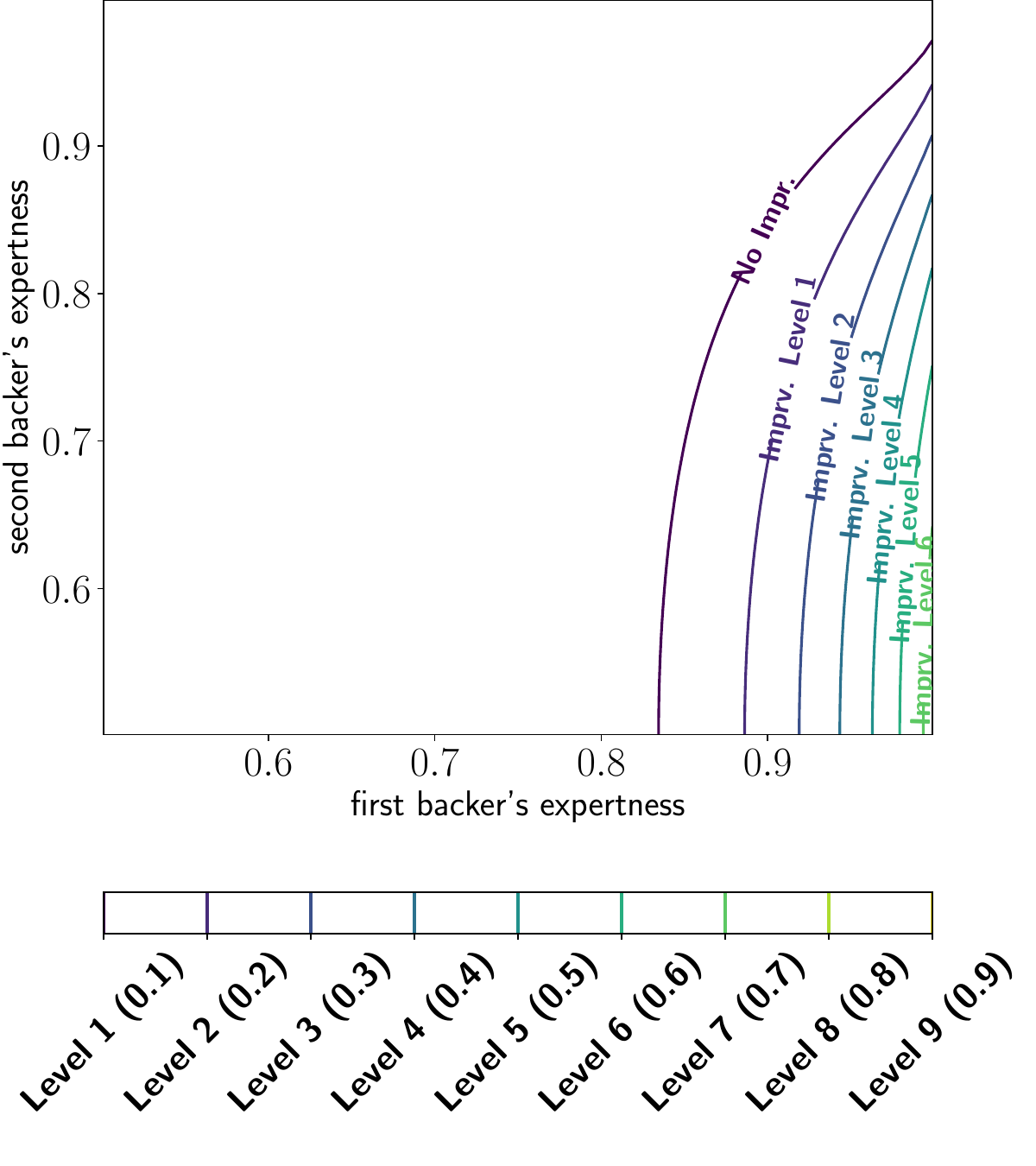}
            \caption{Both projects are low-quality}
            \label{fig:reg_10_projects}
        \end{subfigure}
        \caption{Positive Impact of OL on Backer Contentedness, $\content^{OL}(\bVp) - \content^{NL}(\bVp)$.}
    \label{fig:reg_projects}
\end{figure}

To better observe the impact of OL on backer contentedness, we present Figure \ref{fig:reg_projects}, which shows the magnitude of backer contentedness improvement that can be attributed to OL, for all project quality states and backer expertness conditions. It can be observed that OL achieves the largest improvement in backer contentedness when the early and late backers have high and low expertness levels, respectively. Under such expertness conditions, the late backer lacks the means to accurately evaluate the projects by using his own private signal. However, when there is OL, he can rely on the early backer's decision with high confidence. This way, OL can mitigate potential regret that might be experienced by a non-expert backer. Our results indicate that participation of expert backers during the early stages of crowdfunding would be important for subsequent backers. This is in line with the empirically-demonstrated results reported in \citet{expertness1}, which emphasizes the role of expert backers in online crowdfunding projects proposing mobile app development.

\subsection{Impact of Observational Learning on Creators' Crowdfunding Success} \label{sec:suc_prob_effect}

In this section, we assess the effect of OL on the success probability of projects in AoN reward-based crowdfunding campaigns. Without loss of generality, we focus on the success probability of the first project, $\success_1$, and denote it simply as $\success$ for brevity of exposition. We characterize the project's success probability in Proposition \ref{prop:success} and observe that the outcomes are highly dependent on the competition type. 

\begin{proposition} \label{prop:success}
In \textbf{tight competition}, observational learning has no potential to harm the creator's success probability in any of the four quality states, i.e., $\imp_{\success^-} (\bV) = \emptyset, \;\;\; \forall \bV \in \cV$, and it strictly increases the success probability, on average, in all four quality states, i.e., $\dave{\success}{\bV} > 0, \;\;\; \forall \bV \in \cV$. The maximum improvement potential and average impact of observational learning on success probability in each quality state are as follows: \vspace{-15pt}
\begin{multicols}{2}
\begin{equation}
\dplus{\success}{V}=\begin{cases}
			0.45, & V=(1,0)\\
                0.23, & V=(1,1)\\
                0.05, & V= (0,0)\\
                0.03, & V= (0,1)\\
		 \end{cases}
   \nonumber
\end{equation}

\begin{equation}
\dave {\success}{V}=\begin{cases}
			0.07, & V=(0,0)\\
                0.04, & V=(1,1)\\
                0.02, & V= (0,0)\\
                0.01, & V=(0,1)\\
		 \end{cases} \nonumber
\end{equation}
\end{multicols}

In \textbf{relaxed competition}, observational learning decreases the creator's crowdfunding success in the majority of expertness conditions in all of the four quality states, i.e., $| \imp_{\success^-} (\bV) | > | 
\imp_{\success^+} (\bV) |, \;\;\; \forall \bV \in \cV$, and it strictly decreases the success probability, on average, in all quality states, i.e., $\dave{\success}{\bV} < 0, \;\;\; \forall \bV \in \cV$. The maximum harm potential and the average impact of observational learning on success probability in each quality state are as follows: \vspace{-15pt}
\begin{multicols}{2}
\begin{equation}
\dminus{\success}{V}=\begin{cases}
			-0.05, & V=(1,0)\\
                -0.10, & V=(1,1)\\
                -0.21, & V= (0,1)\\
                -0.37, & V= (0,0)\\
		 \end{cases}
   \nonumber
\end{equation}

\begin{equation}
\dave {\success}{V}=\begin{cases}
			-0.02, & V=(1,0)\\
                -0.03, & V=(1,1)\\
                -0.04, & V= (0,1)\\
                -0.06, & V= (0,0)\\
		 \end{cases} \nonumber
\end{equation}
\end{multicols}

\end{proposition}

In the tight competition case, where projects require support from both backers to reach their funding targets, OL always has potential to benefit the creators. Nevertheless, the actual quality state of the proposed crowdfunding projects is an important determinant of the magnitude of this benefit. We find that OL increases the project's success probability the most when the project is high-quality, especially when it is competing against a low-quality project. 
On the contrary, in the relaxed competition case, projects have better success chances when there is no OL and each backer decides according to their own private signal. Although OL negatively impacts success probability in all quality states, we find, similarly to the tight competition case, that creators of high-quality projects will be better off, since they will be least affected by any potential harm that OL can inflict on their crowdfunding success. 
The reason for the remarkably different findings in the tight and relaxed competition cases is the \textit{herding effect} induced by OL on crowdfunding dynamics.  In tight competition, where the available funding potential is divided among projects, it is difficult for any of the projects to reach their funding targets. In this case, creators can benefit from the herding behavior, which concentrates available funding potential on one of the projects. Conversely, in relaxed competition, herding can be detrimental to project success. In this case, projects' success chances are higher if the available funding potential is split among projects rather than being directed towards overfunding one of the projects.

Our findings provide insights for project creators in making informed decisions prior to launching their projects. For instance, they can strategically choose between crowdfunding platforms that offer full OL, such as Crowdfunder.co.uk, or other less-transparent platforms, such as Indiegogo.com. To make this decision, creators should first analyze the funding landscape on the platforms, including factors such as timing of the campaign (e.g., simultaneous launch with similar high-quality products), the backer pool size (e.g., number of members), and the specific product category, as these can affect competition intensity. Historical data, such as success rates of similar products, can be used to determine the availability of funding. Based on our analysis in Proposition \ref{prop:success}, creators in tight competition, particularly those launching high-quality projects, should prefer platforms that facilitate OL and take actions to support learning during their campaigns. On the other hand, creators in relaxed competition, especially those with low-quality projects, may prefer platforms with limited OL. 

\begin{figure}[h]
    \centering
         \begin{subfigure}[b]{0.4\textwidth}
            \centering
            \includegraphics[width=0.8\textwidth]{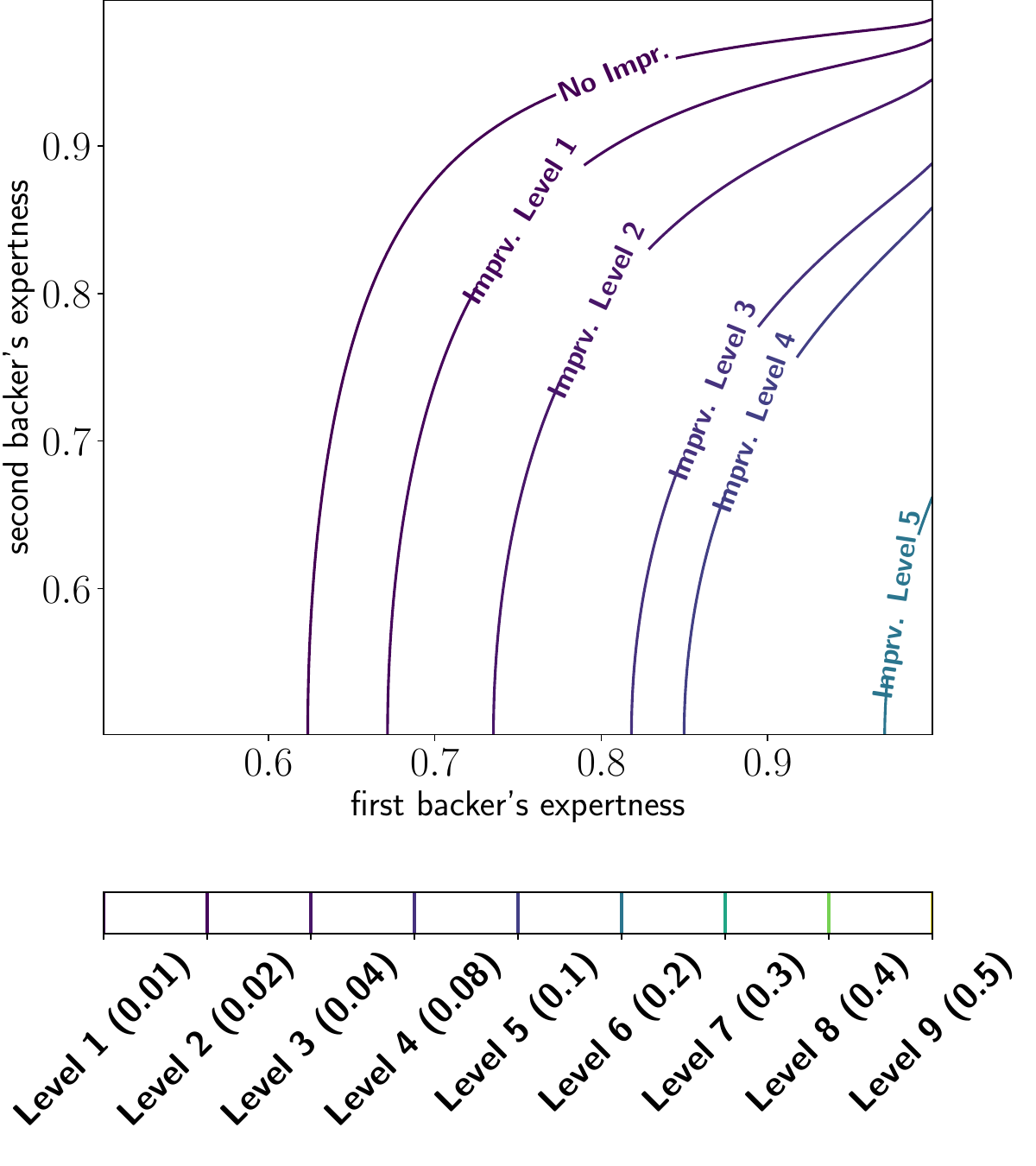} \vspace{-5pt}
            \caption[Network2]%
            {High-quality competing against high-quality} 
            \label{fig:reg_11_projects}
        \end{subfigure}
          \begin{subfigure}[b]{0.4\textwidth}
            \centering 
            \includegraphics[width=0.8\textwidth]{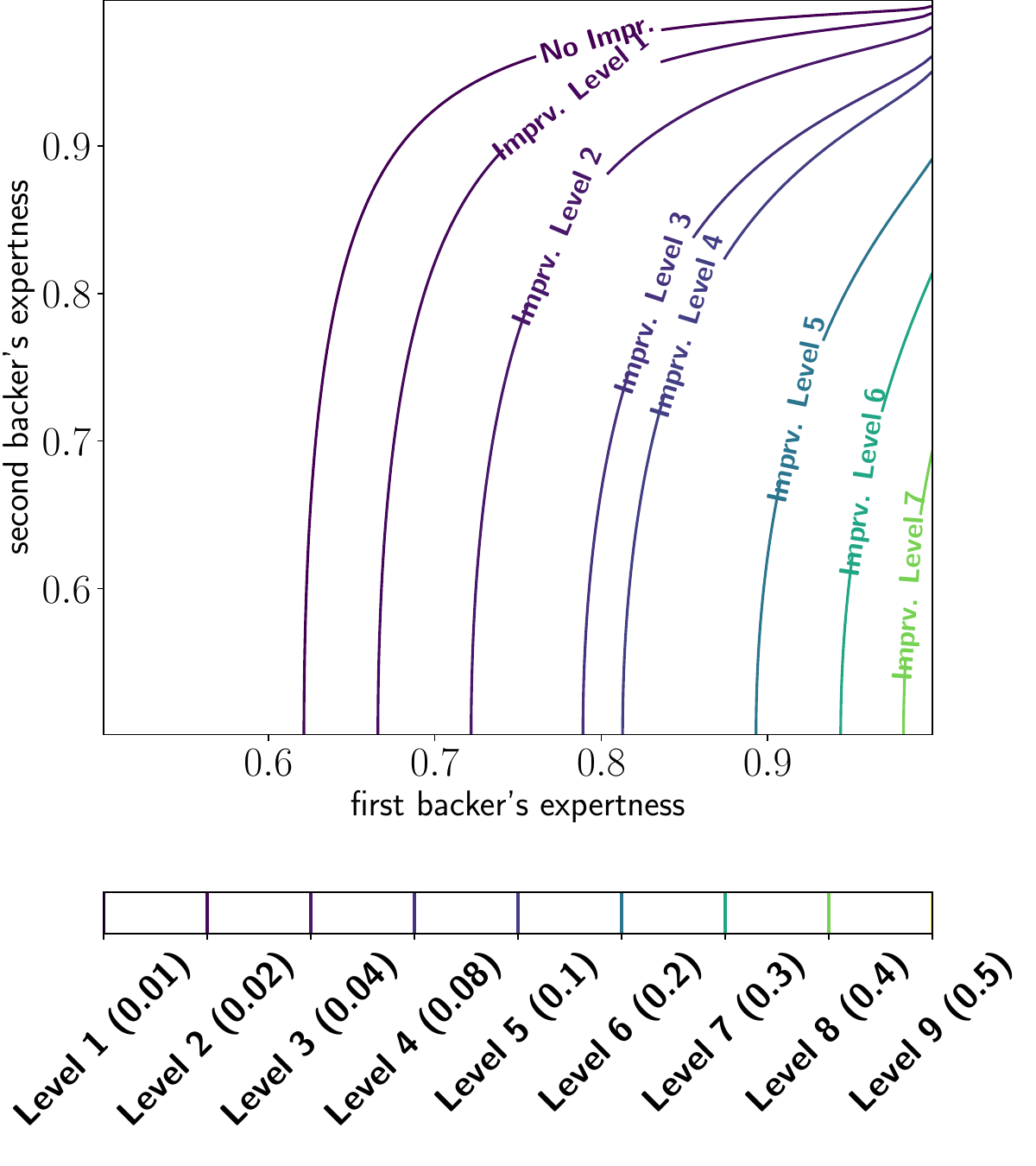}
            \vspace{-5pt}
            \caption{High-quality competing against low-quality} 
            \label{fig:reg_10_projects}
        \end{subfigure} \vspace{10pt}
        
        \begin{subfigure}[b]{0.4\textwidth}
            \centering 
            \includegraphics[width=0.8\textwidth]{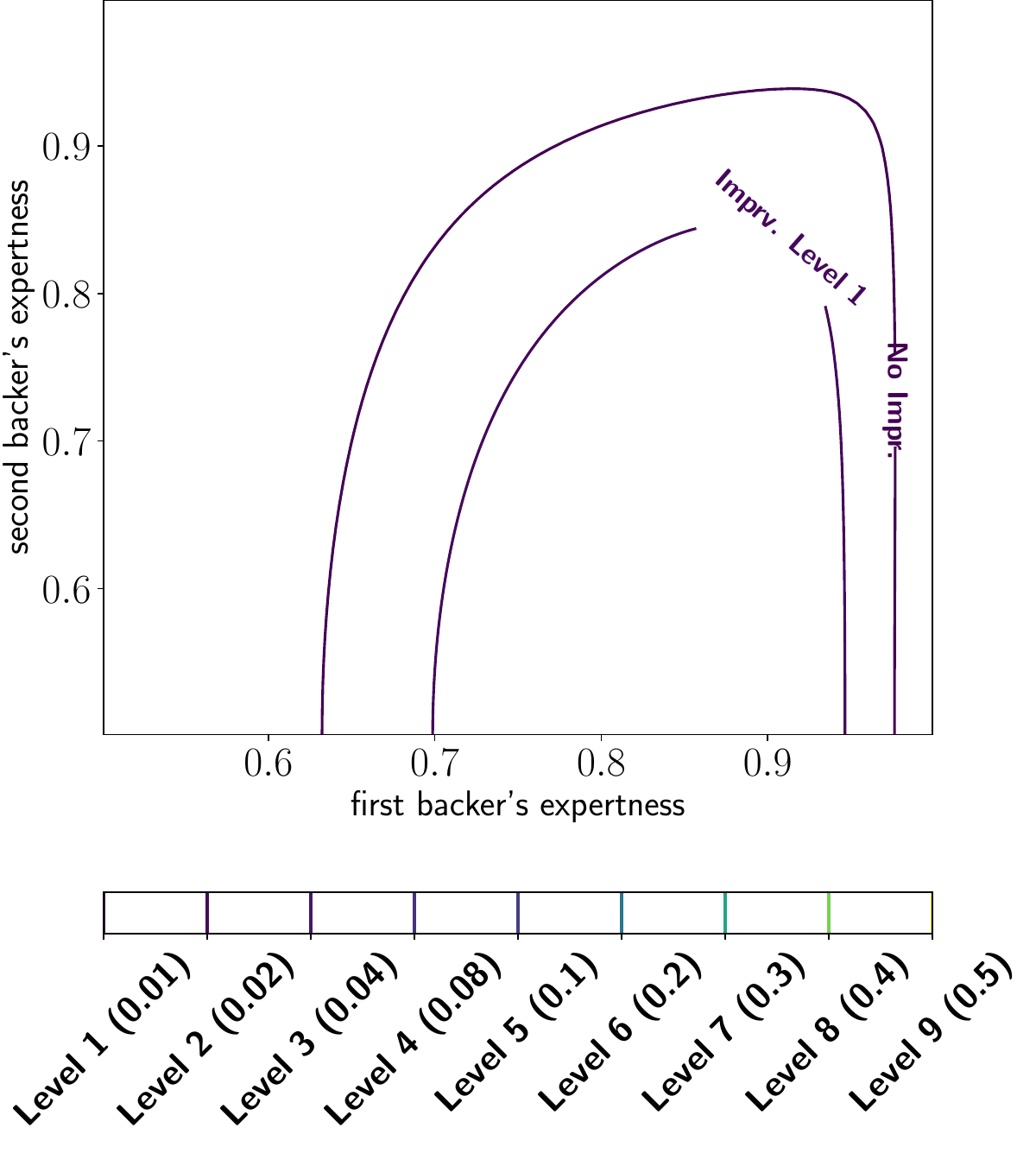} \vspace{-5pt}
            \caption{Low-quality competing against high-quality} 
            \label{fig:reg_10_projects}
        \end{subfigure}
         \begin{subfigure}[b]{0.4\textwidth}
            \centering 
            \includegraphics[width=0.8\textwidth]{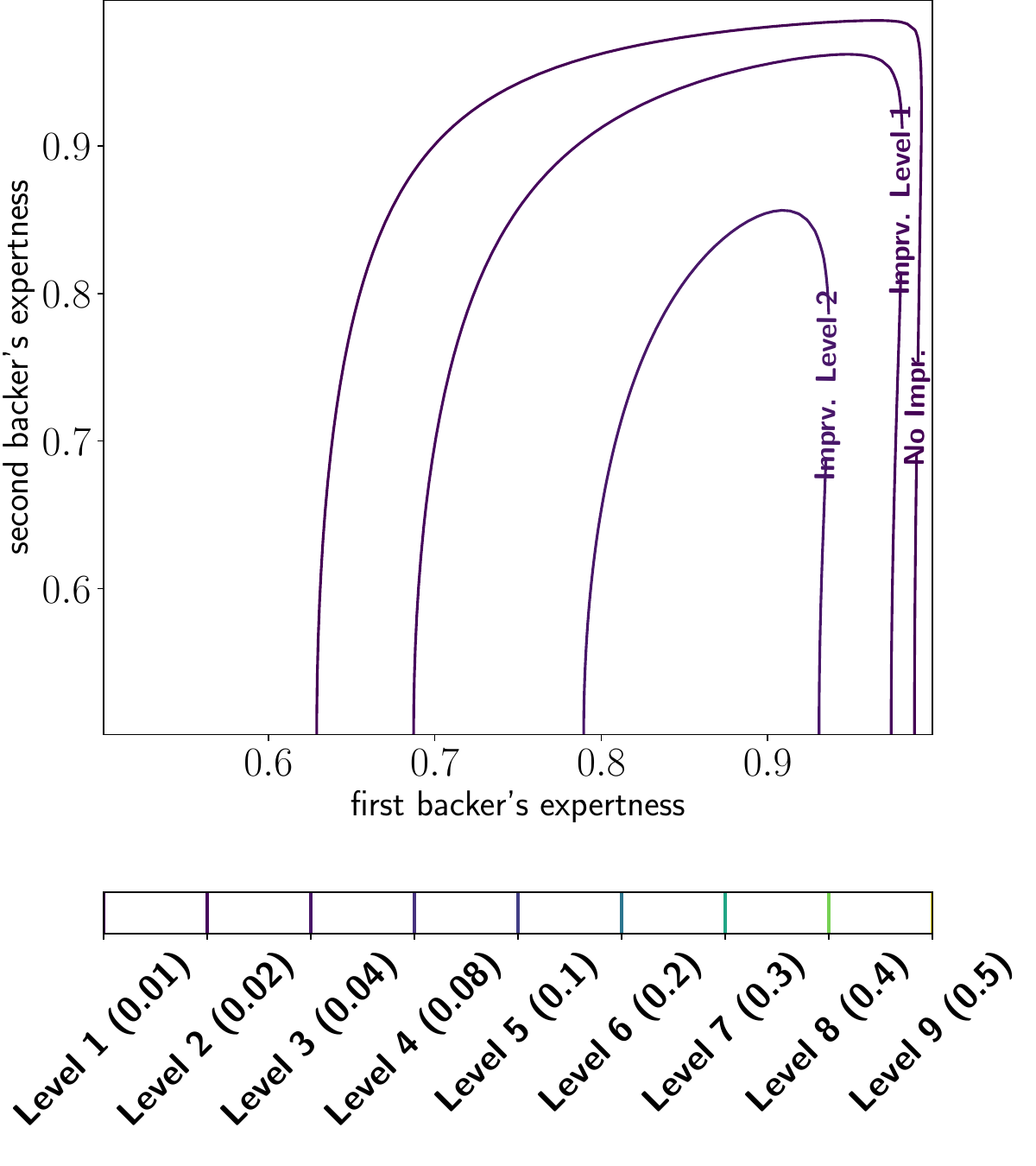} \vspace{-5pt}
            \caption{Low-quality competing against low-quality} 
            \label{fig:reg_10_projects}
        \end{subfigure}
        \caption{Improvement Potential of OL on Project Success Probability, $\success^{OL}(\bVp) - \success^{NL} (\bVp)$, in Tight Competition.}
    \label{fig:succ_projects_tight}
\end{figure}

\begin{figure}[h]
    \centering
         \begin{subfigure}[b]{0.4\textwidth}
            \centering
            \includegraphics[width=0.8\textwidth]{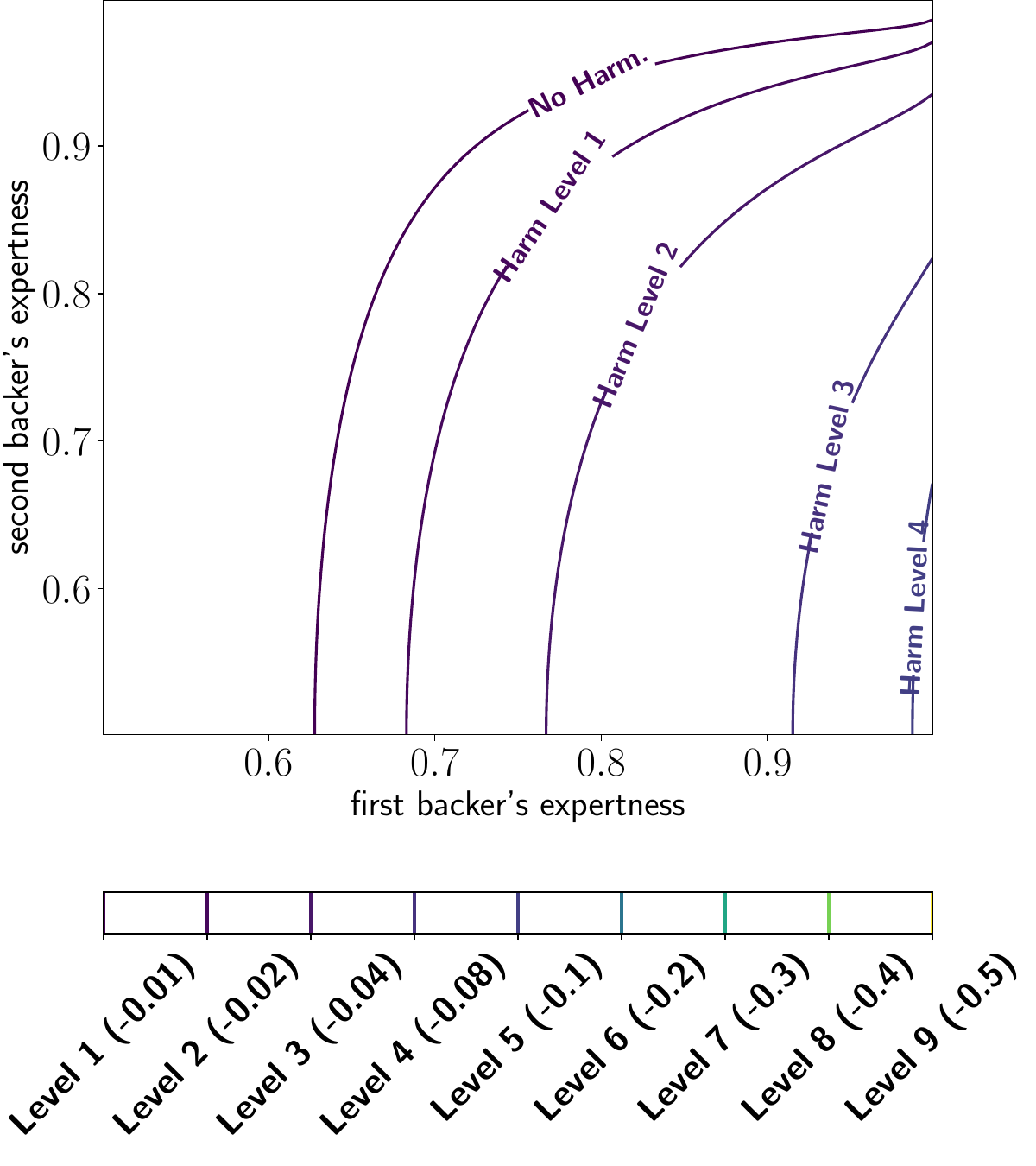} \vspace{-5pt}
            \caption{High-quality competing against high-quality} 
            \label{fig:reg_11_projects}
        \end{subfigure}
          \begin{subfigure}[b]{0.4\textwidth}
            \centering 
            \includegraphics[width=0.8\textwidth]{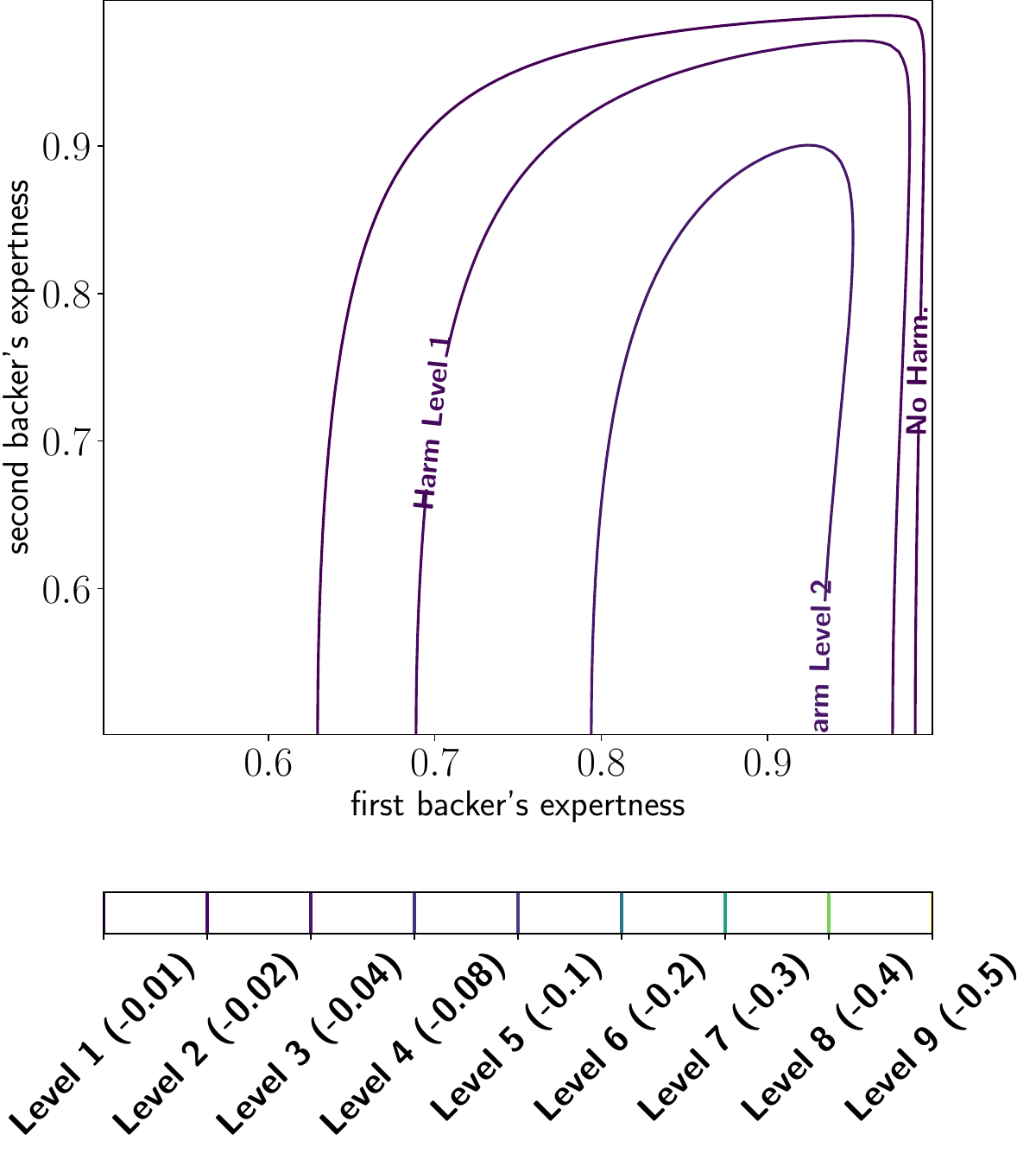} \vspace{-5pt}
            \caption{High-quality competing against low-quality} 
            \label{fig:reg_10_projects}
        \end{subfigure}
        \vspace{10pt}
        
        \begin{subfigure}[b]{0.4\textwidth}
            \centering 
            \includegraphics[width=0.8\textwidth]{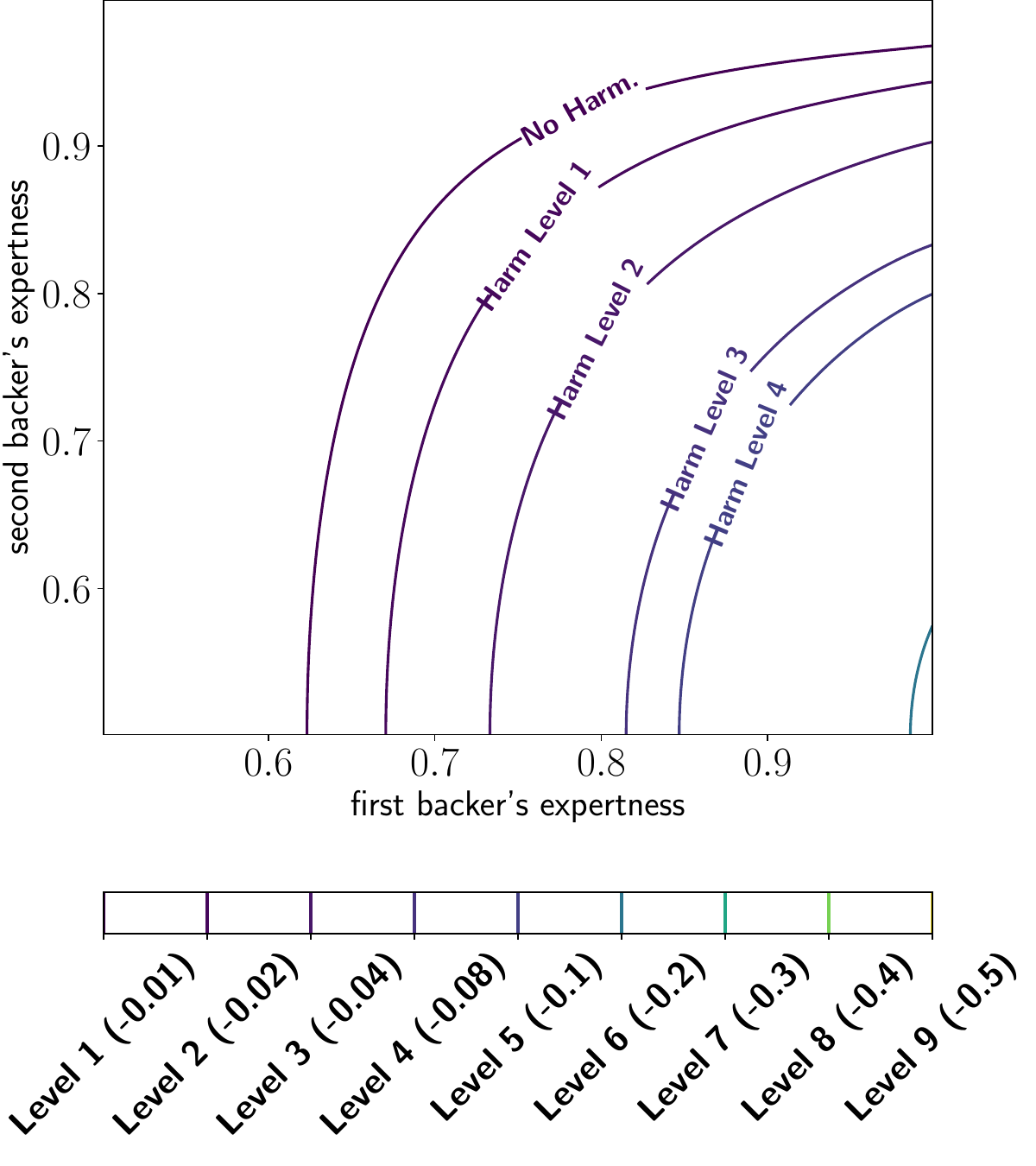} \vspace{-5pt}
            \caption{Low-quality competing against high-quality} 
            \label{fig:reg_10_projects}
        \end{subfigure}
         \begin{subfigure}[b]{0.4\textwidth}
            \centering 
            \includegraphics[width=0.8\textwidth]{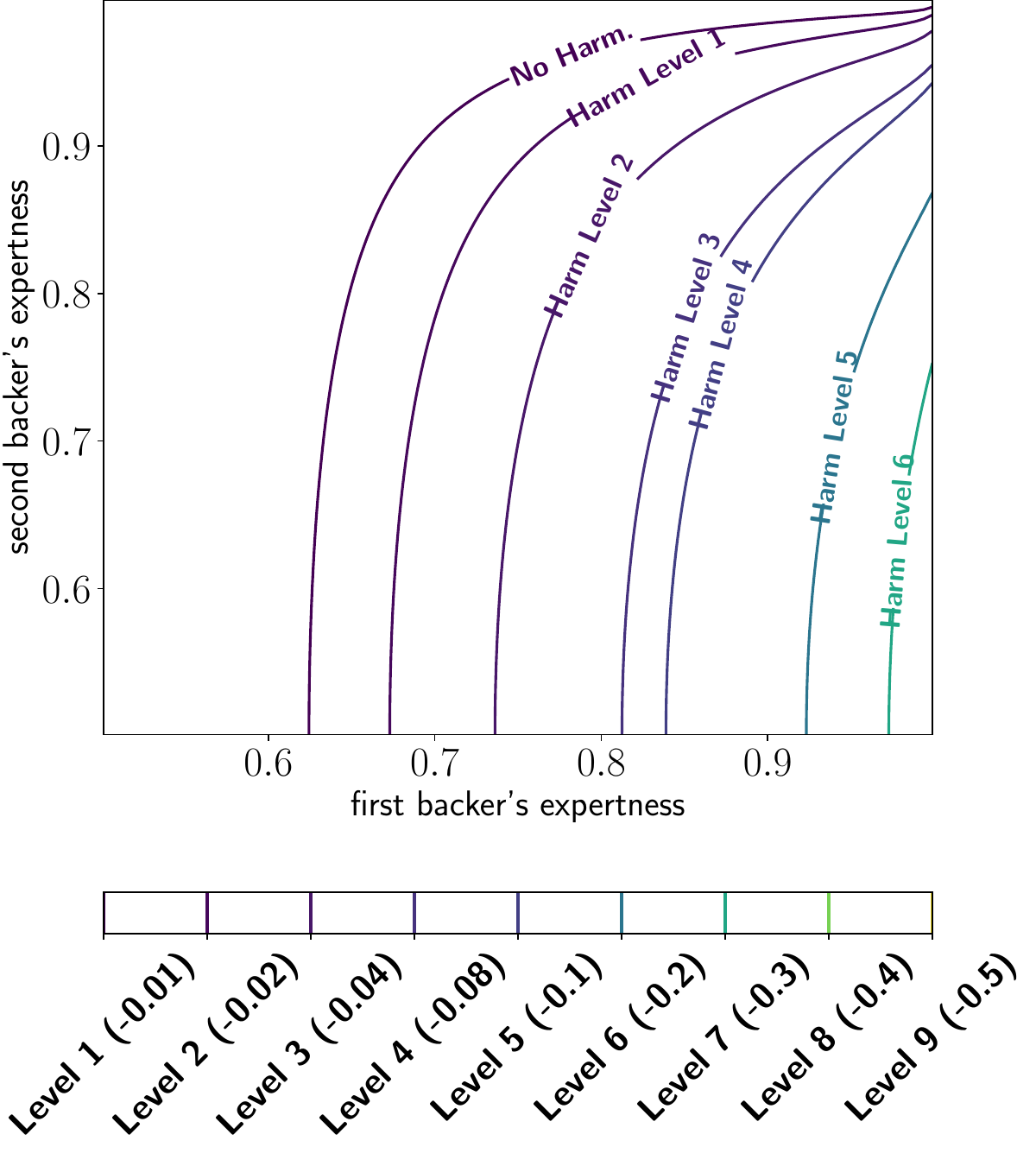} \vspace{-5pt}
            \caption{Low-quality competing against low-quality} 
            \label{fig:reg_10_projects}
        \end{subfigure}
        \caption{\small Harm Potential of OL on Project Success Probability, $\success^{OL}(\bVp) - \success^{NL}(\bVp)$, in Relaxed Competition}
    \label{fig:suc_probs_relaxed}
\end{figure}

To quantify the improvement and harm potential of OL to improve the success probability in tight and relaxed competitions, we present Figures \ref{fig:succ_projects_tight} and \ref{fig:suc_probs_relaxed}, respectively. 
 
When the early backer has a significantly higher expertness level than the late backer, the late backer tends to have higher confidence in the early backer's decision rather than relying solely on his own private signal. In such cases, if the early backer's expertness level is high enough, the late backer may follow the early backer, creating a herding effect. Our findings under tight competition, as illustrated in Figure~\ref{fig:succ_projects_tight}, indicate that the expertness conditions that result in the highest improvement potential of observational learning (OL) are those that facilitate herding. However, we also observe a notable difference between high- and low-quality projects in this regard. Specifically, for high-quality projects, OL is most beneficial when the early backer has a high expertness level, while for low-quality projects, the highest improvement is obtained when the early backer's expertness level is relatively moderate. High-quality projects benefit more from highly expert backers who can accurately identify their true quality and influence late backers in their favor. On the other hand, low-quality projects may not necessarily want backers to be overly accurate in identifying the true quality of their products. However, in tight competition, low-quality project creators can still rely on herding as a strategy to gain an advantage. Therefore, low-quality projects may be better off when the early backer's expertness level is sufficiently high to facilitate herding, but not high enough to accurately infer the low quality of their project.

In the relaxed competition case, as can be observed in  Figure \ref{fig:suc_probs_relaxed}, the expertness conditions leading to herding become unfavorable and harm the projects' success probabilities. In contrast to the tight competition case, in relaxed competition, project creators prefer the early backers to have a low expertness level, especially when their project is low-quality.

\subsection{Impact of Observational Learning on Platform's Profit}\label{sec:profit_short}

In this section, we evaluate the effect of observational learning on the platform's profit and show that it is heavily dependent on whether the competition is tight or relaxed.

\begin{proposition} \label{prop:phi}
In \textbf{tight competition}, unless both projects are low-quality, observational learning has no potential to harm the platform's profit, i.e., $\imp_{\profit^-} (\bV) = \emptyset, \;\;\; \forall \bV \in \cV \setminus \{ (0,0) \}$, and when both projects are low-quality, observational learning improves the platform's profit in the majority of expertness conditions, i.e., $| \imp_{\profit^+} (0,0) | > | \imp_{\profit^-} (0,0) |$. The maximum improvement potential and average impact of observational learning on the crowdfunding platform's profit in each quality state are as follows: \vspace{-25pt}
\begin{multicols}{2}
\begin{equation}
\dplus{\profit}{V}=\begin{cases}
			1.00, & V=(1,0), (0,1)\\
                0.05, & V=(1,1)\\
                0.04, & V= (0,0)\\
		 \end{cases}
   \nonumber
\end{equation}

\begin{equation}
\dave {\profit}{V}=\begin{cases}
			0.13, & V=(1,0), (0,1)\\
                0.08, & V=(1,1)\\
                0.01, & V= (0,0)\\
		 \end{cases} \nonumber
\end{equation}
\end{multicols}

In \textbf{relaxed competition}, observational learning harms the platform's profit in the majority of expertness conditions across all quality states, i.e., $| \imp_{\profit^-} (\bV) | > | \imp_{\profit^+} (\bV) |, \;\;\; \bV \in \cV$. The maximum harm potential and average impact of observational learning on the platform's profit are given below for each quality state. \vspace{-25pt}
\begin{multicols}{2}
\begin{equation}
\dminus{\profit}{V}=\begin{cases}
			-0.06, & V=(1,0), (0,1)\\
                -0.10, & V=(1,1)\\
                -0.32, & V= (0,0)\\
		 \end{cases}
   \nonumber
\end{equation}

\begin{equation}
\dave {\profit}{V}=\begin{cases}
			-0.002, & V=(1,0), (0,1)\\
                -0.03, & V=(1,1)\\
                -0.07, & V= (0,0)\\
		 \end{cases} \nonumber
\end{equation}
\end{multicols}

\end{proposition}

Since the platform collects fees only from successful projects, and in a tight competition only one project has a chance to succeed, herding helps to accumulate the available funding in the favor of one of the projects and increases the platform's profit. Hence, in tight competition, OL improves the platform's profit on average, even when both projects are low-quality. This benefit of herding is particularly apparent when a high-quality project competes against a low-quality project. On the other hand, under relaxed competition, herding due to OL decreases the platform's profit, since it decreases the success probability of one of the projects.

\begin{figure}[h]
    \centering
         \begin{subfigure}[b]{0.32\textwidth}
            \centering
            \includegraphics[width=\textwidth]{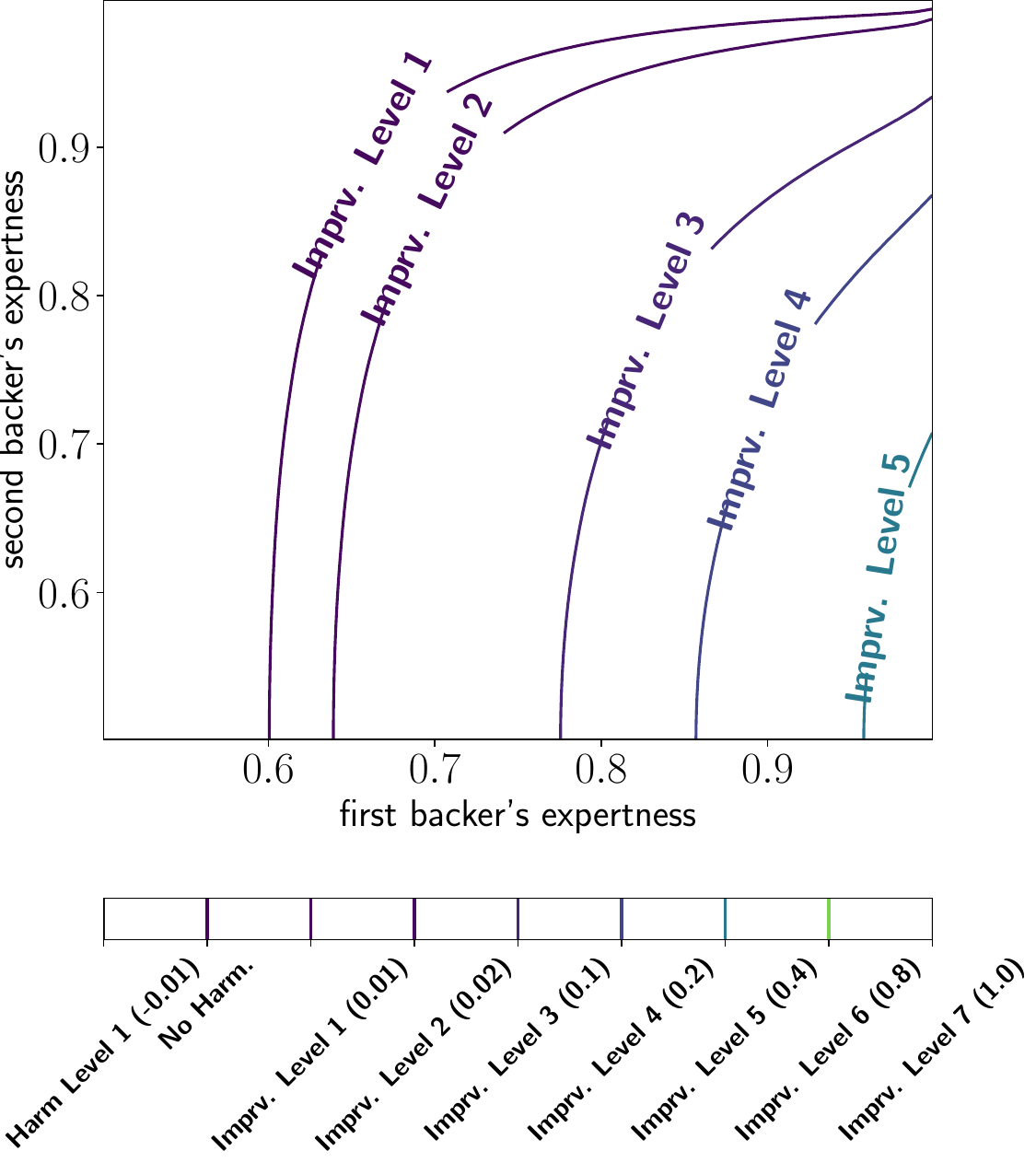}
            \caption{Both projects are high-quality.}    
            \label{fig:reg_11_projects}
        \end{subfigure}
         \begin{subfigure}[b]{0.32\textwidth}
            \centering 
            \includegraphics[width=\textwidth]{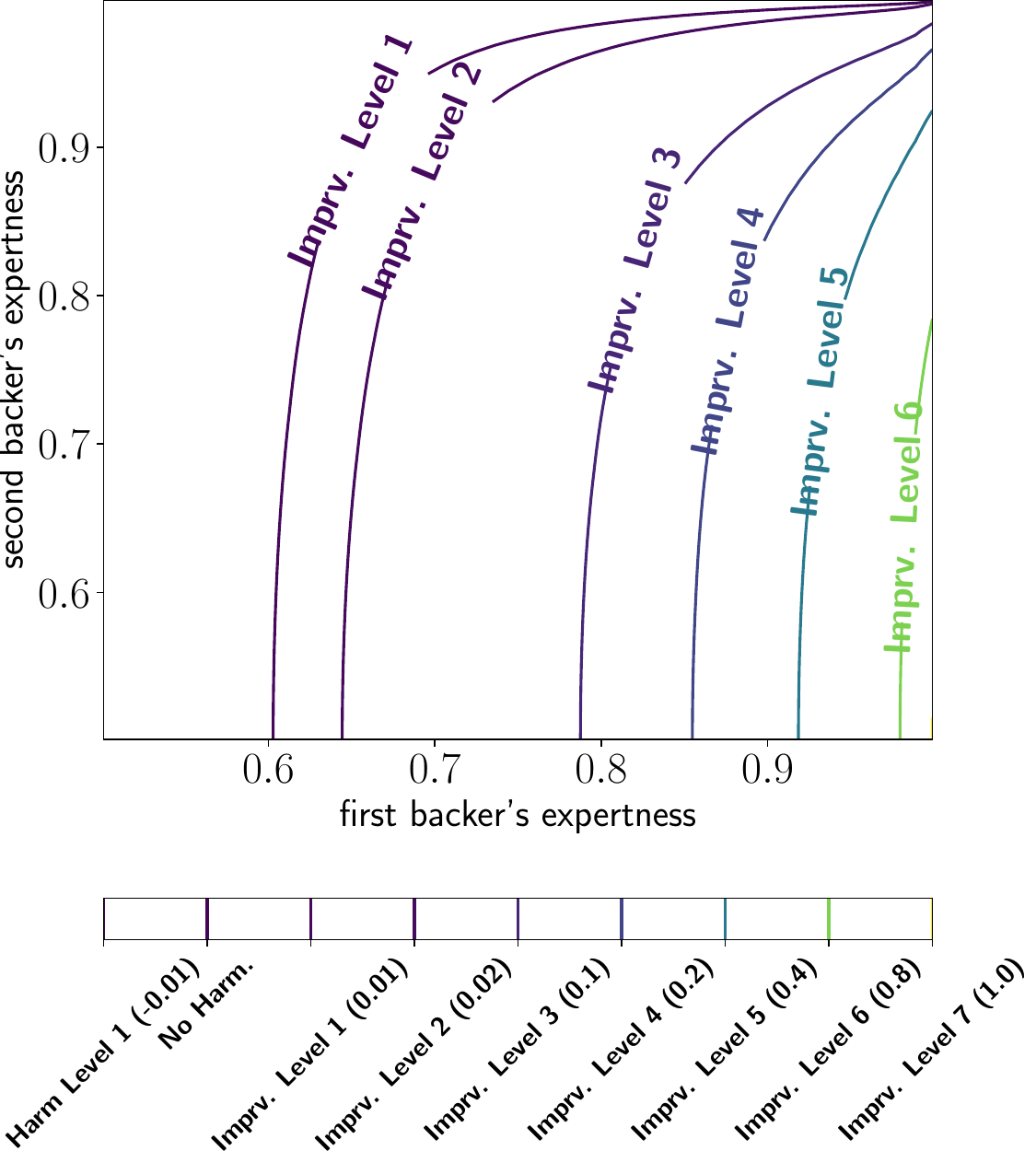}
            \caption{Only one project is high-quality.}    
            \label{fig:reg_10_projects}
        \end{subfigure}
         \begin{subfigure}[b]{0.32\textwidth}
            \centering 
            \includegraphics[width=\textwidth]{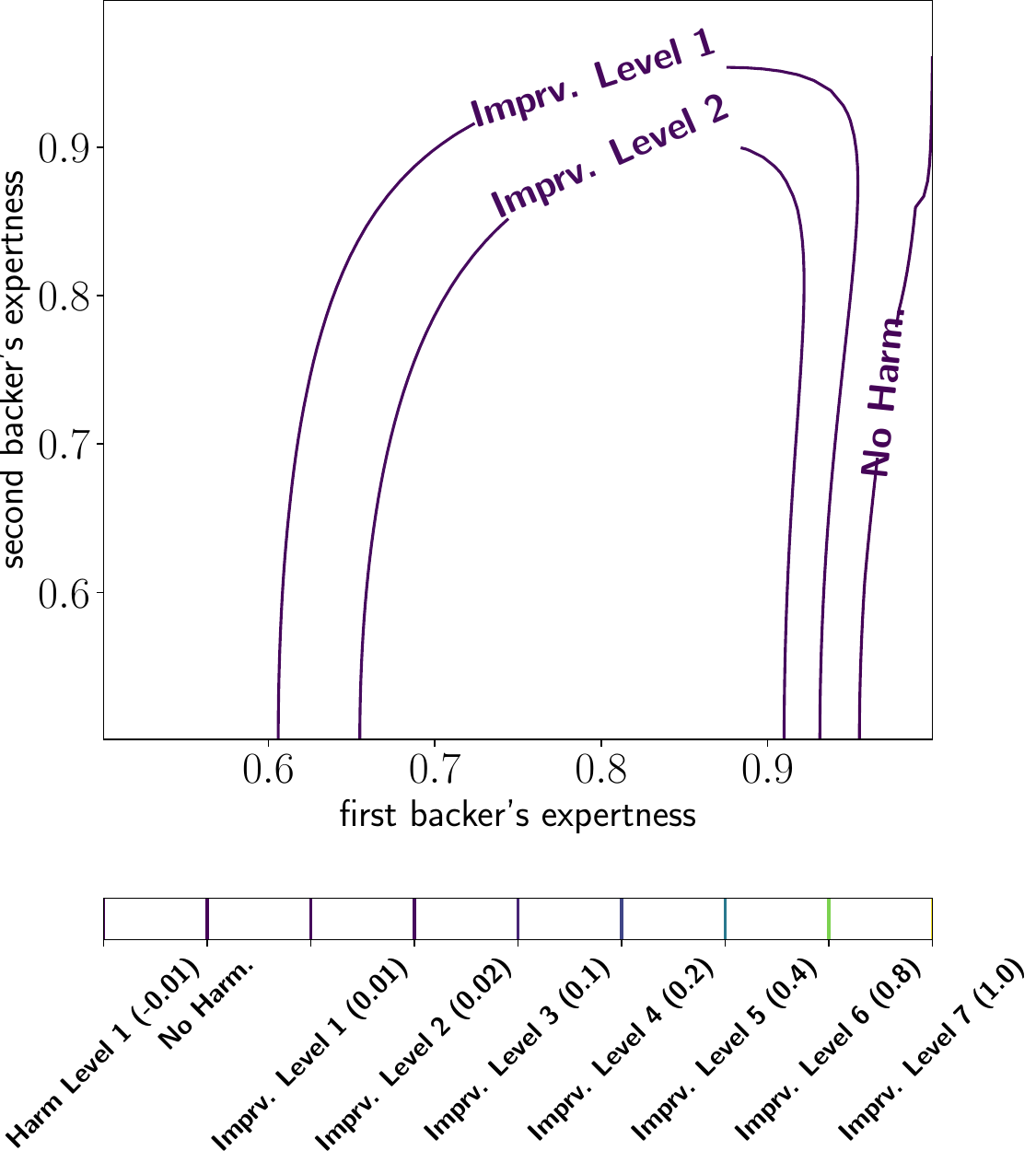}
            \caption{Both projects are low-quality.}    
            \label{fig:reg_10_projects}
        \end{subfigure}
        \caption{Effect of OL on the Platform's Profit, $\profit^{OL} (\bVp) - \profit^{NL} (\bVp)$, in Tight Competition.}
    \label{fig:platform_profits_tight}
\end{figure}

\begin{figure}[h]
    \centering
         \begin{subfigure}[b]{0.32\textwidth}
            \centering
            \includegraphics[width=\textwidth]{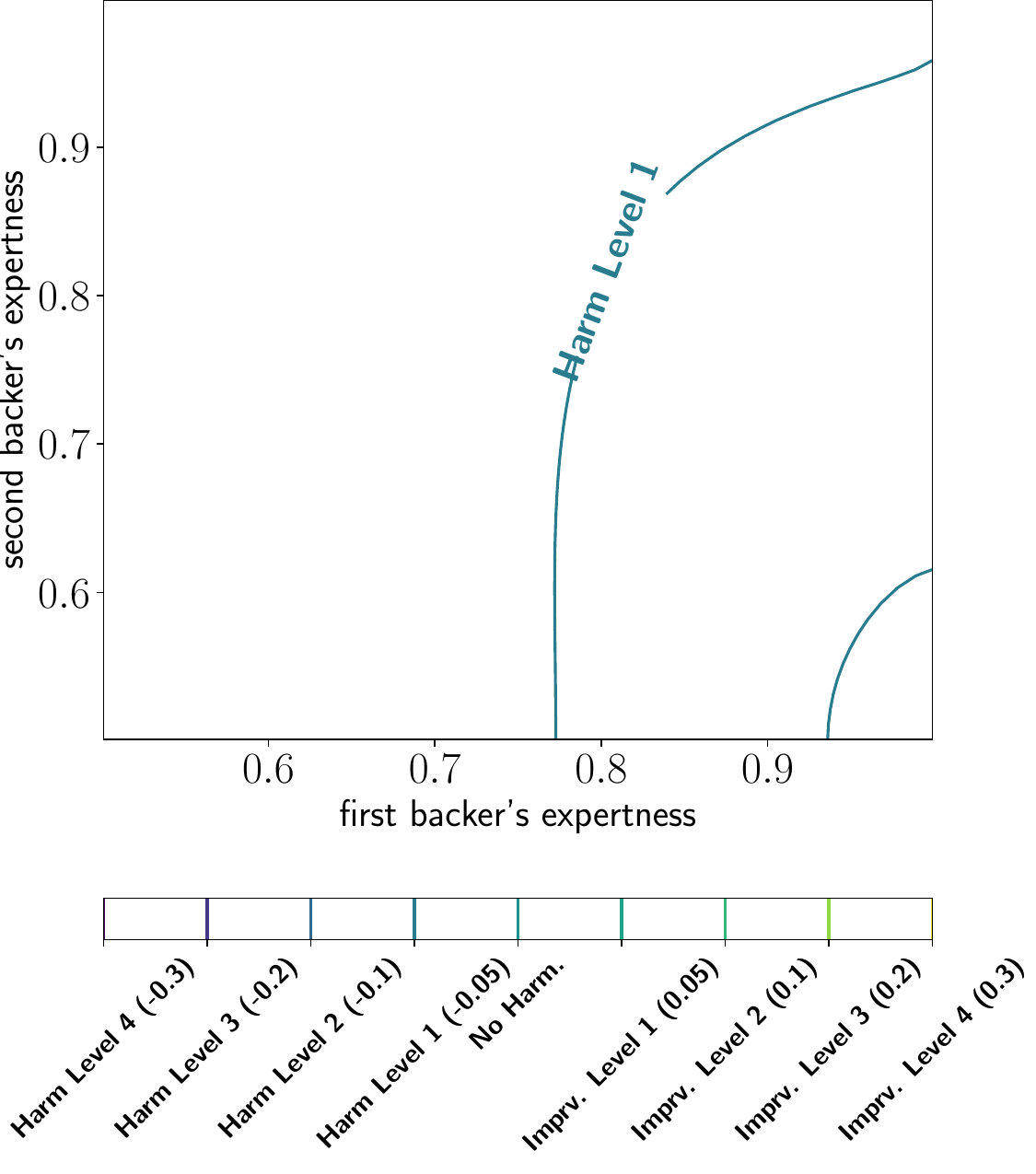}
            \caption{Both projects are high-quality.}
            \label{fig:reg_11_projects}
        \end{subfigure}
          \begin{subfigure}[b]{0.32\textwidth}
            \centering 
            \includegraphics[width=\textwidth]{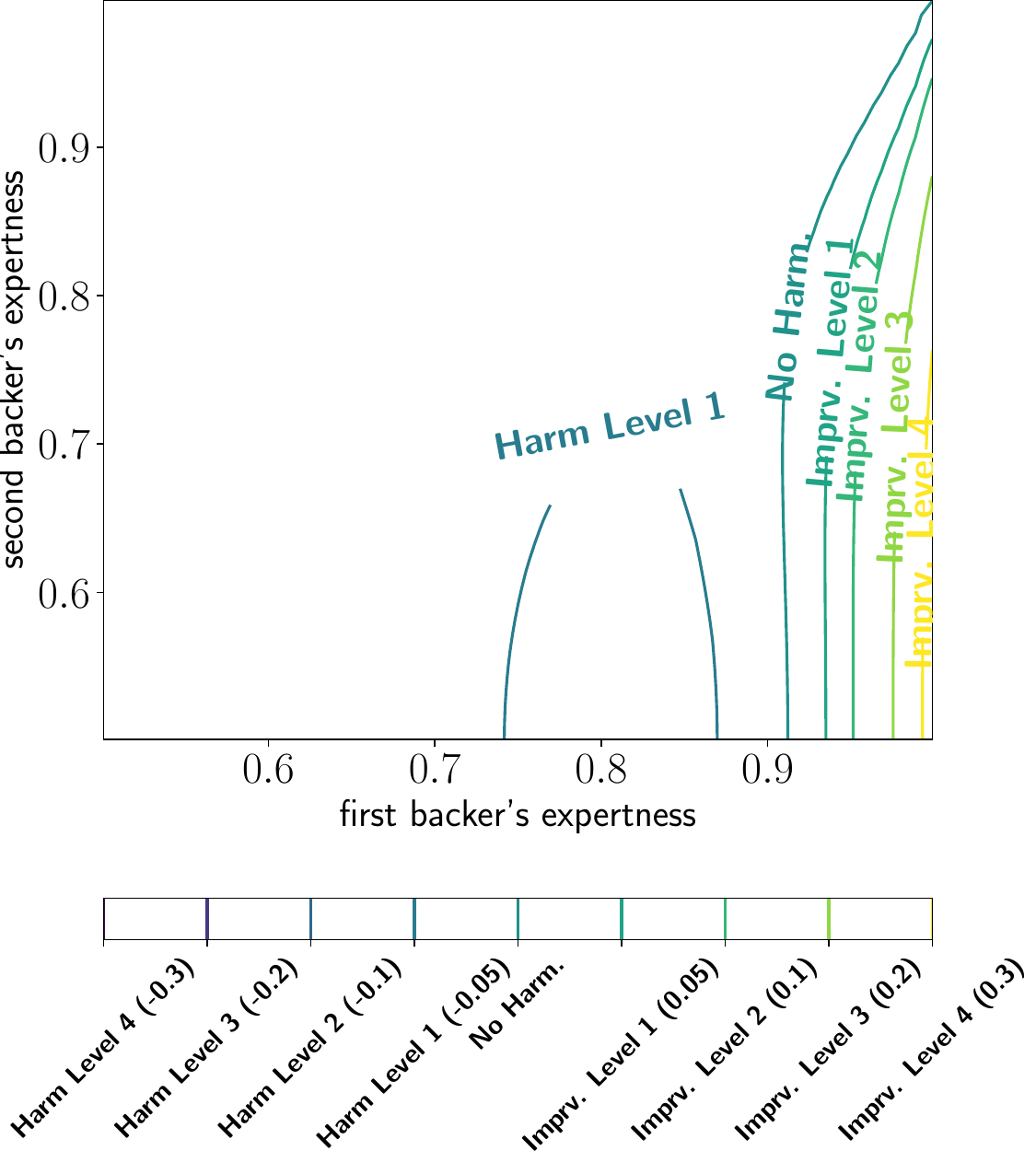}
            \caption{Only one project is high-quality.}
            \label{fig:reg_10_projects}
        \end{subfigure}
         \begin{subfigure}[b]{0.32\textwidth}
            \centering 
            \includegraphics[width=\textwidth]{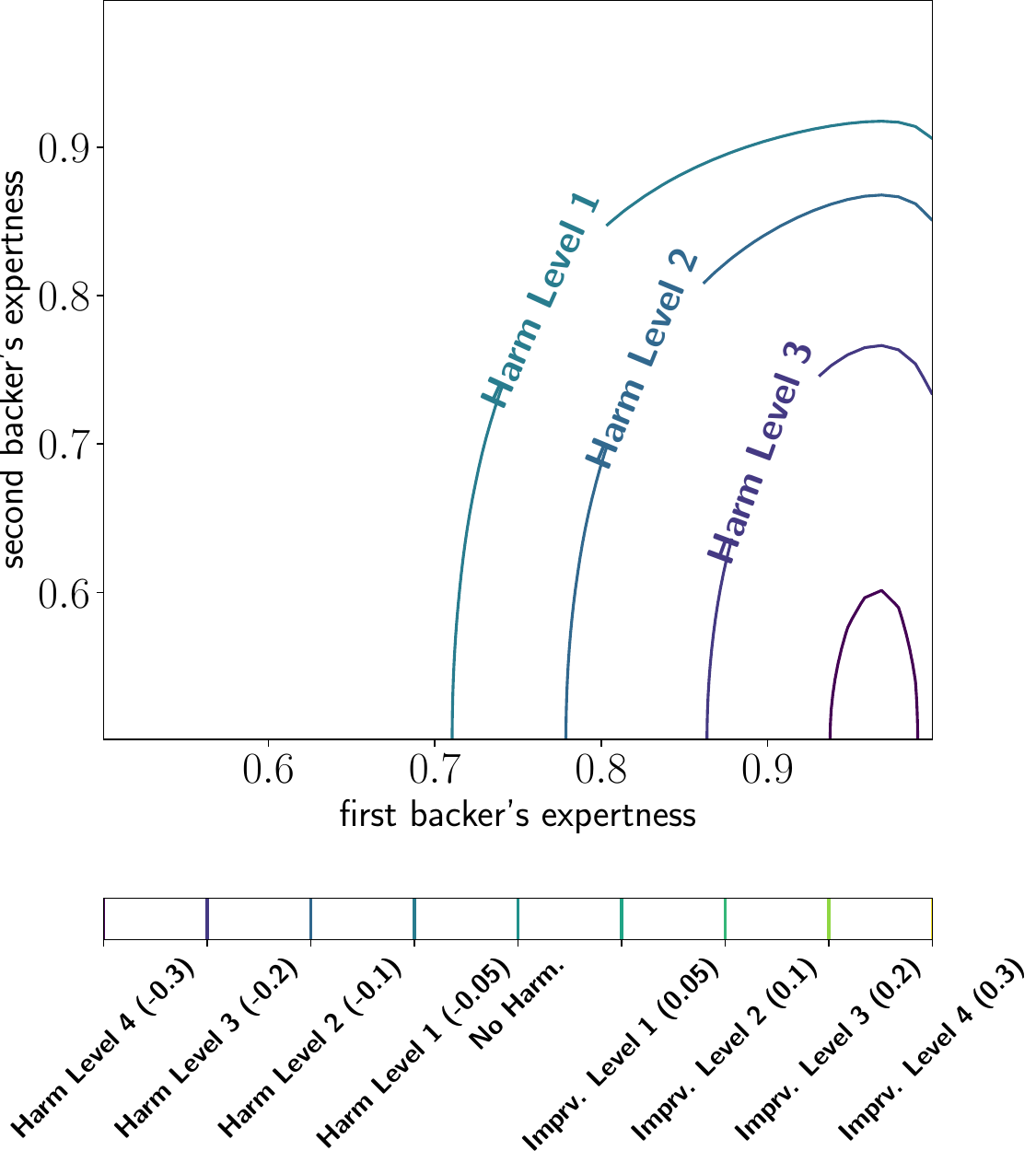}
            \caption{Both projects are low-quality.}    
            \label{fig:reg_10_projects}
        \end{subfigure}
        \caption{Effect of OL on the Platform's Profit, $\profit^{OL} (\bVp) - \profit^{NL} (\bVp)$, in Relaxed Competition}
    \label{fig:platform_profits_relaxed}
\end{figure}

To illustrate the impact of OL on the platform's profit under all expertness conditions, we present Figures \ref{fig:platform_profits_tight} and \ref{fig:platform_profits_relaxed} for the tight and relaxed competition cases, respectively. When Figures \ref{fig:platform_profits_tight}--\ref{fig:platform_profits_relaxed} are examined in parallel with Figures \ref{fig:succ_projects_tight}--\ref{fig:suc_probs_relaxed}, regarding the projects' success probabilities, the relation between the two performance metrics can be identified. Particularly in the tight competition case, we observe that the expertness conditions where OL improves the platform's profit the most coincide with those where it improves the projects' success chances. The only exception to this observation arises when both projects are low-quality. OL does not harm the success probability of low-quality projects under any expertness conditions, but it may have a negative impact on platform profit when the early backer has a high expertness level. In such a scenario, herding can occur in the favor of not pledging to either of the two projects, as the early backer may accurately identify that both projects are low-quality, leading to a reduction in total funding realized, and therefore in the platform's profit. 

Upon examining the impact of OL on the platform's profit in the case of relaxed competition (Figure \ref{fig:platform_profits_relaxed}), we observe that the relation between success probabilities and the platform's profit is prominent when both projects are of high or low quality. In these cases, the expertness conditions that are favorable for herding result in OL harming both the success probabilities and the platform's profit. On the other hand, when one project is high-quality and the other is low-quality, the impact of OL on the platform's profit becomes more complex. We find that when the two competing projects have different quality states, the expertness conditions that facilitate herding (i.e., when the early backer has a very high expertness level) allow OL to improve the platform's profit. However, under expertness conditions where the early backer has an expertness level that is sufficient to influence the late backer but not enough to accurately identify the projects' true quality states, OL has a potential to harm the platform's profit. The positive impact observed under the first herding scenario (i.e., when the early backer is very highly expert) can be attributed to increased prospects for the high-quality project, while the negative impact observed under the second herding scenario (i.e., when the early backer is moderately expert) can be explained by herding in favor of not backing either of the projects.

 \vspace{-15pt}

Our findings on how OL affects the platform's profit indicate that the platform should allow OL to benefit from herding when the available funding is scarce, especially when there is a mix of high- and low-quality projects on the platform. On the other hand, when funding is plentiful, OL can harm the platform profits, particularly when projects are similar in quality. 

\subsection{Impact of Observational Learning on Crowdfunding Effectiveness} \label{sec:profit_long}

The effectiveness of a crowdfunding platform increases with increasing funds collected by high-quality projects and decreases with increasing funds collected by low-quality projects. The analysis of this particular performance metric highlights two cases where the average impact of OL will be negative.

\begin{proposition} \label{prop:eff}
In \textbf{tight competition}, unless both projects are low-quality, observational learning has no potential to harm platform effectiveness, 
i.e., $\imp_{\effect^-} (\bV) = \emptyset, \;\;\; \forall \bV \in \cV \setminus \{(0,0)\}$, however, when both projects are low-quality, observational learning harms platform effectiveness in the majority of expertness conditions, i.e., $| \imp_{\effect^-} (0,0)| > | \imp_{\effect^+} (0,0) |$. The maximum improvement and harm potentials, and average impact of observational learning on platform effectiveness in all four quality states are as follows: \vspace{-20pt}
\footnotesize
\begin{multicols}{3}
\begin{equation}
\dplus{\effect}{V}=\begin{cases}
			1.00, & V=(1,0), (0,1)\\
                0.50, & V=(1,1)\\
                0.005, & V= (0,0)\\
		 \end{cases}
   \nonumber
\end{equation}

\begin{equation}
\dminus {\effect}{V}=\begin{cases}
			0.00, & V=(1,0), (0,1),\\
   &\qquad(1.1)\\
                -0.01, & V= (0,0)\\
		 \end{cases} \nonumber
\end{equation}

\begin{equation}
\dave {\effect}{V}=\begin{cases}
			0.12, & V=(1,0), (0,1)\\
                0.08, & V=(1,1)\\
                -0.01, & V= (0,0)\\
		 \end{cases} \nonumber
\end{equation}
\end{multicols}
\normalsize

In \textbf{relaxed competition}, when both projects are low-quality, observational learning has no potential to harm platform effectiveness, i.e., $\imp_{\effect^-} (0,0) = \emptyset$, and when only one project is high-quality, it improves platform effectiveness in the majority of expertness conditions, i.e., $| \imp_{\effect^+} (\bV) | > | \imp_{\effect^-} (\bV) |, \;\;\; \bV \in \{(0,1), (1,0)\}$. When both projects are high-quality, observational learning harms crowdfunding effectiveness in the majority of expertness conditions, i.e., $| \imp_{\effect^-} (1,1) | > | \imp_{\effect^+} (1,1) |$. The maximum improvement and harm potentials, and average impact of observational learning on platform effectiveness across all four quality states are given below. \vspace{-25pt}
\footnotesize
\begin{multicols}{3}
\begin{equation}
\dplus{\effect}{V}=\begin{cases}
			0.56, & V=(1,0), (0,1)\\
                0.43, & V=(0.0)\\
                0.00, & V= (1,1)\\
		 \end{cases}
   \nonumber
\end{equation}

\begin{equation}
\dminus {\effect}{V}=\begin{cases}
			0.00, & V=(0,0)\\
   			-0.03, & V=(1,0), (0,1)\\
			-0.10, & V=(1,1)\\
		 \end{cases} \nonumber
\end{equation}

\begin{equation}
\dave {\effect}{V}=\begin{cases}
			0.07, & V=(0,0)\\
                0.04, & V=(1,0), (0,1)\\
                -0.03, & V= (1,1)\\
		 \end{cases} \nonumber
\end{equation}
\end{multicols}
\normalsize

 \end{proposition}

\begin{figure}[h]
    \centering
         \begin{subfigure}[b]{0.32\textwidth}
            \centering
            \includegraphics[width=\textwidth]{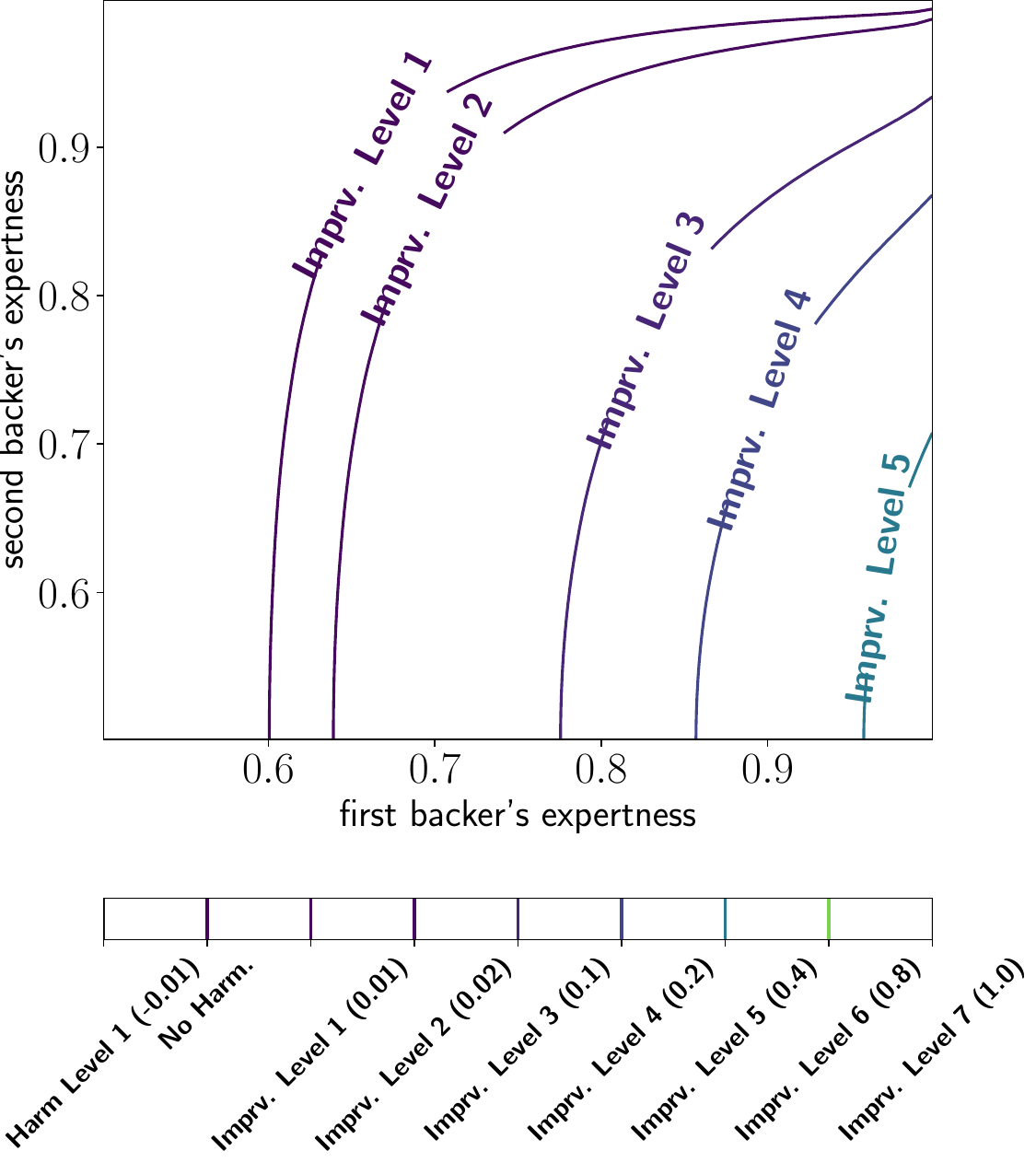}
            \caption[Network2]{Both projects are high-quality.}    
            \label{fig:reg_11_projects}
        \end{subfigure}
         \begin{subfigure}[b]{0.32\textwidth}
            \centering 
            \includegraphics[width=\textwidth]{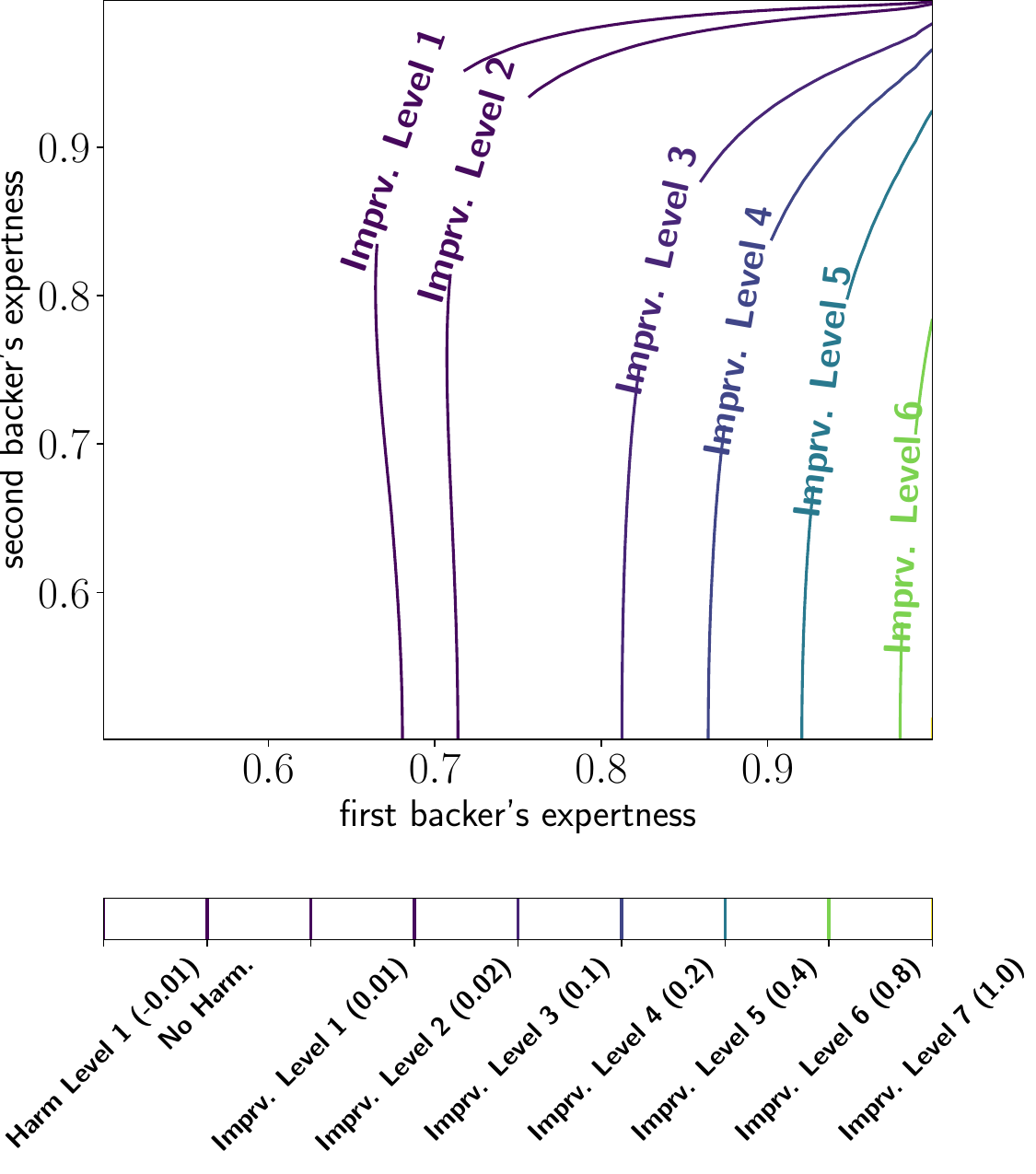}
            \caption{Only one project is high-quality.}  
            \label{fig:reg_10_projects}
        \end{subfigure}
         \begin{subfigure}[b]{0.32\textwidth}
            \centering 
            \includegraphics[width=\textwidth]{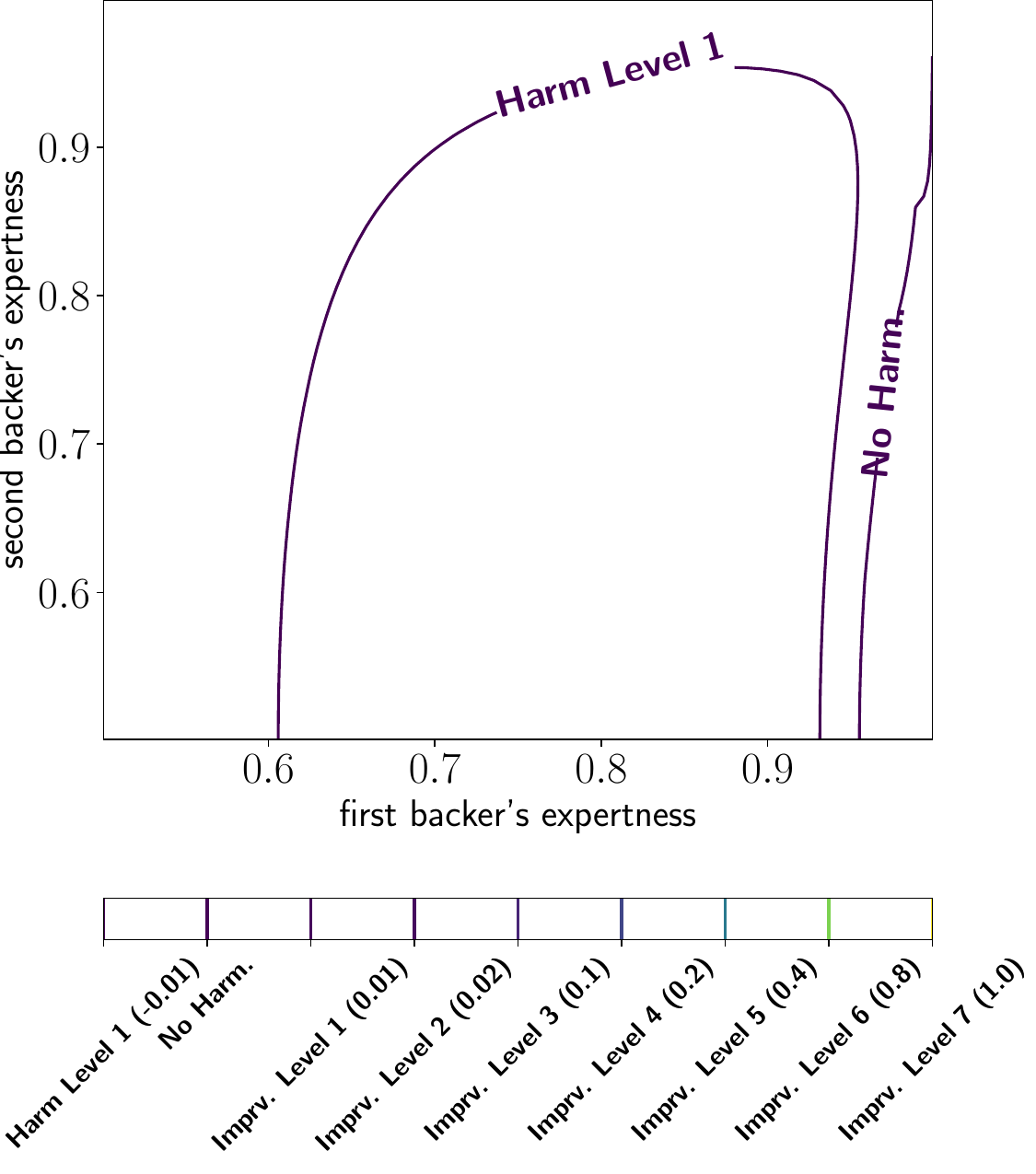}
            \caption{Both projects are low-quality.}
            \label{fig:reg_10_projects}
        \end{subfigure}
        \caption{Effect of OL on Effectiveness, $\effect^{OL} (\bVp) - \effect^{NL} (\bVp)$, in Tight Competition}
    \label{fig:effectiveness_tight}
\end{figure}

\begin{figure}[h]
    \centering
         \begin{subfigure}[b]{0.32\textwidth}
            \centering
            \includegraphics[width=\textwidth]{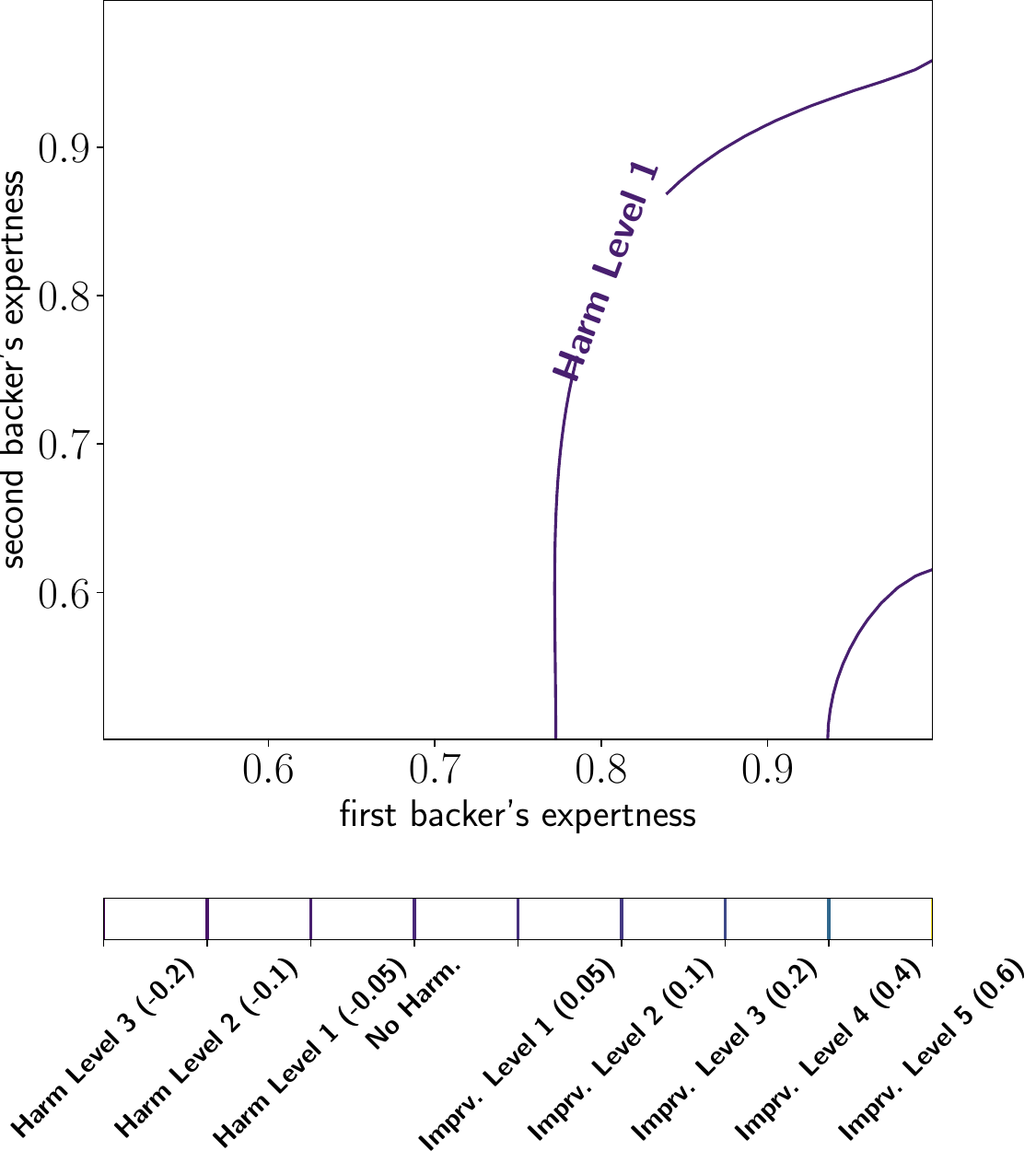}
            \caption[Network2]{Both projects are high-quality.}  
            \label{fig:reg_11_projects}
        \end{subfigure}
         \begin{subfigure}[b]{0.32\textwidth}
            \centering 
            \includegraphics[width=\textwidth]{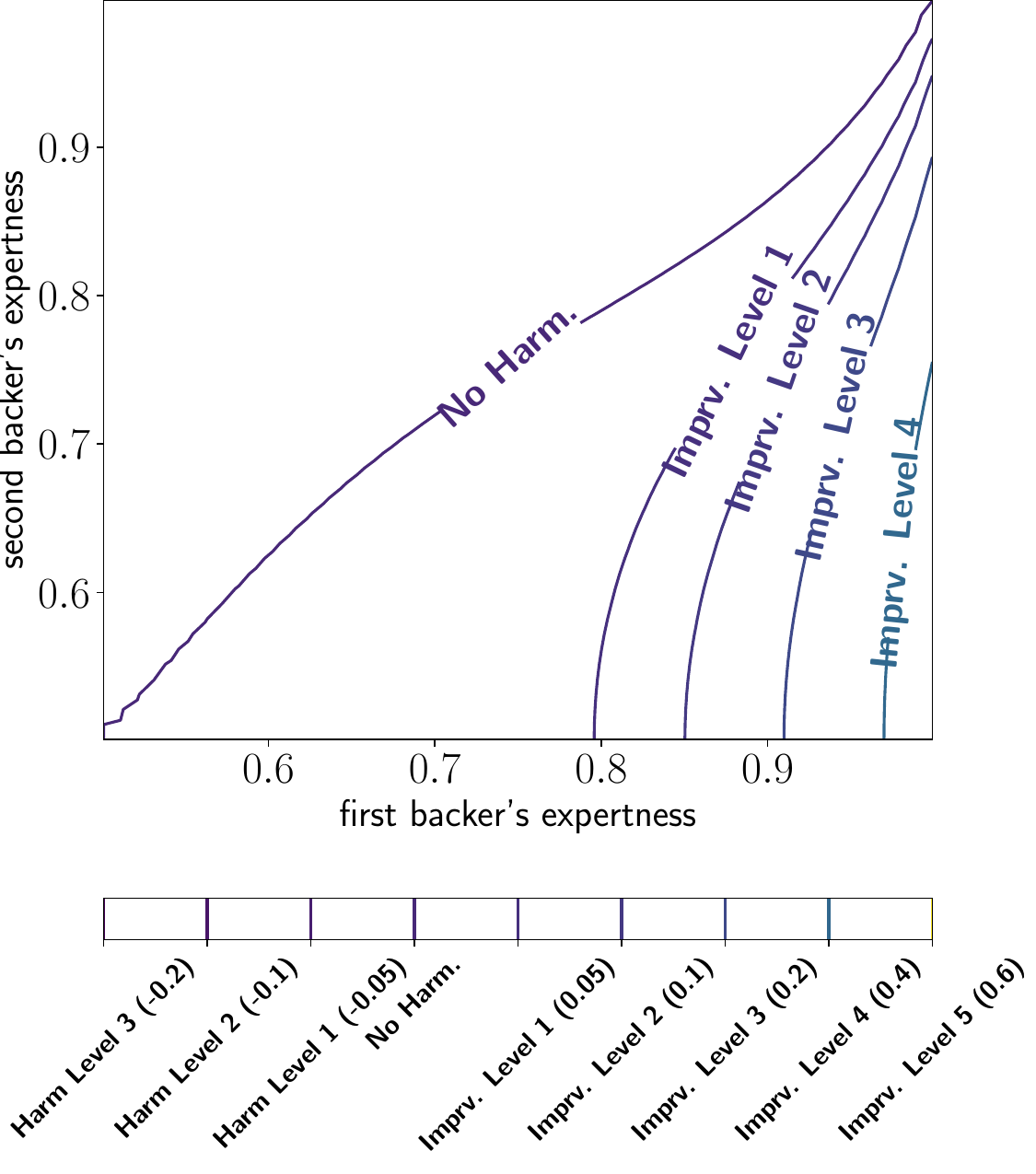}
            \caption{Only one project is high-quality.}  
            \label{fig:reg_10_projects}
        \end{subfigure}
         \begin{subfigure}[b]{0.32\textwidth}
            \centering 
            \includegraphics[width=\textwidth]{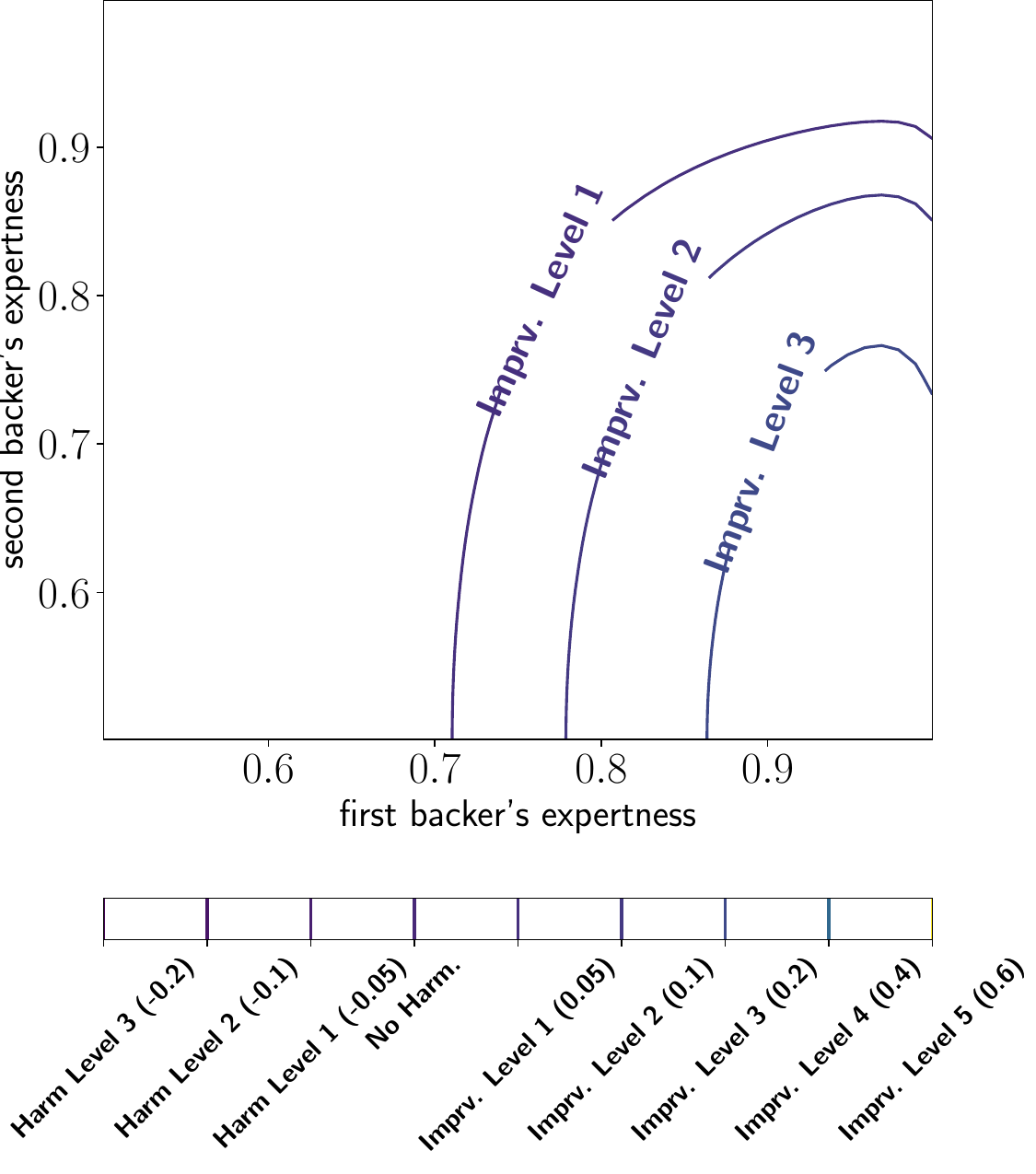}
            \caption{Both projects are low-quality.}  
            \label{fig:reg_10_projects}
        \end{subfigure}
        \caption{Effect of OL on Effectiveness, $\effect^{OL} (\bVp) - \effect^{NL} (\bVp)$, in Relaxed Competition}
    \label{fig:effectiveness_relaxed}
\end{figure}

To illustrate the impact of OL on effectiveness under all expertness conditions, we present Figures \ref{fig:effectiveness_tight} and \ref{fig:effectiveness_relaxed} for the tight and relaxed competition cases, respectively. 

The findings in Proposition \ref{prop:eff} indicate that observational learning improves the platform effectiveness in a wide range of backer expertness, with the exception of two cases: \begin{enumerate*}[label = \itshape(\roman*)] \item when two high-quality projects are in a relaxed competition, and \item when two low-quality projects are in a tight competition.\end{enumerate*} In the former case, OL creates a herding behavior and leads to one of the projects collecting more funds than it requires, unnecessarily decreasing the success chance of the other. In the latter case, the herding behavior induced by OL may result in a low-quality project being funded, which is not a desirable outcome for platform effectiveness. This result supports what has also been uncovered in Section \ref{sec:profit_short}: OL is particularly preferable for the crowdfunding platforms when there is a mix of high- and low-quality projects competing.

\section{Impact of Observational Learning on Product Development Decisions} \label{sec:decs}

So far, we have analyzed how OL affects crowdfunding dynamics for a given project quality state. In this section, we shift our attention to the creators' decision regarding the project quality, which takes place before the launch of the crowdfunding campaign. Our main goal is to understand the additional incentive OL creates for developing high-quality projects.

Developing a high-quality product is costlier than a low-quality one \citep{Chakraborty2021}. In this section, we assume that developing a high-quality product, as opposed to a low-quality product, incurs an additional cost, $c \in [0,1]$. We refer to this as the \textit{development cost}. Let the success probability of a project $i$, $\success_i$, define its \textit{profitability}. Then, the \textit{net profitability} of a low-quality project $i$ will be $\success_i$, and that of a high-quality project will be $\success_i - c$.

To capture the full potential of OL, we assume that project creators have the related information the net profitability function, market size (i.e., the number of backers, $N$) and the expertness levels of all backers before the launch of the project. Since each creator has two quality strategies, high ($V^i = 1$) or low ($V^i = 0$), the equilibrium quality strategy, denoted as $\bV^{*} (\bp)$, for any given expertness scenario, $\bp$, is as shown in equation \eqref{Eq:equilibrium}.
\begin{equation}
\label{Eq:equilibrium}
   \bV^{*} (\bp)= \begin{cases} 
      (1,1), & \success_1 \big( (1,1), \bp \big) - c > \success_1 \big( (0,1), \bp \big) , \;\; \success_2 \big( (1,1), \bp \big) - c > \success_2 \big( (1,0), \bp \big) \\
      (1,0), &  \success_1 \big( (1,0), \bp \big) - c > \success_1 \big( (0,0), \bp \big) ,\;\; \success_2 \big( (1,0), \bp \big) \geq \success_2 \big( (1,1), \bp \big) - c \\
      (0,1), &  \success_1 \big( (0,1), \bp \big) \geq \success_1 \big( (1,1), \bp \big) - c ,\;\; \success_2 \big( (0,1), \bp \big) - c > \success_2 \big( (0,0), \bp \big) \\
      (0,0), & \text{otherwise} 
   \end{cases}
\end{equation}

We analyze the impact of OL on incentivizing high-quality projects under two different expertness schemes. In Section \ref{sec:decs_backer}, we assume, as we did throughout the paper, that expertness differs per backer, but a backer's expertness level for inferring the true project quality is the same for both projects (\textit{backer-asymmetric expertness}). In Section \ref{sec:decs_project}, we conduct our analysis with a different assumption that expertness levels differ per project, but all backers have the same expertness level for the same project (\textit{project-asymmetric expertness}). By examining these two expertness schemes separately, we provide a better understanding for the impact of OL on product quality decisions under information asymmetry.

\subsection{Product Quality Decisions under Backer-Asymmetric Expertness} \label{sec:decs_backer}
In this section, we assume that the expertness levels can differ among backers, but a backer's expertness for both projects is the same. Recall that in a two-backer system, $\bp = (p_1, p_2)$ denotes the expertness conditions, where $p_1$ and $p_2$ denote the expertness levels of the first (early) and the second (late) backer, respectively.

Under backer-asymmetric expertness, the backers' expertness levels do not differ per project, and therefore the net profitability functions are symmetric for the two projects. As a result, the equilibrium strategy that is determined by equation \eqref{Eq:equilibrium} will be either $(1,1)$ or $(0,0)$. We refer to the equilibrium quality states $(1,1)$ and $(0,0)$ as the \textit{high-quality equilibrium} and the \textit{low-quality equilibrium}, respectively.

\begin{figure}[h]
    \centering
     \begin{subfigure}[b]{0.4\textwidth}
            \centering
            \includegraphics[width=0.8\textwidth]{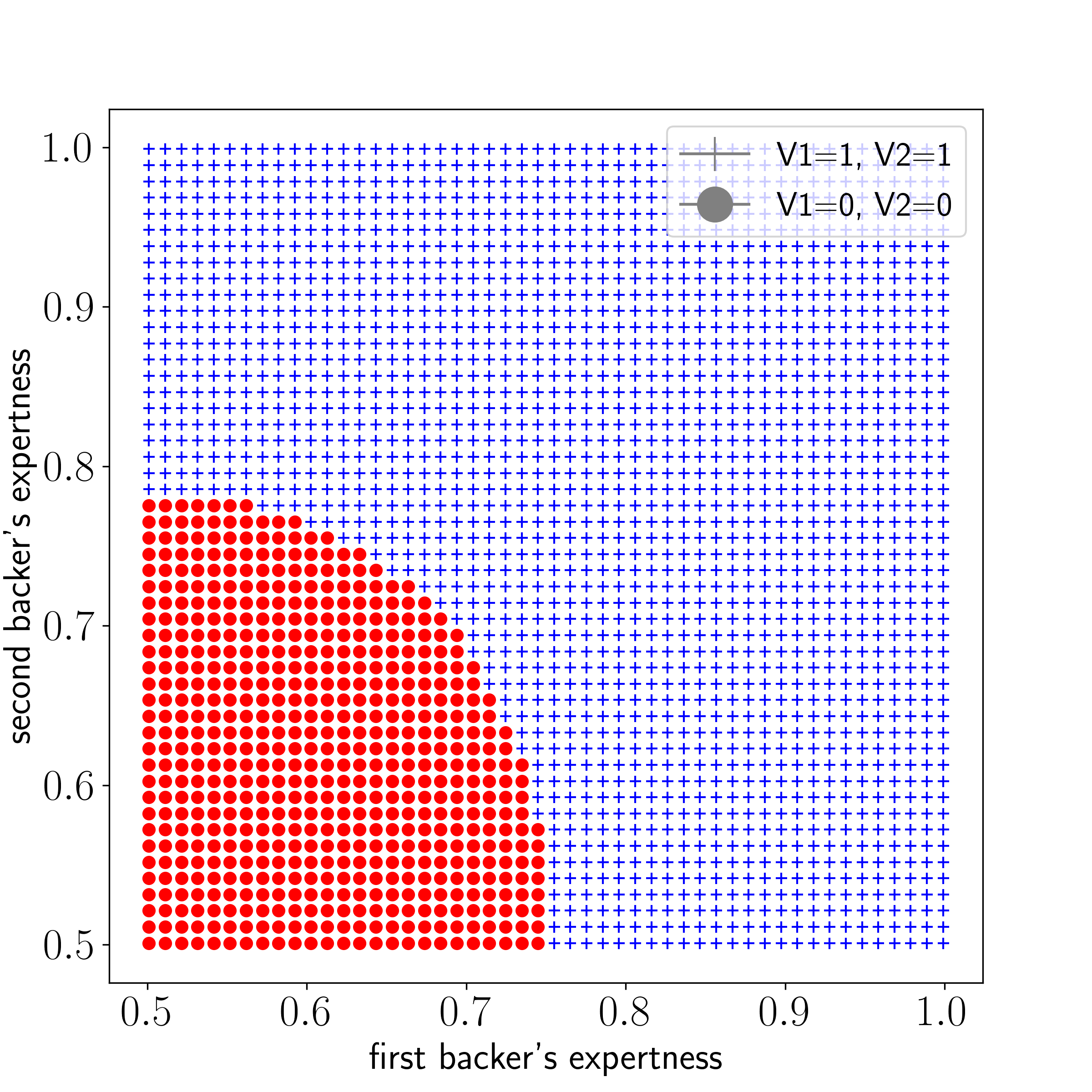}
            \caption{Tight competition, $c=0.1$}   
            \label{fig:tight_c01}
        \end{subfigure}
         \begin{subfigure}[b]{0.4\textwidth}
            \centering
            \includegraphics[width=0.8\textwidth]{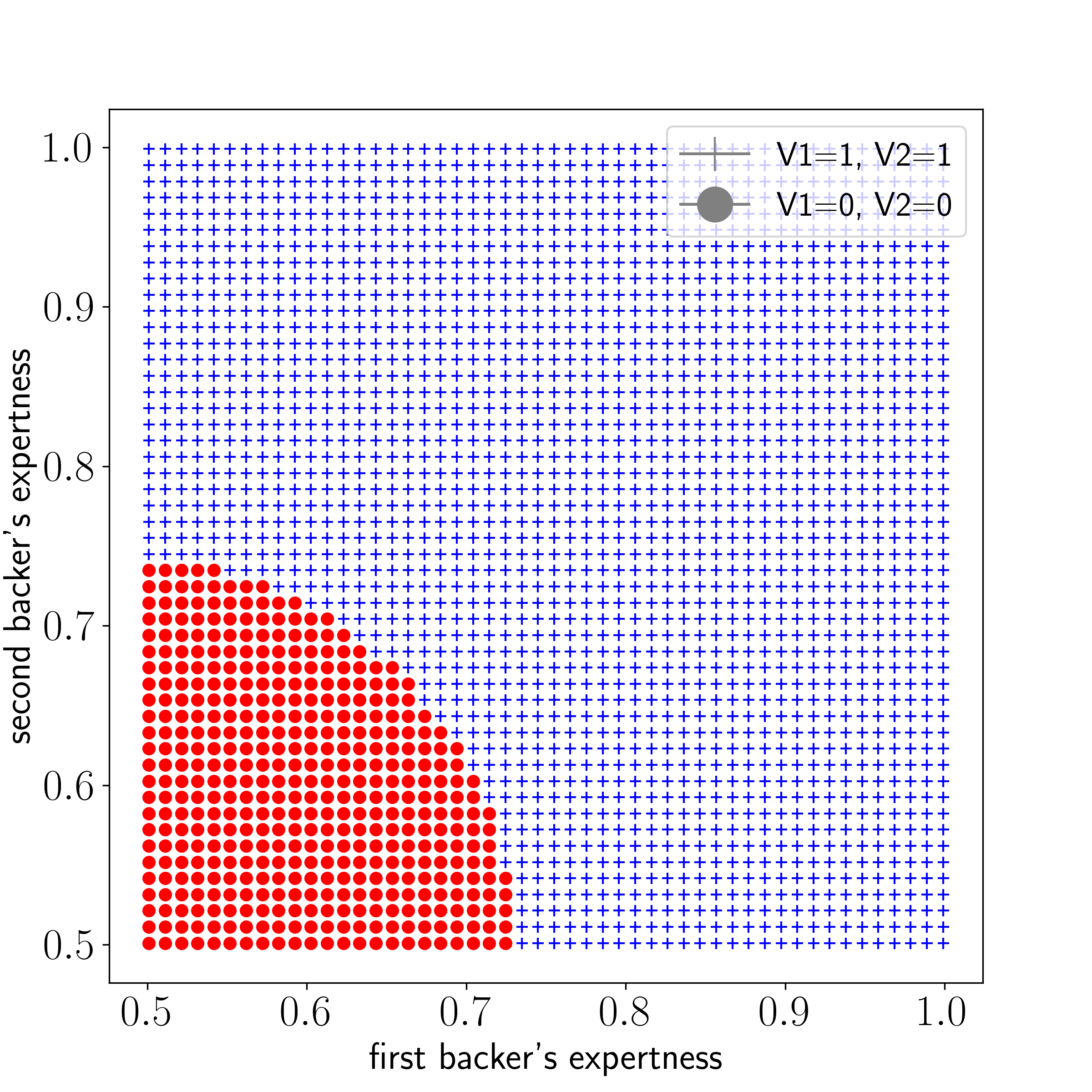}
            \caption{Relaxed competition, $c=0.1$}   
            \label{fig:rel_c01}
            \end{subfigure}
        \begin{subfigure}[b]{0.4\textwidth}
            \centering
            \includegraphics[width=0.8\textwidth]{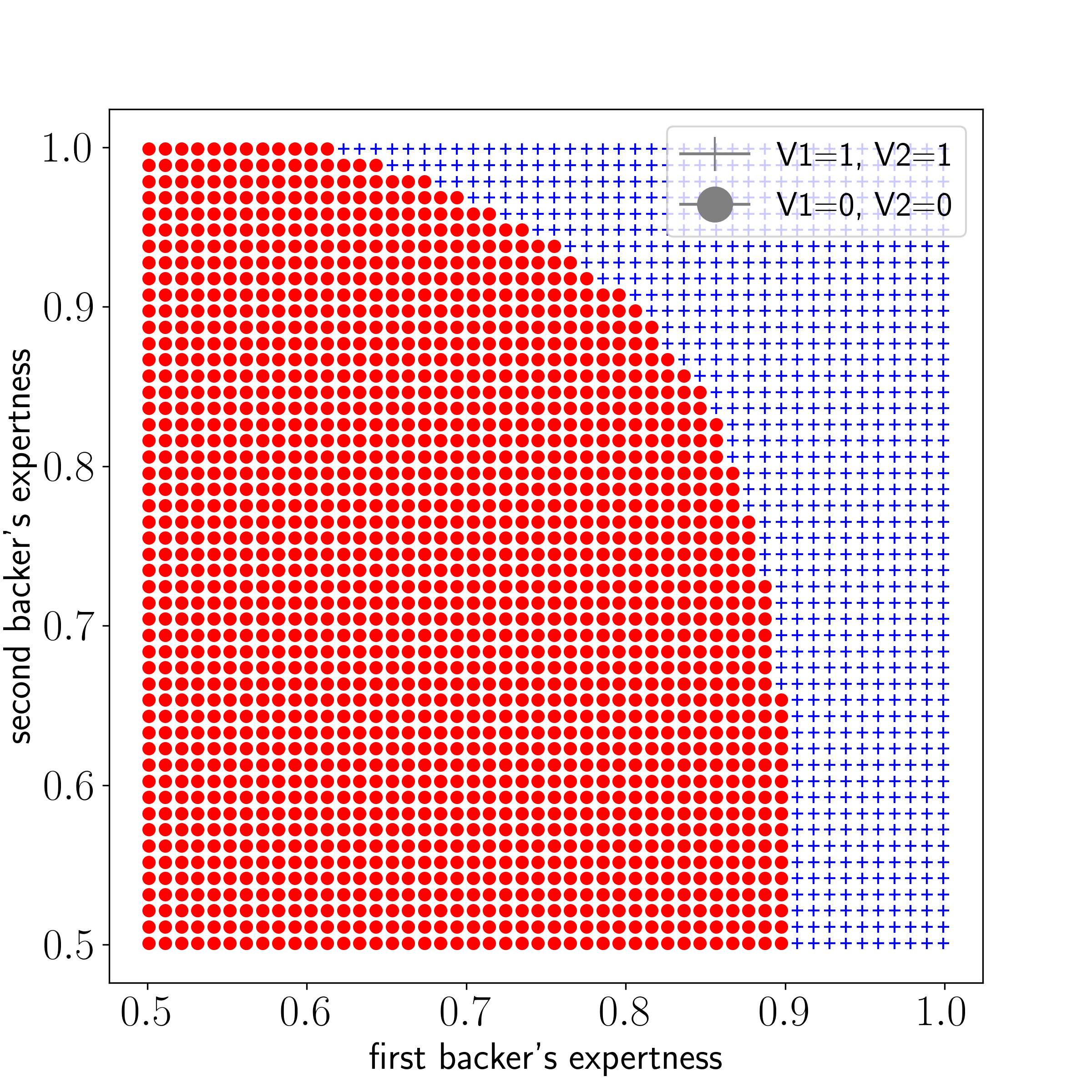}
            \caption{Tight competition, $c=0.4$}   
            \label{fig:tight_c04}
        \end{subfigure}
         \begin{subfigure}[b]{0.4\textwidth}
            \centering
            \includegraphics[width=0.8\textwidth]{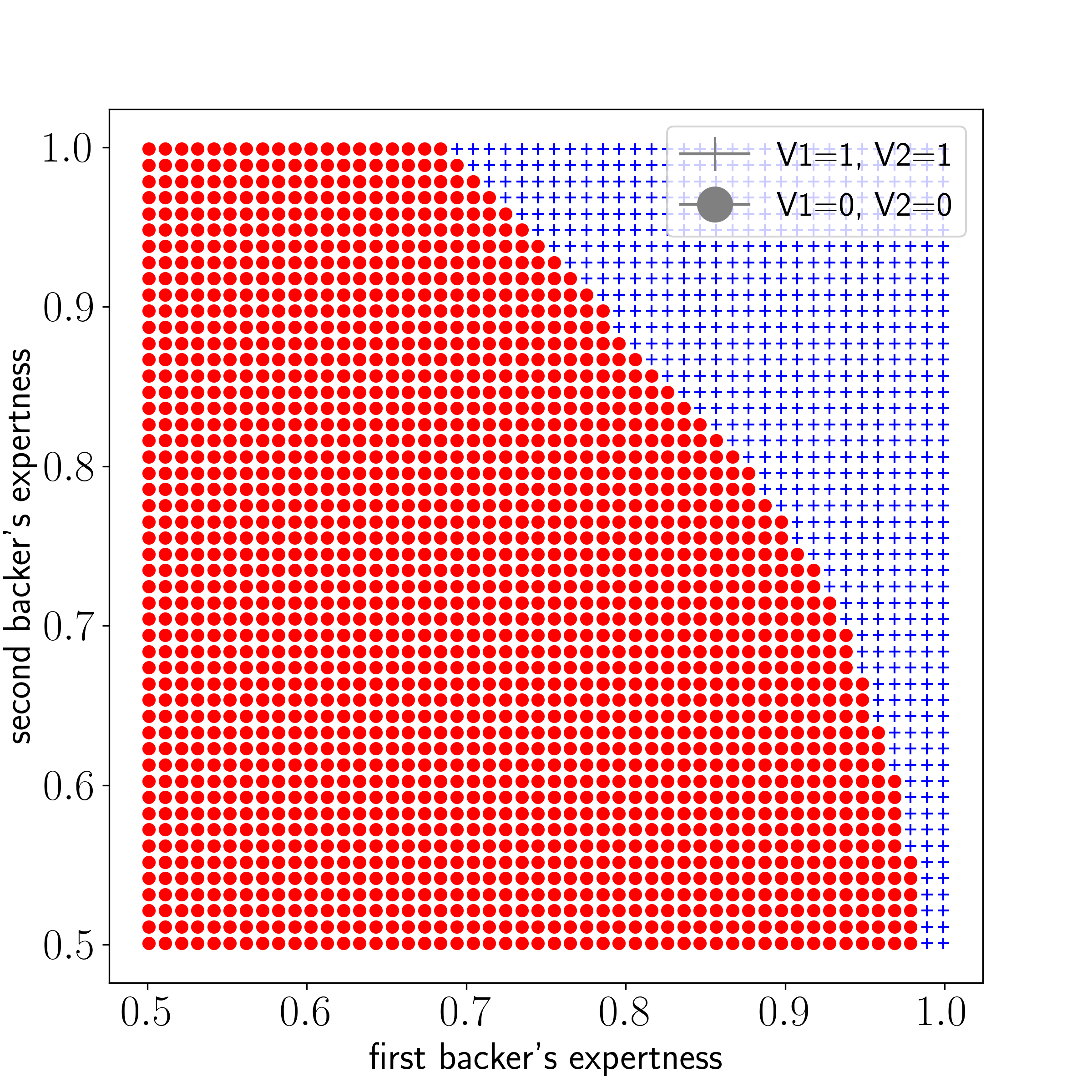}
            \caption{Relaxed competition, $c=0.4$}
            \label{fig:rel_c04}
        \end{subfigure}        
       \caption{Equilibrium Quality Strategies under Backer-asymmetric Expertness} 
    \label{fig:qual_decs_backers} \vspace{-15pt}
\end{figure}

We present the equilibrium quality strategies for all backer expertness levels in Figure \ref{fig:qual_decs_backers} under tight and relaxed competition, with two development cost levels, $c = 0.1$ and $c = 0.4$. 
A common observation in all four cases is that low backer expertness levels lead to a low-quality project equilibrium. This is expected, since backers with low expertness are not likely to infer and appreciate high-quality projects. Creators have more incentives to develop costly high-quality products only when backer expertness levels are high enough to identify the true quality states of the proposed projects. This finding bears some similarity to the results in \citet{Chakraborty2021}, where the importance of backer expertness on the creators' incentive to develop high-quality products is demonstrated in a single-project setting.

Figure \ref{fig:qual_decs_backers} also provides insights about how the nature of competition and early backer expertness impact project quality decisions. We observe that when the development cost is small, the relaxed competition case leads to high-quality projects in a wider range of expertness conditions. Conversely, in tight competition, it is risky for both projects to offer high-quality projects at the same time, especially when the expertness levels are low. In this case, OL prompts herding behavior and increases the success chances of projects in general. However, in most expertness conditions, this is not enough to make the high-quality equilibrium beneficial, since tight competition does not allow both projects to be successful at the same time.

The role of the early backer can be observed more clearly when the development cost is high under tight competition. 
In this scenario, creators choose to develop costly high-quality products only when the expertness of the early backer is high enough to create a very strong herding behavior, and therefore benefit the high-quality projects. As we explore in Section \ref{sec:main_results}, herding increases the projects' success chances in tight competition, however, not to the same extent in relaxed competition, especially when the leader's expertness is high enough to facilitate herding, but not high enough to accurately infer the true quality. This is one reason why the threshold expertness level that the early backer must have in order to incentivize high-quality projects is significantly higher in relaxed competition (Figure \ref{fig:rel_c04}), compared to tight competition (Figure \ref{fig:tight_c04}).

\begin{figure}[h] 
    \centering
         \begin{subfigure}[b]{0.35\textwidth}
            \centering
            \includegraphics[width=\textwidth]{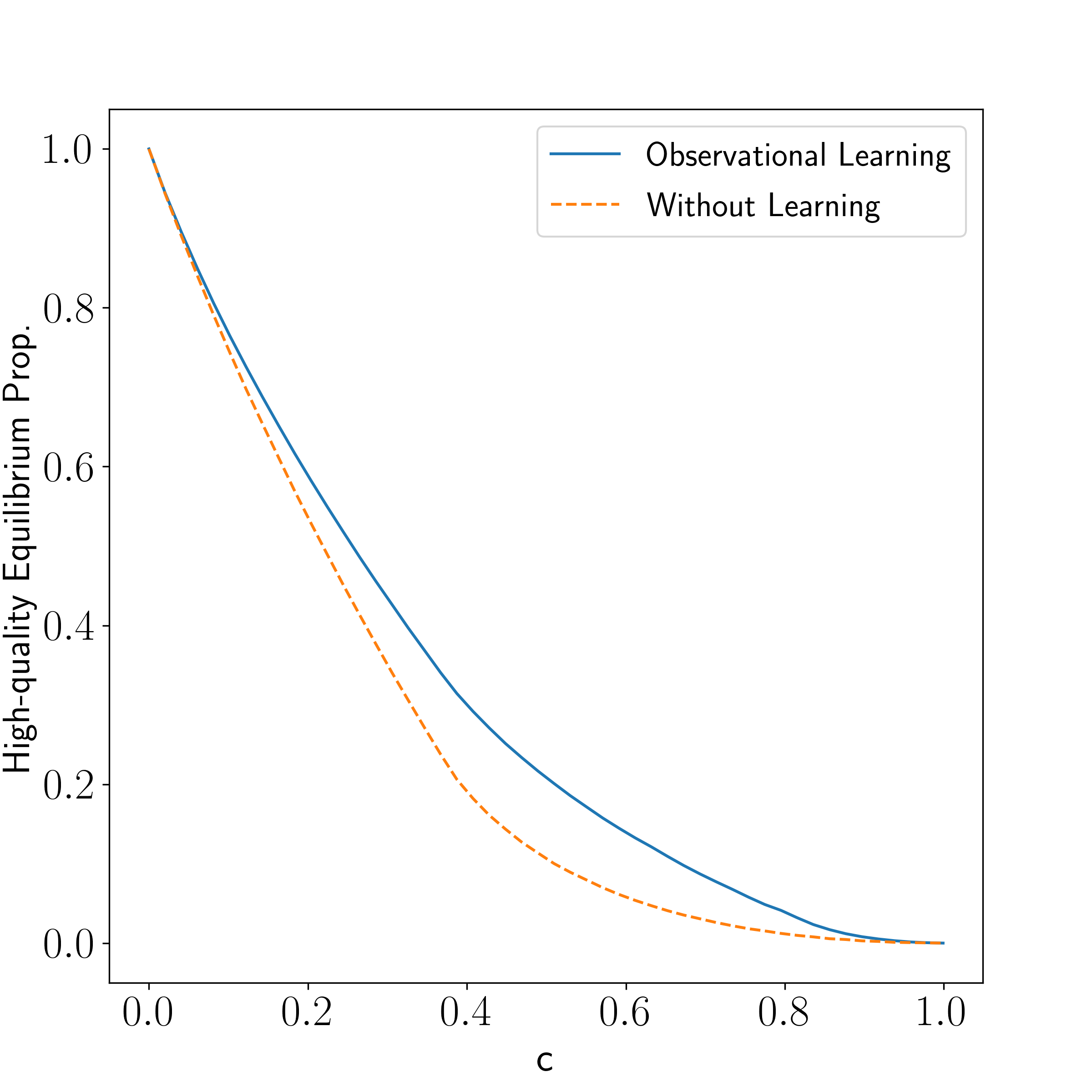}
            \caption{Tight competition}   
            \label{fig:dec_overall_backers_tight}
        \end{subfigure}
         \begin{subfigure}[b]{0.35\textwidth}
            \centering
            \includegraphics[width=\textwidth]{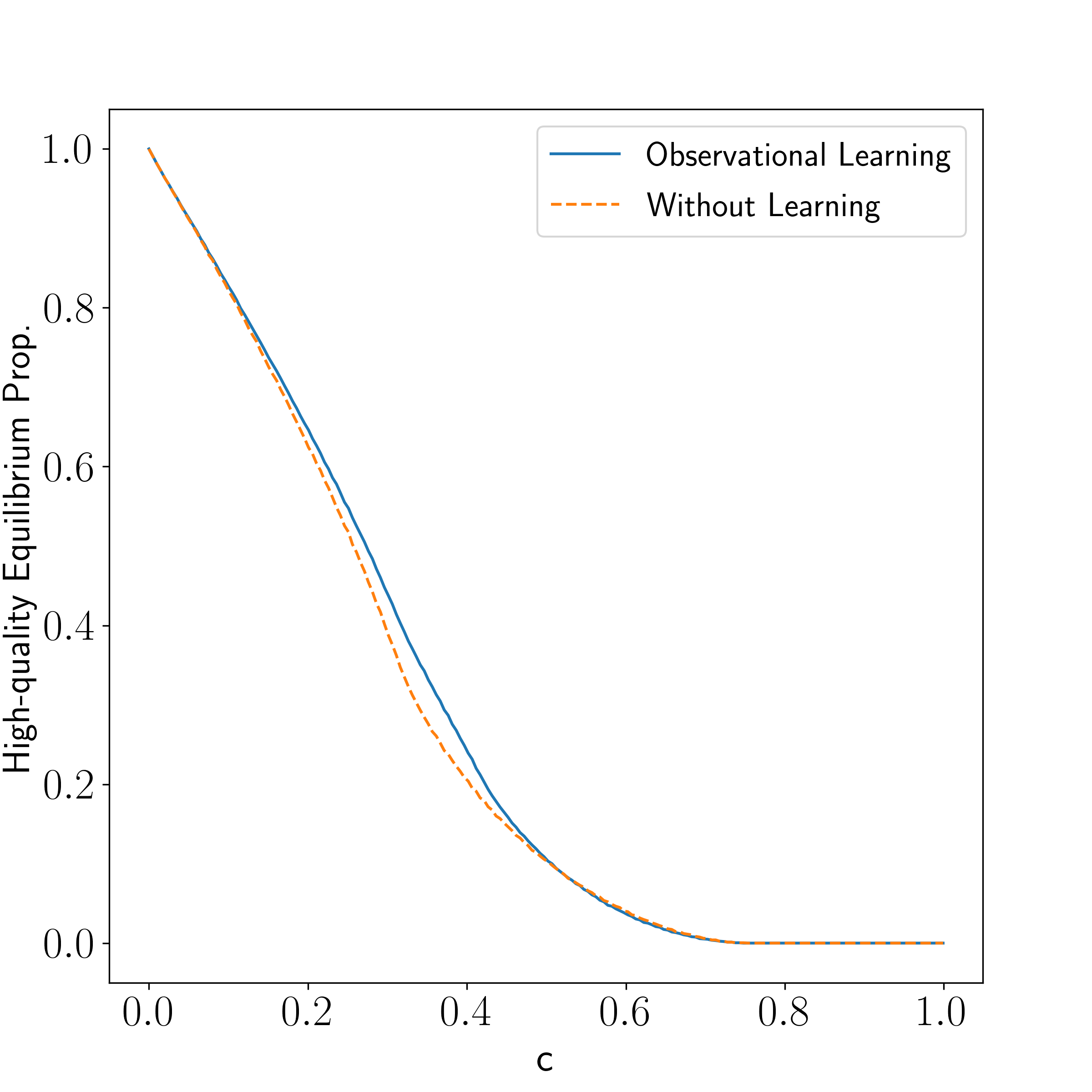}
            \caption{Relaxed competition}  
            \label{fig:dec_overall_backers_rel}
        \end{subfigure}
        \caption{Proportion of Expertness Conditions under which the High-quality Equilibrium is Reached} 
        \label{fig:dec_overall_effect} \vspace{-10pt}
        \end{figure}
     
In order to isolate the effect of OL, we also compare the equilibrium quality strategies when there is OL to those when no learning takes place. As a measure of how frequently OL increases the incentive for developing high-quality products, we use the proportion of expertness conditions under which the high-quality equilibrium is reached. In Figure \ref{fig:dec_overall_effect}, we present this measure with and without OL, in tight and relaxed competition, and for various development costs. It can be observed that in tight competition, OL encourages the creators to develop high-quality products under a wider set of expertness conditions compared to the no-learning setting, regardless of the development cost. In the case of relaxed competition, the outcome is similar when the development cost is not extremely low or high. More generally, these findings imply that crowdfunding platforms that facilitate OL may encourage creators to develop more high-quality products than they would on a platform without OL support, especially when the competition among the projects is tight.

\subsection{Product Quality Decisions under Project-Asymmetric Expertness} \label{sec:decs_project}

So far, we assumed that each backer's expertness level is the same for both projects. However, there may be many practical situations where the backers' expertness levels differ based on the project, for example, when one project creator is more experienced or well-known than the other in that particular product category, or when one project description is much more informative than the other. To understand the impact of OL on such situations, in this section, we extend our model to a case where backers have different expertness levels on the two crowdfunding projects. To observe the impact of expertness asymmetry among the two projects in isolation, we neglect the differences among backers and assume that all backers have the same expertness level for the same project.

\begin{figure}[h]
    \centering
     \begin{subfigure}[b]{0.4\textwidth}
            \centering
            \includegraphics[width=0.8\textwidth]{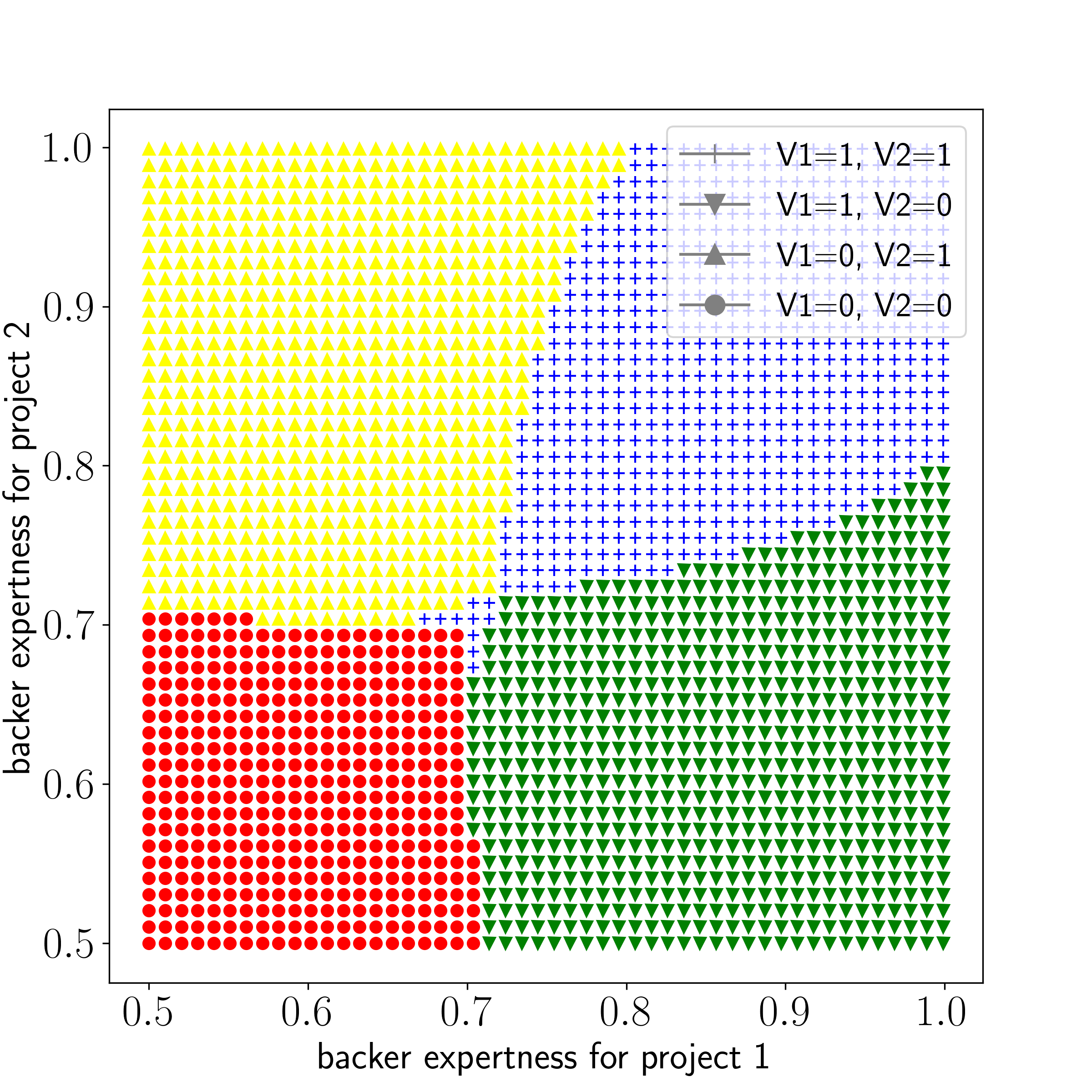}
            \caption{Tight competition, $c=0.1$} 
            \label{fig:tight_c01}
        \end{subfigure}
         \begin{subfigure}[b]{0.4\textwidth}
            \centering
            \includegraphics[width=0.8\textwidth]{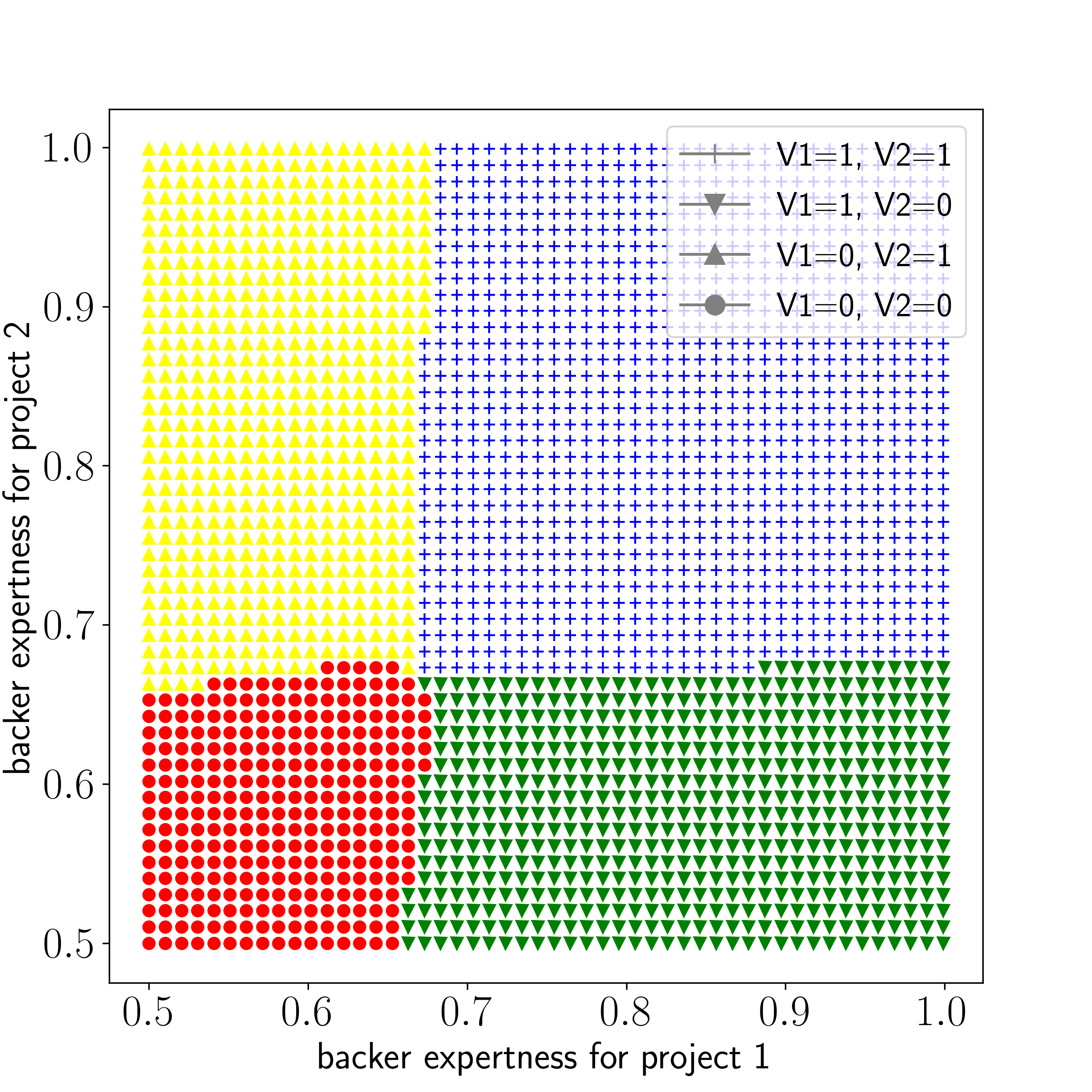}
            \caption{Relaxed competition, $c=0.1$}    
            \label{fig:rel_c01}
            \end{subfigure}
            
             \begin{subfigure}[b]{0.4\textwidth}
            \centering
            \includegraphics[width=0.8\textwidth]{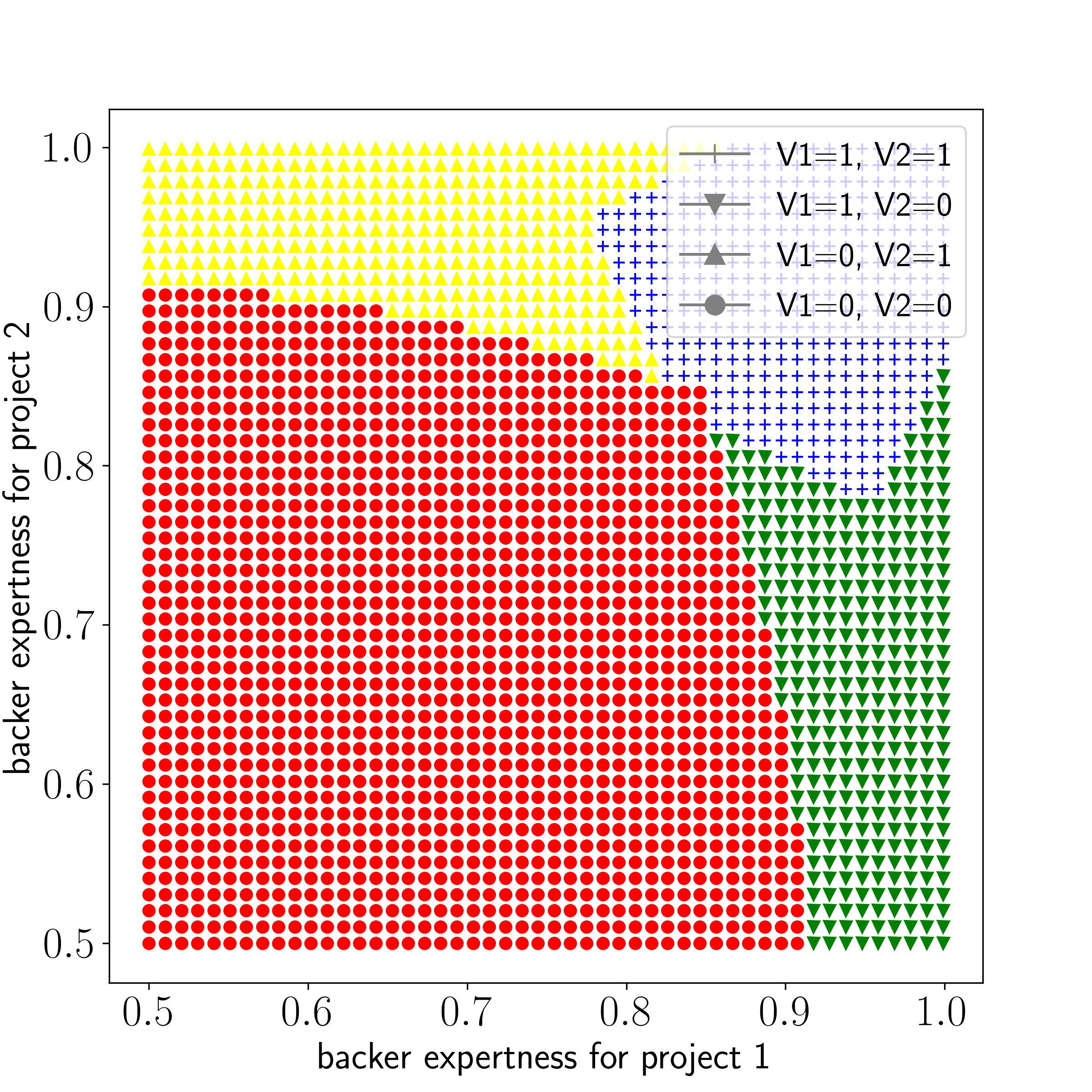}
            \caption{Tight competition, $c=0.4$}
            \label{fig:tight_c01}
        \end{subfigure}
         \begin{subfigure}[b]{0.4\textwidth}
            \centering
            \includegraphics[width=0.8\textwidth]{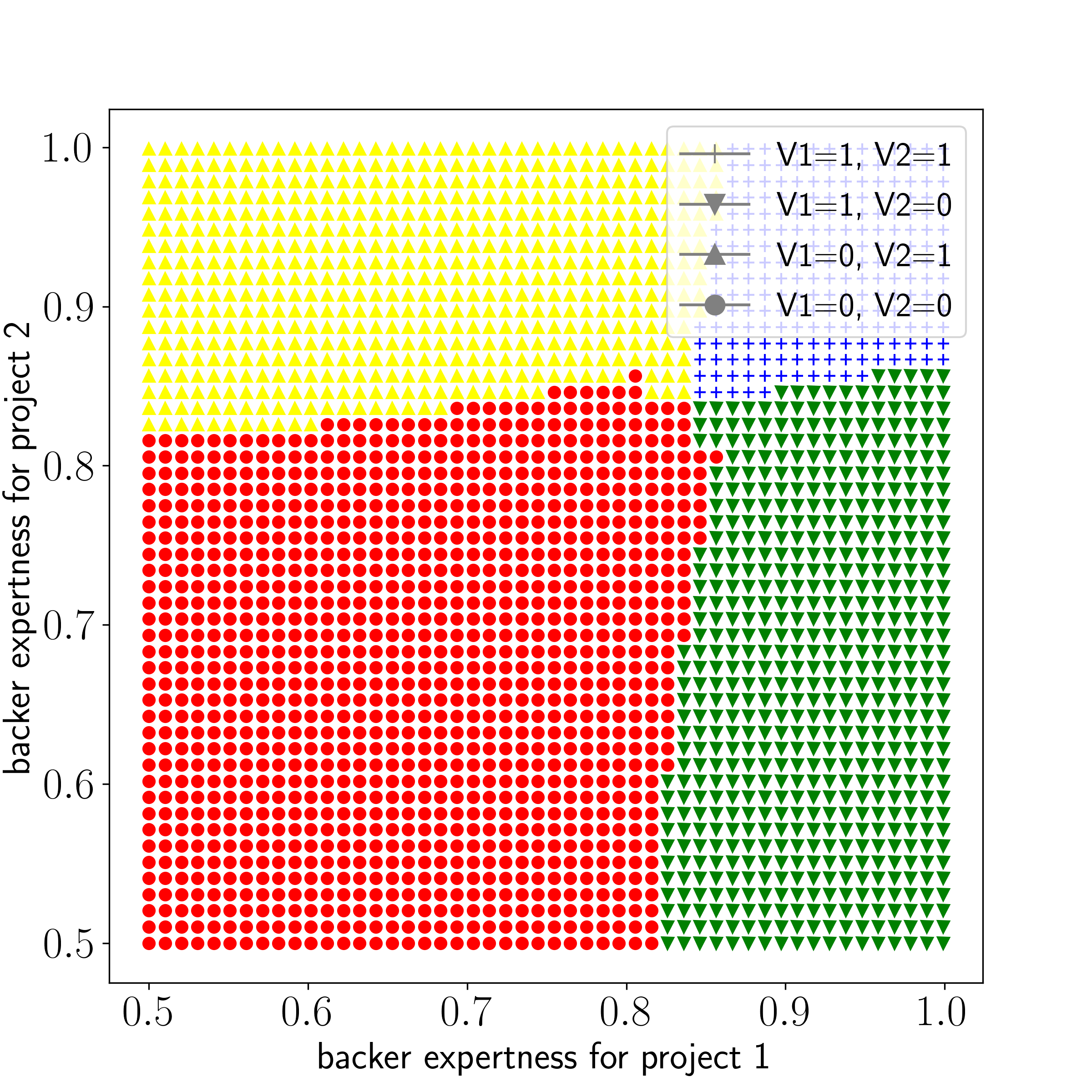}
            \caption{Relaxed competition, $c=0.4$} 
            \label{fig:rel_c01}
        \end{subfigure}        
       \caption{Equilibrium Quality Strategies under Project-asymmetric Expertness.} 
    \label{fig:qual_decs_projects}
\end{figure}

In our investigation, we let $p^1$ and $p^2$ denote the backers' expertness levels for the first and second projects, respectively. Similarly to earlier notation, we let $\bp = (p^1, p^2)$ denote the expertness conditions. Because of the information asymmetry among projects, the net profitability functions are not necessarily symmetric, and thus strategies that are asymmetric in terms of project quality, i.e., $\bV^*(\bp) = (1,0)$ and $\bV^*(\bp) = (0,1)$ can also be the equilibrium.

In Figure \ref{fig:qual_decs_projects}, we present the equilibrium quality strategies under various project-asymmetric expertness levels for two different development costs, $c = 0.1$ and $c = 0.4$, in tight and relaxed competition. As before, we observe that increased development cost results in less incentive for creators to propose high-quality projects. The development of high-quality projects becomes viable only when backers have high levels of expertness to identify the true quality of the projects. 
A common observation in all four parameter settings is that when backers have low expertness for both projects, creators tend to develop low-quality projects, and vice versa. When the expertness level is low for both projects, backers do not have the means to understand and appreciate the value of the project, hence, the investment needed to develop a high-quality project does not provide enough return. Therefore, higher backer expertness levels incentivize high project quality. These effects become more pronounced with tighter competition and higher development cost.

In situations where the expertness of backers differ significantly between two projects, a mixed-equilibrium strategy may arise where one project opts for high-quality and the other for low-quality. Specifically, as the expertness level for project $i$, $p^i$, becomes larger compared to that for project $j$, $p^j$, creator $i$ is more incentivized to develop a high-quality project, while creator $j$ becomes less motivated to do so. In such situations, developing a costly high-quality project can be riskier for the creator for whom the backers have less information, because the backers are not likely to recognize the true quality state of this project. On the other hand, the creator who has the expertness advantage can dominantly use the high-quality strategy, because the backers can understand and appreciate the high-quality status of this project.

Our analysis on the impact of OL in small systems highlights the significance of the nature of crowdfunding competition and the quality difference between projects. For example, we discover that OL improves the crowdfunding platform's profit in tight competition, while decreasing it in relaxed competition. Additionally, we find that in tight competition, OL is beneficial for platform effectiveness except when both projects are low-quality, while in relaxed competition, OL improves effectiveness except when both projects are high-quality. These insights imply that the strategy for allowing OL among backers depends on these two crucial factors. In the next section, we demonstrate that these findings remain valid in larger systems with $N>>2$. 

\section{Validation in Larger Systems} 

\label{sec:valid}
In this section, we analyze crowdfunding platforms with a large number of backers to show the robustness of our earlier observations. We specifically focus on platform effectiveness for brevity, as it embodies other metrics such as project success chances.

As argued in \citet{Zhang2015}, when the information set grows, the number of paths in the history grows exponentially, which makes it difficult to compute posterior probabilities and the resulting decision probabilities. This is also the case in our system when the number of backers, $N$, increases. Therefore, in this section, we use simulations to evaluate platform effectiveness in larger systems for a variety of parameter settings.

The level of competition depends on \begin{enumerate*}[label = \itshape (\roman*)] \item the number of pledges required for a project to succeed (henceforth referred to as the \textit{target pledge count} for brevity), $\Tilde{n}$, and \item the total number of backers, $N$, who will visit the crowdfunding platform and review the projects.\end{enumerate*} Competition gets tighter as the target pledge count increases, and the backer population size decreases. To examine platform effectiveness in large systems, we conduct simulations with two backer population sizes, $N = 50$ and $N = 100$, and 10 different target pledge counts, $\Tilde{n}$. We display the impact of OL on platform effectiveness for these scenarios in Figure \ref{fig:largeN}.

When two projects have high quality, OL is the most effective when both projects need around half of the backers to pledge. This is because, in this case, the population has just enough funds for the two high-quality projects to succeed. OL ensures that these funds can be used to support both projects which increases the platform effectiveness significantly. 

We observe in Figure \ref{fig:largeN} that in both backer population sizes, the impact of OL on platform effectiveness changes as target pledge count becomes sufficiently large. When projects require more than 15 (30) pledges from a population of 50 (100) backers to succeed, the impact of OL is in line with the behavior we observed in our earlier results with two-backer systems in tight competition: OL improves platform effectiveness when at least one project is high-quality, but can harm platform effectiveness when both projects are low-quality. In situations where projects require fewer pledges to succeed, the impact of OL resembles what we observed in two-backer systems under relaxed competition: OL can harm platform effectiveness when both projects are high-quality, but improves it in all other quality states. 

It is worth mentioning that although the two larger systems with 50 and 100 backers exhibit similar behavior, the improvement potential of OL is much higher in the larger system ($N = 100$), particularly when one of the projects is high-quality and the other is not. This can be attributed to the increased number of pledging decisions involved in the larger system, which provides more opportunities for learning to have an impact on the outcome. Overall, these findings validate our results from the two-backer systems, indicating that the impact of OL on platform effectiveness remains consistent as the number of backers increases. Furthermore, our results suggest that OL can be even more effective in larger systems, highlighting the potential of OL in modern crowdfunding platforms such as Kickstarter.com, where campaigns are often exposed to a large number of backers who actively interact.

\vspace{10pt}
\begin{figure}[h]
     \begin{subfigure}[b]{0.45\textwidth}
        \centering
           \begin{tikzpicture}[scale=0.6]
\begin{axis}[
    axis lines = left,
    xlabel = \(\Tilde{n}\),
    ylabel = \(\effect^{OL}(\bVp) - \effect^{NL}(\bVp)\),
    xtick={5,10,...,50},
    xmin= 2.5,
    xmax=52.5,
    ymin=-10,
    ymax=45,
    legend style={anchor=north west}
]
\addplot[
    color=red,
    mark=o,
    ]
    coordinates {
    (5, 2.73)(10,1.68)(15, 0.25)(20,4.18)(25,9.89) (30,6.05)(35, 2.17)(40, 0.87)(45, 0.34)
      };
    \addplot[
    color=blue,
    mark=square,
    ]
    coordinates {
    (5, 13.52)(10,11.30)(15, 7.05)(20,4.23)(25,6.48) (30,8.98)(35, 9.74)(40, 9.81)(45, 8.13)
      };
     \addplot[
    color=green,
    mark=triangle,
    ]
    coordinates {
    (5, 13.22)(10,12.00)(15, 5.61)(20,-4.48)(25,-8.29) (30,-8.35)(35, -8.21)(40, -8.14)(45, -8.12)
      };
 \end{axis}
\end{tikzpicture}
    \caption{$N = 50$} \label{fig:50}
    \end{subfigure} \hspace{-20pt}
        \begin{subfigure}[b]{0.45\textwidth}
            \centering
            \begin{tikzpicture}[scale=0.6]
\begin{axis}[
    axis lines = left,
    xlabel = \(\Tilde{n}\),
    ylabel = \(\effect^{OL}(\bVp) - \effect^{NL}(\bVp)\),
    xtick={10,20,...,100},
    xmin= 7.5,
    xmax=102.5,
    ymin=-10,
    ymax=45,
    legend style={anchor=north west},
    legend cell align={left}
]
\addplot[
    color=red,
    mark=o,
    ]
    coordinates {
    (10, 6.83)(20,5.35)(30,2.13)(40,12.07)(50,25.67)(60,9.13)(70,2.75)(80,1.16)(90,0.47)
      };
   \addlegendentry{Both projects are high-quality} 
    \addplot[
    color=blue,
    mark=square,
    ]
    coordinates {
   (10,  35.98)(20,34.90)(30,29.36)(40,25.19)(50,34.02)(60,40.37)(70,41.99)(80,41.95)(90,38.45)
      };
    \addlegendentry{Only one project is high-quality} 
     \addplot[
    color=green,
    mark=triangle,
    ]
    coordinates {
    (10,  33.98)(20,34.96)(30,25.16)(40,2.48)(50,-1.43)(60,-0.63)(70,-0.26)(80,-0.10)(90,-0.03)
      };
    \addlegendentry{Both projects are low-quality} 
 \end{axis}
\end{tikzpicture}
    \caption{\small $N=100$} \label{fig:100}
        \end{subfigure} 

\vspace{10pt}
\begin{minipage}{\linewidth}
{\scriptsize \textit{Note:} Each plotted value is an average over 25 distinct systems with a wide variety of expertness conditions. The expertness of the first half of the backers is fixed to $p_1$ and that of the second half is fixed to $p_2$, and all $(p_1, p_2)$ pairs with $p_1, p_2 \in \{0.55, 0.65, 0.75, 0.85, 0.95\}$ are considered. For each system, 50\,000 simulations are performed.} \vspace{3pt}
\end{minipage}
       \caption{\small Effect of OL on Platform Effectiveness, $\effect^{OL}(\bVp) - \effect^{NL}(\bVp)$, in Large Systems}
    \label{fig:largeN}
\end{figure}

\section{Conclusion and Discussions}\label{sec:conc}
Reward-based crowdfunding platforms, such as Kickstarter.com, have gained popularity among creators seeking funding for their products by offering backers rewards for pledging to their projects. Despite their success, these platforms have challenges related to uncertainty during the post-campaign delivery phase and product quality. This uncertainty can affect backer decisions. While backers can review project descriptions to form opinions about product quality, they also rely on observational learning (OL) to gather additional information and infer the quality of products through the decisions of early backers. This can help to reduce the information asymmetry between backers and creators.

In this paper, we explore the impact of OL on crowdfunding dynamics, from the perspectives of backers, creators, and platforms compared to the situation without learning. We use a model in which backers arrive sequentially and make pledging decisions for two available crowdfunding projects without knowing the true quality states of the proposed products. Unlike previous studies, our model considers competition between projects and examines the impact of funding scarcity on outcomes. Additionally, we differentiate between backers based on their expertness in identifying product quality, as explored in the empirical study by \citet{expertness1}, providing new insights on how this affects crowdfunding outcomes. This research offers a novel approach by combining theoretical OL models with previous empirical studies.

We discover that OL has the potential to enhance backers' contentedness with their pledging decisions regardless of their expertness levels and the true quality of the proposed projects. However, the improvement that OL provides on the backers' contentedness is most pronounced when early backers have high expertness, and late backers have low expertness. In this scenario, which encourages herding, non-expert late backers can follow the decisions of expert early backers and decrease their regret. This finding indicates that the participation of expert backers in the early stages of the crowdfunding campaign is crucial for non-expert backers, supporting the conclusions from previous empirical studies.

From the creators' perspective, our findings on the projects' success probabilities highlight the importance of competition in crowdfunding outcomes, showing how the availability of potential funding affects the impact of OL on crowdfunding dynamics. We find that when competition is tight between projects due to scarce funding, OL improves the success probabilities, even for low-quality projects. However, in relaxed competition, OL can impede success, even for high-quality projects. The preference for OL in tight competition is due to the herding behavior it creates, which is critical for success in tight competition as projects require a large proportion of the total available funding to meet their targets. However, herding can be detrimental in relaxed competition, where both projects can be successful if funding is distributed evenly, rather than focused on one project.

Our findings on the impact of OL on the short-term (platform profit) and long-term performance (platform effectiveness) of the crowdfunding platforms reveal another important dimension. We find that crowdfunding outcomes are influenced by both the availability of potential funding and the quality difference between the proposed projects. 
The benefit of OL is found to be prevalent when there are differences in the quality of the proposed projects. For example, we find that OL improves platform effectiveness in tight competition, except when both projects are low-quality, and in relaxed competition, OL is shown to improve platform effectiveness, except when both projects are high-quality.

Additionally, we explore the impact of OL on the quality decisions of creators competing for funding from backers. Specifically, we examine whether OL provides an incentive for project creators to develop high-quality projects. We find that OL incentivizes high-quality products, especially in crowdfunding campaigns subject to tight competition.

One potential future research direction is the extension of our model to a case with more than two crowdfunding projects, and also considering that expertness might vary between projects, as well as backers. It is also an interesting avenue to explore different OL models. In this paper, we use OL under the rationality assumption. There could be cases where backers may have bounded rationality. Alternatively, behavioral OL models, as proposed in \citet{Qiu2017}, can be used. 


\bibliographystyle{pomsref}

 \let\oldbibliography\thebibliography
 \renewcommand{\thebibliography}[1]{%
    \oldbibliography{#1}%
    \baselineskip14pt 
    \setlength{\itemsep}{10pt}
 }
\bibliography{main}

\begin{thebibliography}{30}
\expandafter\ifx\csname natexlab\endcsname\relax\def\natexlab#1{#1}\fi
\expandafter\ifx\csname url\endcsname\relax
  \def\url#1{{\tt #1}}\fi
\expandafter\ifx\csname urlprefix\endcsname\relax\def\urlprefix{URL }\fi
\expandafter\ifx\csname urlstyle\endcsname\relax
  \expandafter\ifx\csname doi\endcsname\relax
  \def\doi#1{doi:\discretionary{}{}{}#1}\fi \else
  \expandafter\ifx\csname doi\endcsname\relax
  \def\doi{doi:\discretionary{}{}{}\begingroup \urlstyle{rm}\Url}\fi \fi

\bibitem[{Acemoglu et~al.(2011)Acemoglu, Dahleh, Lobel, and Ozdaglar}]{Acemoglu2011}
Acemoglu, D., M.~A. Dahleh, I.~Lobel, A.~Ozdaglar. 2011.
\newblock Bayesian learning in social networks.
\newblock {\it The Review of Economic Studies\/}, { 78} 1201-1236.

\bibitem[{Alaei et~al.(2016)Alaei, Malekian, and Mostagir}]{RW}
Alaei, Saeed, Azarakhsh Malekian, Mohamed Mostagir. 2016.
\newblock Anticipatory random walks: {A} dynamic model of crowdfunding.
\newblock Working Paper.

\bibitem[{Astashkina and Marinesi(2022)}]{efficiency}
Astashkina, Ekaterina, Simone Marinesi. 2022.
\newblock All-or-nothing vs keep-it-all: Comparing campaign designs in rewards-based crowdfunding platforms.
\newblock Working Paper.

\bibitem[{Banerjee(1992)}]{Banerjee1992}
Banerjee, A.~V. 1992.
\newblock A simple model of herd behavior.
\newblock {\it The Quarterly Journal of Economics\/}, { 107} 797-817.

\bibitem[{Barnett(2015)}]{Barnett2015}
Barnett, Chance. 2015.
\newblock Trends show crowdfunding to surpass vc in 2016.
\newblock {\it Forbes\/}, \urlprefix\url{https://www.forbes.com/sites/chancebarnett/2015/06/09/trends-show-crowdfunding-to-surpass-vc-in-2016/}.

\bibitem[{Belleflamme et~al.(2015)Belleflamme, Omrani, and Peitz}]{Belleflamme2015}
Belleflamme, Paul, Nessrine Omrani, Martin Peitz. 2015.
\newblock The economics of crowdfunding platforms.
\newblock {\it Information Economics and Policy\/}, { 33} 11-28.

\bibitem[{Besbes and Sauré(2016)}]{Besbes16}
Besbes, Omar, Denis Sauré. 2016.
\newblock Product assortment and price competition under multinomial logit demand.
\newblock {\it Production and Operations Management\/}, { 25} (1), 114-127.

\bibitem[{Bi et~al.(2019)Bi, Geng, and Liu}]{bi}
Bi, Gongbing, Botao Geng, Lindong Liu. 2019.
\newblock On the fixed and flexible funding mechanisms in reward-based crowdfunding.
\newblock {\it European Journal of Operational Research\/}, { 279} (1), 168-183.

\bibitem[{Bikhchandani et~al.(1998)Bikhchandani, Hirshleifer, and Welch}]{Bikhchandani1998}
Bikhchandani, Sushil, David Hirshleifer, Ivo Welch. 1998.
\newblock Learning from the behavior of others: Conformity, fads, and informational cascades.
\newblock {\it Journal of Economic Perspectives\/}, { 12} 151-170.

\bibitem[{Candoğan et~al.(2021)Candoğan, Cornelius, Gokpinar, Körpeoğlu, and Tang}]{can}
Candoğan, Sıdıka~Tunç, Philipp~B. Cornelius, Bilal Gokpinar, Ersin Körpeoğlu, Christopher~S. Tang. 2021.
\newblock Product development in crowdfunding: Theoretical and empirical analysis.
\newblock Working Paper.

\bibitem[{Chakraborty et~al.(2023)Chakraborty, Ma, and Swinney}]{chak2}
Chakraborty, Soudipta, Anyi Ma, Robert Swinney. 2023.
\newblock Designing rewards-based crowdfunding campaigns for strategic (but distracted) contributors.
\newblock {\it Naval Research Logistics (NRL)\/}, { 70} (1), 3-20.

\bibitem[{Chakraborty and Swinney(2021)}]{Chakraborty2021}
Chakraborty, Soudipta, Robert Swinney. 2021.
\newblock Signaling to the crowd: Private quality information and rewards-based crowdfunding.
\newblock {\it Manufacturing and Service Operations Management\/}, { 23} 155-169.

\bibitem[{Cong and Xiao(2021)}]{Cong2017}
Cong, Lin~William, Yizhou Xiao. 2021.
\newblock Up-cascaded wisdom of the crowd.
\newblock Working Paper.

\bibitem[{Cornelius and Gokpinar(2020)}]{cornel}
Cornelius, Philipp~B., Bilal Gokpinar. 2020.
\newblock The role of customer investor involvement in crowdfunding success.
\newblock {\it Management Science\/}, { 66} (1), 452-472.

\bibitem[{Du et~al.(2022)Du, Hu, and Wu}]{du}
Du, Longyuan, Ming Hu, Jiahua Wu. 2022.
\newblock Contingent stimulus in crowdfunding.
\newblock {\it Production and Operations Management\/}, { 31} (9), 3543-3558.

\bibitem[{Hu et~al.(2015)Hu, Li, and Shi}]{MarkCF}
Hu, Ming, Xi~Li, Mengze Shi. 2015.
\newblock Product and pricing decisions in crowdfunding.
\newblock {\it Marketing Science\/}, { 34} (3), 331-345.

\bibitem[{Huang et~al.(2013)Huang, Leng, and Parlar}]{Huang2013}
Huang, Jian, Mingming Leng, Mahmut Parlar. 2013.
\newblock Demand functions in decision modeling: A comprehensive survey and research directions.
\newblock {\it Decision Sciences\/}, { 44} 557-609.

\bibitem[{Kim and Viswanathan(2019)}]{expertness1}
Kim, Keongtae, Siva Viswanathan. 2019.
\newblock The experts in the crowd: The role of experienced investors in a crowdfunding market.
\newblock {\it MIS Quarterly\/}, { 43} (2), 347–372.

\bibitem[{Li and Cao(2021)}]{Li21}
Li, He, Erbao Cao. 2021.
\newblock Competitive crowdfunding under asymmetric quality information.
\newblock {\it Annals of Operations Research\/}, \doi{https://doi.org/10.1007/s10479-021-03939-y}.

\bibitem[{Li et~al.(2020)Li, Duan, and Ransbotham}]{Li20}
Li, Zhuoxin, Jason~A. Duan, Sam Ransbotham. 2020.
\newblock Coordination and dynamic promotion strategies in crowdfunding with network externalities.
\newblock {\it Production and Operations Management\/}, { 29} (4), 1032-1049.

\bibitem[{Liu et~al.(2022)Liu, Liu, and Shen}]{Liu2022}
Liu, Jue, Xiaofeng Liu, Houcai Shen. 2022.
\newblock Reward-based crowdfunding: The role of information disclosure.
\newblock {\it Decision Sciences\/}, { 53} 390-422.

\bibitem[{Liu et~al.(2021)Liu, Liu, and Balachander}]{Liu21}
Liu, Qiang, Xiaofeng Liu, Subramanian Balachander. 2021.
\newblock Crowdfunding project design: Optimal product menu and funding target.
\newblock {\it Production and Operations Management\/}, { 30} (10), 3800-3811.

\bibitem[{Manuel(2022)}]{Manuel2022}
Manuel, Ryan. 2022.
\newblock Crowdfunding market size worth \$1.30 billion by 2028 - million insights.
\newblock {\it Bloomberg\/}, \urlprefix\url{https://www.bloomberg.com/press-releases/2022-03-16/crowdfunding-market-size-worth-1-30-billion-by-2028-million-insights}.

\bibitem[{Metz(2016)}]{Metz2016}
Metz, Rachel. 2016.
\newblock How pebble is killing it on kickstarter.
\newblock {\it MIT Technology Review\/}, \urlprefix\url{https://www.technologyreview.com/2016/06/02/159865/how-pebble-is-killing-it-on-kickstarter/}.

\bibitem[{Newman(2014)}]{Newman2014}
Newman, Jared. 2014.
\newblock Kreyos smartwatch: We’re not a scam, just a trainwreck.
\newblock {\it PCWorld\/}, \urlprefix\url{https://www.pcworld.com/article/435268/kreyos-smartwatch-were-not-a-scam-just-a-trainwreck.html}.

\bibitem[{Peng et~al.(2020)Peng, Zhang, Nie, Zhu, and Du}]{peng}
Peng, Jing, Jianghua Zhang, Tengfei Nie, Yangguang Zhu, Shaofu Du. 2020.
\newblock Pricing and package size decisions in crowdfunding.
\newblock {\it Transportation Research Part E: Logistics and Transportation Review\/}, { 143} 102091.

\bibitem[{Qiu and Whinston(2017)}]{Qiu2017}
Qiu, Liangfei, Andrew~B. Whinston. 2017.
\newblock Pricing strategies under behavioral observational learning in social networks.
\newblock {\it Production and Operations Management\/}, { 26} 1249-1267.

\bibitem[{Xiao et~al.(2021)Xiao, Ho, and Che}]{momentum}
Xiao, Shengsheng, Yi-Chun~(Chad) Ho, Hai Che. 2021.
\newblock Building the momentum: Information disclosure and herding in online crowdfunding.
\newblock {\it Production and Operations Management\/}, { 30} (9), 3213-3230.

\bibitem[{Zhang et~al.(2022)Zhang, Savin, and Veeraraghavan}]{revenue}
Zhang, Jiding, Sergei Savin, Senthil Veeraraghavan. 2022.
\newblock Revenue management in crowdfunding.
\newblock {\it Manufacturing \& Service Operations Management\/}, \doi{https://doi.org/10.1287/msom.2022.1147}.

\bibitem[{Zhang et~al.(2015)Zhang, Liu, and Chen}]{Zhang2015}
Zhang, Jurui, Yong Liu, Yubo Chen. 2015.
\newblock Social learning in networks of friends versus strangers.
\newblock {\it Marketing Science\/}, { 34} 573-589.

\end{thebibliography}






%
%
%

\end{document}


\maketitle

To assess crowdfunding performance measures in a two-backer system, we present the derivation of posterior beliefs using Bayesian rationality. In our analysis, we primarily focus on backer-asymmetric expertness for the backers, as detailed in Section 5. Nevertheless, we also explore results that are applicable when considering project-asymmetric expertness, as discussed in Section 6.1. To ensure that our derivations can accommodate both scenarios, we take into account project and backer asymmetric expertness. We denote the expertness level of backer $j$ for project $i$ as $p_j^{i}$, where $i=1,2$ (project index) and $j=1,2$ (backer index). As in the main manuscript, for clarity of exposition, we
henceforth use the ``she'' pronoun to refer to the early backer, and ``he'' to refer to the late backer.

\section{Derivation of Posterior Beliefs in the 2-backers System}
Let us start with the first backer. The first backer only has information from her private signals $s_1^1$ and $s_1^2$. Her posterior probabilities are derived from 

{\scriptsize
\begin{equation}
      P(V^1, V^2| s_1^1, s_1^2)=P(s_1^1|V^1)P(s_1^2| V^2),
    \end{equation}}
which leads to the following depending on the signals. 
{\scriptsize
\begin{align*}
    P(V^1=1, V^2=1| s_1^1=H, s_1^2=H)&=P(V^1=0, V^2=0| s_1^1=L, s_1^2=L)=P(V^1=1, V^2=0| s_1^1=H, s_1^2=L)\\
    &=P(V^1=0, V^2=1| s_1^1=L, s_1^2=H)=p_1^1 p_1^2\\
    P(V^1=1, V^2=1| s_1^1=H, s_1^2=L)&=P(V^1=0, V^2=0| s_1^1=L, s_1^2=H)=P(V^1=1, V^2=0| s_1^1=H, s_1^2=H)\\
    &=P(V^1=0, V^2=1| s_1^1=L, s_1^2=L)=p_1^1 (1-p_1^2)\\
    P(V^1=1, V^2=1| s_1^1=L, s_1^2=H)&=P(V^1=0, V^2=0| s_1^1=H, s_1^2=L)=P(V^1=1, V^2=0| s_1^1=L, s_1^2=L)\\
    &=P(V^1=0, V^2=1| s_1^1=H, s_1^2=H)=(1-p_1^1) p_1^2\\
    P(V^1=1, V^2=1| s_1^1=L, s_1^2=L)&=P(V^1=0, V^2=0| s_1^1=H, s_1^2=H)=P(V^1=1, V^2=0| s_1^1=L, s_1^2=H)\\
    &=P(V^1=0, V^2=1| s_1^1=H, s_1^2=L)=(1-p_1^1) (1-p_1^2).
\end{align*}
}

Based on these beliefs, she will make her pledging decision according to the following probabilities conditional on the signals she receives.
\begin{changemargin}{-1.2cm}{-1.2cm} 
{\scriptsize
\begin{align}
    P(x_1=1| s_1^1=H, s_1^2=H)=&P(V^1=1, V^2=0| s_1^1=H, s_1^2=H)+0.5 P(V^1=1, V^2=1| s_1^1=H, s_1^2=H)=p_1^1(1-\frac{p_1^2}{2})\\
    P(x_1=2| s_1^1=H, s_1^2=H)=&P(V^1=0, V^2=1| s_1^1=H, s_1^2=H)+0.5 P(V^1=1, V^2=1| s_1^1=H, s_1^2=H)=p_1^2(1-\frac{p_1^1}{2})\\
    P(x_1=0| s_1^1=H, s_1^2=H)=&P(V^1=0, V^2=0| s_1^1=H, s_1^2=H)=(1-p_1^1) (1-p_1^2)\\
     P(x_1=1| s_1^1=H, s_1^2=L)=&P(V^1=1, V^2=0| s_1^1=H, s_1^2=L)+0.5 P(V^1=1, V^2=1| s_1^1=H, s_1^2=L)=p_1^1(\frac{1+p_1^2}{2})\\
     P(x_1=2| s_1^1=H, s_1^2=L)=&P(V^1=0, V^2=1| s_1^1=H, s_1^2=L)+0.5 P(V^1=1, V^2=1| s_1^1=H, s_1^2=L)=(1-p_1^2)(1-\frac{p_1^1}{2})\\
     P(x_1=0| s_1^1=H, s_1^2=L)=&P(V^1=0, V^2=0| s_1^1=H, s_1^2=L)=p_1^2(1-p_1^1)\\
     P(x_1=1| s_1^1=L, s_1^2=H)=&P(V^1=1, V^2=0| s_1^1=L, s_1^2=H)+0.5 P(V^1=1, V^2=1| s_1^1=L, s_1^2=H)=(1-p_1^1)(1-\frac{p_1^2}{2})\\
     P(x_1=2| s_1^1=L, s_1^2=H)=&P(V^1=0, V^2=1| s_1^1=L, s_1^2=H)+0.5 P(V^1=1, V^2=1| s_1^1=L, s_1^2=H)=p_1^2 (\frac{1+p_1^1}{2})\\
     P(x_1=0| s_1^1=L, s_1^2=H)=&P(V^1=0, V^2=0| s_1^1=L, s_1^2=H)=p_1^1(1-p_1^2)\\
     P(x_1=1| s_1^1=L, s_1^2=L)=&P(V^1=1, V^2=0| s_1^1=L, s_1^2=L)+0.5 P(V^1=1, V^2=1| s_1^1=L, s_1^2=L)=(1-p_1^1)(1+\frac{p_1^2}{2})\\
     P(x_1=2| s_1^1=L, s_1^2=L)=&P(V^1=0, V^2=1| s_1^1=L, s_1^2=L)+0.5 P(V^1=1, V^2=1| s_1^1=L, s_1^2=L)=(1-p_1^2)(1+\frac{p_1^1}{2})\\
     P(x_1=0| s_1^1=L, s_1^2=L)=&P(V^1=0, V^2=0| s_1^1=L, s_1^2=L)=p_1^1p_1^2
\end{align}
}%
\end{changemargin}

Next, using these we derive the following probability

{\scriptsize
\begin{equation}
    P(x_1| V^1, V^2)=\sum_{s_1^1, s_1^2} P(x_1|s_1^1, s_1^2) P(s_1^1|V^1)P(s_1^2|V^2),
\end{equation}}

which is required for the inference of the second backer, as he can observe her decision ($x_1$) but not the signals she received ($s_1^1, s_1^2$). Below, we calculate this probability for all possible true quality states and pledge decisions of the first backer.  
\begin{changemargin}{-1.2cm}{-1.2cm} 
{\scriptsize
\begin{align}
    P(x_1=1| V^1=1, V^2=1)=& p_1^1(1-\frac{p_1^2}{2}) p_1^1p_1^2+p_1^1 (\frac{1+p_1^2}{2}) p_1^1 (1-p_1^2)+ (1-p_1^1)(1-\frac{p_1^2}{2})(1-p_1^1)p_1^2+ (1-p_1^1)(\frac{1+p_1^2}{2})(1-p_1^1)(1-p_1^2)\\
    P(x_1=2| V^1=1, V^2=1)=& p_1^2(1-\frac{p_1^1}{2}) p_1^1p_1^2+(1-p_1^2) (1-\frac{p_1^1}{2}) p_1^1 (1-p_1^2)+ p_1^2 (\frac{1+p_1^1}{2})(1-p_1^1)p_1^2+ (1-p_1^2)(\frac{1+p_1^1}{2})(1-p_1^1)(1-p_1^2)\\
    P(x_1=0| V^1=1, V^2=1)=& (1-p_1^1)(1-p_1^2) p_1^1p_1^2+p_1^2 (1-p_1^1) p_1^1 (1-p_1^2)+ p_1^1 (1-p_1^2)(1-p_1^1)p_1^2+ p_1^1p_1^2(1-p_1^1)(1-p_1^2)\\
    P(x_1=1| V^1=1, V^2=0)=& p_1^1(1-\frac{p_1^2}{2}) p_1^1(1-p_1^2)+p_1^1 (\frac{1+p_1^2}{2}) p_1^1 p_1^2+ (1-p_1^1)(1-\frac{p_1^2}{2})(1-p_1^1)(1-p_1^2)+ (1-p_1^1)(\frac{1+p_1^2}{2})(1-p_1^1)p_1^2\\
        P(x_1=2| V^1=1, V^2=0)=& p_1^2(1-\frac{p_1^1}{2}) p_1^1(1-p_1^2)+(1-p_1^2) (1-\frac{p_1^1}{2}) p_1^1 p_1^2+ p_1^2 (\frac{1+p_1^1}{2})(1-p_1^1)(1-p_1^2)+ (1-p_1^2)(\frac{1+p_1^1}{2})(1-p_1^1)p_1^2
    \end{align}
    \begin{align}
    P(x_1=0| V^1=1, V^2=0)=& (1-p_1^1)(1-p_1^2)  p_1^1(1-p_1^2)+p_1^2 (1-p_1^1) p_1^1 p_1^2+ p_1^1 (1-p_1^2)(1-p_1^1)(1-p_1^2)+ p_1^1p_1^2(1-p_1^1)p_1^2\\
    P(x_1=1| V^1=0, V^2=1)=& p_1^1(1-\frac{p_1^2}{2}) (1-p_1^1)p_1^2+p_1^1 (\frac{1+p_1^2}{2}) (1-p_1^1) (1-p_1^2)+ (1-p_1^1)(1-\frac{p_1^2}{2})p_1^1p_1^2+ (1-p_1^1)(\frac{1+p_1^2}{2})p_1^1(1-p_1^2)\\
    P(x_1=2| V^1=0, V^2=1)=& p_1^2(1-\frac{p_1^1}{2})(1-p_1^1)p_1^2+(1-p_1^2) (1-\frac{p_1^1}{2}) (1-p_1^1) (1-p_1^2)+ p_1^2 (\frac{1+p_1^1}{2})p_1^1p_1^2+ (1-p_1^2)(\frac{1+p_1^1}{2})p_1^1(1-p_1^2)\\
     P(x_1=0| V^1=0, V^2=1)=& (1-p_1^1)(1-p_1^2) (1-p_1^1)p_1^2+p_1^2 (1-p_1^1) (1-p_1^1) (1-p_1^2)+ p_1^1 (1-p_1^2)p_1^1p_1^2+ p_1^1p_1^2p_1^1(1-p_1^2)\\
     P(x_1=1| V^1=0, V^2=0)=& p_1^1(1-\frac{p_1^2}{2}) (1-p_1^1)(1-p_1^2)+p_1^1 (\frac{1+p_1^2}{2}) (1-p_1^1) p_1^2+ (1-p_1^1)(1-\frac{p_1^2}{2})p_1^1(1-p_1^2)+ (1-p_1^1)(\frac{1+p_1^2}{2})p_1^1p_1^2\\
     P(x_1=2| V^1=0, V^2=0)=& p_1^2(1-\frac{p_1^1}{2})(1-p_1^1)(1-p_1^2)+(1-p_1^2) (1-\frac{p_1^1}{2}) (1-p_1^1) p_1^2+ p_1^2 (\frac{1+p_1^1}{2})p_1^1(1-p_1^2)+ (1-p_1^2)(\frac{1+p_1^1}{2})p_1^1p_1^2\\
     P(x_1=0| V^1=0, V^2=0)=& (1-p_1^1)(1-p_1^2) (1-p_1^1)(1-p_1^2)+p_1^2 (1-p_1^1) (1-p_1^1) p_1^2+ p_1^1 (1-p_1^2)p_1^1(1-p_1^2)+ p_1^1p_1^2p_1^1p_1^2
\end{align}
}%
\end{changemargin}

The probability of the early backer to fund the first project (i.e., the case $x_1=1$) for different quality cases is simplified as the following.

{\scriptsize
\begin{align}
     P(x_1=1| V^1=1, V^2=1)=&((p_1^1)^2 + (1-p_1^1)^2) (\frac{1}{2}+p_1^2-(p_1^2)^2)\\
     P(x_1=1| V^1=1, V^2=0)=&((p_1^1)^2 + (1-p_1^1)^2) (1-p_1^2+(p_1^2)^2)\\
     P(x_1=1| V^1=0, V^2=1)=&(2 p_1^1 (1-p_1^1)) (\frac{1}{2}+p_1^2-(p_1^2)^2)\\
     P(x_1=1| V^1=0, V^2=0)=&(2 p_1^1 (1-p_1^1)) (1-p_1^2+(p_1^2)^2)
\end{align}}

Note that this probability is the highest when $V^1=1$ and $V^2=0$. For example, $P(x_1=1| V^1=1, V^2=0)\geq  P(x_1=1| V^1=1, V^2=1)$ always because $1-p_1^2+(p_1^2)^2 \geq \frac{1}{2}+p_1^2-(p_1^2)^2$ always (as $1/4 \geq p_1^2 (1-p_1^2)$, $p_1^2 (1-p_1^2)$ attains its maximum at 1/2 and becomes equal to 1/4 there.). For $P(x_1=1| V^1=1, V^2=0)$, the minimum is attained at $p_1^1=p_1^2=0.5$. Likewise, we simplify the expressions for the remaining pledge decision cases. 

{\scriptsize
\begin{align}
     P(x_1=2| V^1=1, V^2=1)=&((p_1^2)^2 + (1-p_1^2)^2) (\frac{1}{2}+p_1^1-(p_1^1)^2)\\
     P(x_1=2| V^1=1, V^2=0)=&(2 p_1^2 (1-p_1^2)) (\frac{1}{2}+p_1^1-(p_1^1)^2)\\
     P(x_1=2| V^1=0, V^2=1)=&((p_1^2)^2 + (1-p_1^2)^2) (1-p_1^1+(p_1^1)^2)\\
     P(x_1=2| V^1=0, V^2=0)=&(2 p_1^2 (1-p_1^2)) (1-p_1^1+(p_1^1)^2)\\
     P(x_1=0| V^1=1, V^2=1)=&(2 p_1^1 (1-p_1^1)) (2 p_1^2 (1-p_1^2))\\
     P(x_1=0| V^1=1, V^2=0)=&(2 p_1^1 (1-p_1^1)) ((p_1^2)^2 + (1-p_1^2)^2)\\
     P(x_1=0| V^1=0, V^2=1)=&(2 p_1^2 (1-p_1^2)) ((p_1^1)^2 + (1-p_1^1)^2)\\
     P(x_1=0| V^1=0, V^2=0)=&((p_1^1)^2 + (1-p_1^1)^2) ((p_1^2)^2 + (1-p_1^2)^2)
     \end{align}}

Note that when there is no learning involved (the NL case we investigate in the main manuscript), decision probabilities of every backer will be the same as the first backer. Therefore, the evaluation of crowdfunding performance measures under NL follows from the posterior beliefs derived above. 

We now consider the posterior beliefs of the second backerunder OL. After observing his private signals $s_2^1$ and $s_2^2$, and the pledge decision of the first backer $x_1$, the second backer forms his posterior probabilities as follows. 

{\scriptsize
\begin{equation}
    P(V^1=1,V^2=1| s_2^1, s_2^2, x_1)=\frac{P(s_2^1|V^1=1) P(s_2^2|V^2=1) P(x_1| V^1=1, V^2=1)}{\sum_{V^1, V^2} P(s_2^1|V^1) P(s_2^2|V^2) P(x_1| V^1, V^2) }
\end{equation}}

Given the four possible cases for the signals ($s_2^1$ and $s_2^2$) and three possible cases for $x_1$, we derive the posterior beliefs for the 12 potential scenarios comprising the information set of the second backer. Below, we present the derivations for all scenarios. 

Suppose that $s_2^1=H, s_2^2=H, x_1=1$.

{\scriptsize
\begin{align}
    P(V^1=1,V^2=1| s_2^1=H, s_2^2=H, x_1=1)=& \frac{p_2^1 p_2^2 ((p_1^1)^2 + (1-p_1^1)^2) (\frac{1}{2}+p_1^2-(p_1^2)^2)}{(p_2^1 +(2 p_1^1(1-p_1^1))(1-2 p_2^1))(1-\frac{p_2^2}{2}+(2 p_1^2 (1-p_1^2))(p_2^2-\frac{1}{2}))}\\
    P(V^1=1,V^2=0| s_2^1=H, s_2^2=H, x_1=1)=& \frac{p_2^1 (1-p_2^2) ((p_1^1)^2 + (1-p_1^1)^2) (1-p_1^2 +(p_1^2)^2)}{(p_2^1 +(2 p_1^1(1-p_1^1))(1-2 p_2^1))(1-\frac{p_2^2}{2}+(2 p_1^2 (1-p_1^2))(p_2^2-\frac{1}{2}))}\\
    P(V^1=0,V^2=1| s_2^1=H, s_2^2=H, x_1=1)=& \frac{(1-p_2^1) p_2^2 (2 p_1^1 (1-p_1^1)) (\frac{1}{2}+p_1^2 -(p_1^2)^2)}{(p_2^1 +(2 p_1^1(1-p_1^1))(1-2 p_2^1))(1-\frac{p_2^2}{2}+(2 p_1^2 (1-p_1^2))(p_2^2-\frac{1}{2}))}\\
    P(V^1=0,V^2=0| s_2^1=H, s_2^2=H, x_1=1)=& \frac{(1-p_2^1) (1-p_2^2) (2 p_1^1 (1-p_1^1)) (1-p_1^2 +(p_1^2)^2)}{(p_2^1 +(2 p_1^1(1-p_1^1))(1-2 p_2^1))(1-\frac{p_2^2}{2}+(2 p_1^2 (1-p_1^2))(p_2^2-\frac{1}{2}))}
\end{align}}

Suppose that $s_2^1=H, s_2^2=H, x_1=2$.
{\scriptsize
\begin{align}
    P(V^1=1,V^2=1| s_2^1=H, s_2^2=H, x_1=2)=& \frac{p_2^1 p_2^2 ((p_1^2)^2 + (1-p_1^2)^2) (\frac{1}{2}+p_1^1-(p_1^1)^2)}{(p_2^2 +(2 p_1^2(1-p_1^2))(1-2 p_2^2))(1-\frac{p_2^1}{2}+(2 p_1^1 (1-p_1^1))(p_2^1-\frac{1}{2}))}\\
    P(V^1=1,V^2=0| s_2^1=H, s_2^2=H, x_1=2)=& \frac{p_2^1 (1-p_2^2) (2 p_1^2 (1-p_1^2)) (\frac{1}{2}+p_1^1-(p_1^1)^2)}{(p_2^2 +(2 p_1^2(1-p_1^2))(1-2 p_2^2))(1-\frac{p_2^1}{2}+(2 p_1^1 (1-p_1^1))(p_2^1-\frac{1}{2}))}\\
    P(V^1=0,V^2=1| s_2^1=H, s_2^2=H, x_1=2)=& \frac{(1-p_2^1) p_2^2 ((p_1^2)^2 + (1-p_1^2)^2) (1-p_1^1 +(p_1^1)^2)}{(p_2^2 +(2 p_1^2(1-p_1^2))(1-2 p_2^2))(1-\frac{p_2^1}{2}+(2 p_1^1 (1-p_1^1))(p_2^1-\frac{1}{2}))}\\
    P(V^1=0,V^2=0| s_2^1=H, s_2^2=H, x_1=2)=& \frac{(1-p_2^1) (1-p_2^2) (2 p_1^2 (1-p_1^2)) (1-p_1^1 +(p_1^1)^2)}{(p_2^2 +(2 p_1^2(1-p_1^2))(1-2 p_2^2))(1-\frac{p_2^1}{2}+(2 p_1^1 (1-p_1^1))(p_2^1-\frac{1}{2}))}
\end{align}}

Suppose that $s_2^1=H, s_2^2=H, x_1=0$.
{\scriptsize
\begin{align}
    P(V^1=1,V^2=1| s_2^1=H, s_2^2=H, x_1=0)=& \frac{p_2^1 p_2^2 (2 p_1^1 (1-p_1^1)) (2 p_1^2 (1-p_1^2))}{(1-p_2^2 + 2 p_1^2 (1-p_1^2)(2 p_2^2-1))(1-p_2^1+ 2 p_1^1 (1-p_1^1)(2 p_2^1-1))}\\
    P(V^1=1,V^2=0| s_2^1=H, s_2^2=H, x_1=0)=& \frac{p_2^1 (1-p_2^2) (2 p_1^1 (1-p_1^1)) ((p_1^2)^2 + (1-p_1^2)^2)}{(1-p_2^2 + 2 p_1^2 (1-p_1^2)(2 p_2^2-1))(1-p_2^1+ 2 p_1^1 (1-p_1^1)(2 p_2^1-1))}\\
    P(V^1=0,V^2=1| s_2^1=H, s_2^2=H, x_1=0)=& \frac{(1-p_2^1) p_2^2 (2 p_1^2 (1-p_1^2)) ((p_1^1)^2 + (1-p_1^1)^2)}{(1-p_2^2 + 2 p_1^2 (1-p_1^2)(2 p_2^2-1))(1-p_2^1+ 2 p_1^1 (1-p_1^1)(2 p_2^1-1))}\\
    P(V^1=0,V^2=0| s_2^1=H, s_2^2=H, x_1=0)=& \frac{(1-p_2^1) (1-p_2^2) ((p_1^1)^2 + (1-p_1^1)^2)) ((p_1^2)^2 + (1-p_1^2)^2)}{(1-p_2^2 + 2 p_1^2 (1-p_1^2)(2 p_2^2-1))(1-p_2^1+ 2 p_1^1 (1-p_1^1)(2 p_2^1-1))}
\end{align}}

Suppose that $s_2^1=H, s_2^2=L, x_1=1$.
{\scriptsize
\begin{align}
    P(V^1=1,V^2=1| s_2^1=H, s_2^2=L, x_1=1)=& \frac{p_2^1 (1-p_2^2) ((p_1^1)^2 + (1-p_1^1)^2) (\frac{1}{2}+p_1^2-(p_1^2)^2)}{(p_2^1+2p_1^1(1-p_1^1)(1-2 p_2^1))(\frac{1+p_2^2}{2} +2 p_1^2 (1-p_1^2) (\frac{1}{2}-p_2^2))}\\
    P(V^1=1,V^2=0| s_2^1=H, s_2^2=L, x_1=1)=& \frac{p_2^1 p_2^2 ((p_1^1)^2 + (1-p_1^1)^2) (1-p_1^2 +(p_1^2)^2)}{(p_2^1+2p_1^1(1-p_1^1)(1-2 p_2^1))(\frac{1+p_2^2}{2} +2 p_1^2 (1-p_1^2) (\frac{1}{2}-p_2^2))}\\
    P(V^1=0,V^2=1| s_2^1=H, s_2^2=L, x_1=1)=& \frac{(1-p_2^1) (1-p_2^2) (2 p_1^1 (1-p_1^1)) (\frac{1}{2}+p_1^2 -(p_1^2)^2)}{(p_2^1+2p_1^1(1-p_1^1)(1-2 p_2^1))(\frac{1+p_2^2}{2} +2 p_1^2 (1-p_1^2) (\frac{1}{2}-p_2^2))}\\
    P(V^1=0,V^2=0| s_2^1=H, s_2^2=L, x_1=1)=& \frac{(1-p_2^1) p_2^2 (2 p_1^1 (1-p_1^1)) (1-p_1^2 +(p_1^2)^2)}{(p_2^1+2p_1^1(1-p_1^1)(1-2 p_2^1))(\frac{1+p_2^2}{2} +2 p_1^2 (1-p_1^2) (\frac{1}{2}-p_2^2))}
\end{align}}

Suppose that $s_2^1=H, s_2^2=L, x_1=2$.
{\scriptsize
\begin{align}
    P(V^1=1,V^2=1| s_2^1=H, s_2^2=L, x_1=2)=& \frac{p_2^1 (1-p_2^2) ((p_1^2)^2 + (1-p_1^2)^2) (\frac{1}{2}+p_1^1-(p_1^1)^2)}{(1-\frac{p_2^1}{2}+2 p_1^1 (1-p_1^1)(p_2^1-\frac{1}{2}))(1-p_2^2 + 2 p_1^2 (1-p_1^2) (2 p_2^2-1 ))}\\
    P(V^1=1,V^2=0| s_2^1=H, s_2^2=L, x_1=2)=& \frac{p_2^1 p_2^2 (2 p_1^2 (1-p_1^2)) (\frac{1}{2}+p_1^1-(p_1^1)^2)}{(1-\frac{p_2^1}{2}+2 p_1^1 (1-p_1^1)(p_2^1-\frac{1}{2}))(1-p_2^2 + 2 p_1^2 (1-p_1^2)(2 p_2^2-1 ))}\\
    P(V^1=0,V^2=1| s_2^1=H, s_2^2=L, x_1=2)=& \frac{(1-p_2^1) (1-p_2^2) ((p_1^2)^2 + (1-p_1^2)^2) (1-p_1^1 +(p_1^1)^2)}{(1-\frac{p_2^1}{2}+2 p_1^1 (1-p_1^1)(p_2^1-\frac{1}{2}))(1-p_2^2 + 2 p_1^2 (1-p_1^2)(2 p_2^2-1 ))}\\
    P(V^1=0,V^2=0| s_2^1=H, s_2^2=L, x_1=2)=& \frac{(1-p_2^1) p_2^2 (2 p_1^2 (1-p_1^2)) (1-p_1^1 +(p_1^1)^2)}{(1-\frac{p_2^1}{2}+2 p_1^1 (1-p_1^1)(p_2^1-\frac{1}{2}))(1-p_2^2 + 2 p_1^2 (1-p_1^2)(2 p_2^2-1 ))}
\end{align}}

Suppose that $s_2^1=H, s_2^2=L, x_1=0$.
{\scriptsize
\begin{align}
    P(V^1=1,V^2=1| s_2^1=H, s_2^2=L, x_1=0)=& \frac{p_2^1 (1-p_2^2) (2 p_1^1 (1-p_1^1)) (2 p_1^2 (1-p_1^2))}{(1-p_2^1 + 2 p_1^1 (1-p_1^1) (2 p_2^1-1))(p_2^2 + 2 p_1^2 (1-p_1^2) (1-2 p_2^2))}\\
    P(V^1=1,V^2=0| s_2^1=H, s_2^2=L, x_1=0)=& \frac{p_2^1 p_2^2 (2 p_1^1 (1-p_1^1)) ((p_1^2)^2 + (1-p_1^2)^2)}{(1-p_2^1 + 2 p_1^1 (1-p_1^1) (2 p_2^1-1))(p_2^2 + 2 p_1^2 (1-p_1^2) (1-2 p_2^2))}\\
    P(V^1=0,V^2=1| s_2^1=H, s_2^2=L, x_1=0)=& \frac{(1-p_2^1) (1-p_2^2) (2 p_1^2 (1-p_1^2)) ((p_1^1)^2 + (1-p_1^1)^2)}{(1-p_2^1 + 2 p_1^1 (1-p_1^1) (2 p_2^1-1))(p_2^2 + 2 p_1^2 (1-p_1^2) (1-2 p_2^2))}\\
    P(V^1=0,V^2=0| s_2^1=H, s_2^2=L, x_1=0)=& \frac{(1-p_2^1) p_2^2 ((p_1^1)^2 + (1-p_1^1)^2)) ((p_1^2)^2 + (1-p_1^2)^2)}{(1-p_2^1 + 2 p_1^1 (1-p_1^1) (2 p_2^1-1))(p_2^2 + 2 p_1^2 (1-p_1^2) (1-2 p_2^2))}
\end{align}}

Suppose that $s_2^1=L, s_2^2=H, x_1=1$.
{\scriptsize
\begin{align}
    P(V^1=1,V^2=1| s_2^1=L, s_2^2=H, x_1=1)=& \frac{(1-p_2^1) p_2^2 ((p_1^1)^2 + (1-p_1^1)^2) (\frac{1}{2}+p_1^2-(p_1^2)^2)}{(1-p_2^1 + 2 p_1^1 (1-p_1^1) (2 p_2^1 -1))(1- \frac{p_2^2}{2}+ 2 p_1^2 (1-p_1^2) (p_2^2 -\frac{1}{2}))}\\
    P(V^1=1,V^2=0| s_2^1=L, s_2^2=H, x_1=1)=& \frac{(1-p_2^1) (1-p_2^2) ((p_1^1)^2 + (1-p_1^1)^2) (1-p_1^2 +(p_1^2)^2)}{(1-p_2^1 + 2 p_1^1 (1-p_1^1) (2 p_2^1 -1))(1- \frac{p_2^2}{2}+ 2 p_1^2 (1-p_1^2) (p_2^2 -\frac{1}{2}))}\\
    P(V^1=0,V^2=1| s_2^1=L, s_2^2=H, x_1=1)=& \frac{p_2^1 p_2^2 (2 p_1^1 (1-p_1^1)) (\frac{1}{2}+p_1^2 -(p_1^2)^2)}{(1-p_2^1 + 2 p_1^1 (1-p_1^1) (2 p_2^1 -1))(1- \frac{p_2^2}{2}+ 2 p_1^2 (1-p_1^2) (p_2^2 -\frac{1}{2}))}\\
    P(V^1=0,V^2=0| s_2^1=L, s_2^2=H, x_1=1)=& \frac{p_2^1 (1-p_2^2) (2 p_1^1 (1-p_1^1)) (1-p_1^2 +(p_1^2)^2)}{(1-p_2^1 + 2 p_1^1 (1-p_1^1) (2 p_2^1 -1))(1- \frac{p_2^2}{2}+ 2 p_1^2 (1-p_1^2) (p_2^2 -\frac{1}{2}))}
\end{align}}

Suppose that $s_2^1=L, s_2^2=H, x_1=2$.
{\scriptsize
\begin{align}
    P(V^1=1,V^2=1| s_2^1=L, s_2^2=H, x_1=2)=& \frac{(1-p_2^1) p_2^2 ((p_1^2)^2 + (1-p_1^2)^2) (\frac{1}{2}+p_1^1-(p_1^1)^2)}{(\frac{1+ p_2^1}{2} + 2 p_1^1 (1-p_1^1)(\frac{1}{2}-p_2^1))(p_2^2 + 2 p_1^2 (1-p_1^2)(1-2 p_2^2))}\\
    P(V^1=1,V^2=0| s_2^1=L, s_2^2=H, x_1=2)=& \frac{(1-p_2^1) (1-p_2^2) (2 p_1^2 (1-p_1^2)) (\frac{1}{2}+p_1^1-(p_1^1)^2)}{(\frac{1+ p_2^1}{2} + 2 p_1^1 (1-p_1^1)(\frac{1}{2}-p_2^1))(p_2^2 + 2 p_1^2 (1-p_1^2)(1-2 p_2^2))}\\
    P(V^1=0,V^2=1| s_2^1=L, s_2^2=H, x_1=2)=& \frac{p_2^1 p_2^2 ((p_1^2)^2 + (1-p_1^2)^2) (1-p_1^1 +(p_1^1)^2)}{(\frac{1+ p_2^1}{2} + 2 p_1^1 (1-p_1^1)(\frac{1}{2}-p_2^1))(p_2^2 + 2 p_1^2 (1-p_1^2)(1-2 p_2^2))}\\
    P(V^1=0,V^2=0| s_2^1=L, s_2^2=H, x_1=2)=& \frac{p_2^1 (1-p_2^2) (2 p_1^2 (1-p_1^2)) (1-p_1^1 +(p_1^1)^2)}{(\frac{1+ p_2^1}{2} + 2 p_1^1 (1-p_1^1)(\frac{1}{2}-p_2^1))(p_2^2 + 2 p_1^2 (1-p_1^2)(1-2 p_2^2))}
\end{align}}

Suppose that $s_2^1=L, s_2^2=H, x_1=0$.
{\scriptsize
\begin{align}
    P(V^1=1,V^2=1| s_2^1=L, s_2^2=H, x_1=0)=& \frac{(1-p_2^1) p_2^2 (2 p_1^1 (1-p_1^1)) (2 p_1^2 (1-p_1^2))}{(p_2^1 + 2 p_1^1 (1-p_1^1)(1-2 p_2^1))(1-p_2^2 + 2 p_1^2 (1-p_1^2)(2 p_2^2-1))}\\
    P(V^1=1,V^2=0| s_2^1=L, s_2^2=H, x_1=0)=& \frac{(1-p_2^1) (1-p_2^2) (2 p_1^1 (1-p_1^1)) ((p_1^2)^2 + (1-p_1^2)^2)}{(p_2^1 + 2 p_1^1 (1-p_1^1) (1-2 p_2^1))(1-p_2^2 + 2 p_1^2 (1-p_1^2)(2 p_2^2-1))}\\
    P(V^1=0,V^2=1| s_2^1=L, s_2^2=H, x_1=0)=& \frac{p_2^1 p_2^2 (2 p_1^2 (1-p_1^2)) ((p_1^1)^2 + (1-p_1^1)^2)}{(p_2^1 + 2 p_1^1 (1-p_1^1) (1-2 p_2^1))(1-p_2^2 + 2 p_1^2 (1-p_1^2)(2 p_2^2-1))}\\
    P(V^1=0,V^2=0| s_2^1=L, s_2^2=H, x_1=0)=& \frac{p_2^1 (1-p_2^2) ((p_1^1)^2 + (1-p_1^1)^2)) ((p_1^2)^2 + (1-p_1^2)^2)}{(p_2^1 + 2 p_1^1 (1-p_1^1) (1-2 p_2^1))(1-p_2^2 + 2 p_1^2 (1-p_1^2)(2 p_2^2-1))}
\end{align}
}
Suppose that $s_2^1=L, s_2^2=L, x_1=1$.
{\scriptsize
\begin{align}
    P(V^1=1,V^2=1| s_2^1=L, s_2^2=L, x_1=1)=& \frac{(1-p_2^1) (1-p_2^2) ((p_1^1)^2 + (1-p_1^1)^2) (\frac{1}{2}+p_1^2-(p_1^2)^2)}{(1-p_2^1 + 2 p_1^1 (1-p_1^1)(2 p_2^1-1))(\frac{1+p_2^2}{2}+ 2 p_1^2 (1-p_1^2)(\frac{1}{2}-p_2^2))}\\
    P(V^1=1,V^2=0| s_2^1=L, s_2^2=L, x_1=1)=& \frac{(1-p_2^1) p_2^2 ((p_1^1)^2 + (1-p_1^1)^2) (1-p_1^2 +(p_1^2)^2)}{(1-p_2^1 + 2 p_1^1 (1-p_1^1)(2 p_2^1-1))(\frac{1+p_2^2}{2}+ 2 p_1^2 (1-p_1^2)(\frac{1}{2}-p_2^2))}\\
    P(V^1=0,V^2=1| s_2^1=L, s_2^2=L, x_1=1)=& \frac{p_2^1 (1-p_2^2) (2 p_1^1 (1-p_1^1)) (\frac{1}{2}+p_1^2 -(p_1^2)^2)}{(1-p_2^1 + 2 p_1^1 (1-p_1^1)(2 p_2^1-1))(\frac{1+p_2^2}{2}+ 2 p_1^2 (1-p_1^2)(\frac{1}{2}-p_2^2))}\\
    P(V^1=0,V^2=0| s_2^1=L, s_2^2=L, x_1=1)=& \frac{p_2^1 p_2^2 (2 p_1^1 (1-p_1^1)) (1-p_1^2 +(p_1^2)^2)}{(1-p_2^1 + 2 p_1^1 (1-p_1^1)(2 p_2^1-1))(\frac{1+p_2^2}{2}+ 2 p_1^2 (1-p_1^2)(\frac{1}{2}-p_2^2))}
\end{align}}

Suppose that $s_2^1=L, s_2^2=L, x_1=2$.
{\scriptsize\begin{align}
    P(V^1=1,V^2=1| s_2^1=L, s_2^2=L, x_1=2)=& \frac{(1-p_2^1) (1-p_2^2) ((p_1^2)^2 + (1-p_1^2)^2) (\frac{1}{2}+p_1^1-(p_1^1)^2)}{(\frac{1+p_2^1}{2} + 2 p_1^1 (1-p_1^1) (\frac{1}{2}-p_2^1))(1-p_2^2 + 2 p_1^2 (1-p_1^2) (2 p_2^2-1))}\\
    P(V^1=1,V^2=0| s_2^1=L, s_2^2=L, x_1=2)=& \frac{(1-p_2^1) p_2^2 (2 p_1^2 (1-p_1^2)) (\frac{1}{2}+p_1^1-(p_1^1)^2)}{(\frac{1+p_2^1}{2} + 2 p_1^1 (1-p_1^1) (\frac{1}{2}-p_2^1))(1-p_2^2 + 2 p_1^2 (1-p_1^2) (2 p_2^2-1))}\\
    P(V^1=0,V^2=1| s_2^1=L, s_2^2=L, x_1=2)=& \frac{p_2^1 (1-p_2^2) ((p_1^2)^2 + (1-p_1^2)^2) (1-p_1^1 +(p_1^1)^2)}{(\frac{1+p_2^1}{2} + 2 p_1^1 (1-p_1^1) (\frac{1}{2}-p_2^1))(1-p_2^2 + 2 p_1^2 (1-p_1^2) (2 p_2^2-1))}\\
    P(V^1=0,V^2=0| s_2^1=L, s_2^2=L, x_1=2)=& \frac{p_2^1 p_2^2 (2 p_1^2 (1-p_1^2)) (1-p_1^1 +(p_1^1)^2)}{(\frac{1+p_2^1}{2} + 2 p_1^1 (1-p_1^1) (\frac{1}{2}-p_2^1))(1-p_2^2 + 2 p_1^2 (1-p_1^2) (2 p_2^2-1))}
\end{align}}

Suppose that $s_2^1=L, s_2^2=L, x_1=0$.
{\scriptsize\begin{align}
    P(V^1=1,V^2=1| s_2^1=L, s_2^2=L, x_1=0)=& \frac{(1-p_2^1) (1-p_2^2) (2 p_1^1 (1-p_1^1)) (2 p_1^2 (1-p_1^2))}{(p_2^1 + 2 p_1^1 (1-p_1^1)(1-2 p_2^1))(p_2^2 + 2 p_1^2 (1-p_1^2)(1-2 p_2^2))}\\
    P(V^1=1,V^2=0| s_2^1=L, s_2^2=L, x_1=0)=& \frac{(1-p_2^1) p_2^2 (2 p_1^1 (1-p_1^1)) ((p_1^2)^2 + (1-p_1^2)^2)}{(p_2^1 + 2 p_1^1 (1-p_1^1)(1-2 p_2^1))(p_2^2 + 2 p_1^2 (1-p_1^2)(1-2 p_2^2))}\\
    P(V^1=0,V^2=1| s_2^1=L, s_2^2=L, x_1=0)=& \frac{p_2^1 (1-p_2^2) (2 p_1^2 (1-p_1^2)) ((p_1^1)^2 + (1-p_1^1)^2)}{(p_2^1 + 2 p_1^1 (1-p_1^1)(1-2 p_2^1))(p_2^2 + 2 p_1^2 (1-p_1^2)(1-2 p_2^2))}\\
    P(V^1=0,V^2=0| s_2^1=L, s_2^2=L, x_1=0)=& \frac{p_2^1 p_2^2 ((p_1^1)^2 + (1-p_1^1)^2)) ((p_1^2)^2 + (1-p_1^2)^2)}{(p_2^1 + 2 p_1^1 (1-p_1^1)(1-2 p_2^1))(p_2^2 + 2 p_1^2 (1-p_1^2)(1-2 p_2^2))}
\end{align}}

The decision probabilities for the second backer under OL can be found by

{\scriptsize\begin{equation}
    P(x_1, x_2| V^1, V^2)=P(x_1| V^1, V^2) \sum_{s_2^1, s_2^2} P(x_2| s_2^1, s_2^2, x_1) P(s_2^1|V^1) P(s_2^2|V^2),
\end{equation}}

for all possible true quality states.  

In the next step, we use these probabilities to calculate the sample path probabilities. Given that each backer can choose between three pledging decisions, for the 2-backers system we have nine possibilities. Also, considering that there are four possible true quality states, in total we evaluate 36 sample path probabilities of the following form. 

{\scriptsize\begin{equation}
    P\{x_1,x_2| V^1, V^2 \}.
\end{equation}}

The expressions for these probabilities are hard to simplify. For this reason, we omit them in this document. We use python scripts to formulate these in functions of the posterior probabilities, which can be provided upon request from the authors.
The last step is to evaluate the performance measures that are defined over the set of all expertness scenarios for the two backers using sample path probabilities. Through numerical evaluations on the region $[0.5,1] \times [0.5, 1]$ for backer expertness, we derive the crowdfunding performance measures.